\documentclass[11pt,amssymb,amsfonts]{article}
\usepackage{latexsym}
\usepackage{amssymb}
\usepackage{amscd}
\usepackage{amsfonts}
\usepackage{stmaryrd}
\usepackage{eufrak}

\usepackage {pstricks,pst-node}
\usepackage [dvips] {graphics, graphicx}

\def\on{\hskip2truemm}
\def\frf#1{\vcenter{\hrule\hbox{\vrule\kern1pt
\vbox{\kern1pt\hbox{$\displaystyle#1$}%
\kern1pt}\kern1pt\vrule}\hrule}}

\newcommand{\be}{\begin{equation}}
\newcommand{\ee}{\end{equation}}
\def\bea{\begin{eqnarray}}
\def\eea{\end{eqnarray}}
\def\beq{\begin{eqnarray*}}
\def\eeq{\end{eqnarray*}}
\def\ba{\begin{array}}
\def\ea{\end{array}}

\def\ss{\scriptstyle}

\newtheorem{theorem}{\sc Theorem}[section]
\newtheorem{lemma}{\sc Lemma}[section]
\newtheorem{definition}{\sc Definition}[section]
\newtheorem{example}{\sc Example}[section]

\newtheorem{corollary}{\sc Corollary}[section]
\newtheorem{proposition}{\sc Proposition}[section]

\def\notni{\hbox{{\kern0pt\raise.7pt\hbox{${\scriptstyle +}$}}
{\kern-9.5pt\raise0pt\hbox{$\supset$}}}}

\def\Cb{\mathbb{C}}

\def\Zb{{\mathbb Z}}

\def\gtL{{\mathfrak L}}

\def\bj{{\bf j}}

\def\bcup{\mathop{\bigcup}}

\def\Aut{{\rm Aut}}
\def\End{{\rm End}}

\def\Mor{{\rm Mor}}

\def\Ob{{\rm Ob}}

\def\im{{\rm im}}
\def\Tr{{\rm Tr}}

\def\Ann{{\rm Ann}}
\def\nully{{\rm null}}
\def\Indef{{\rm Indef}}

\def\min{{\rm min}}
\def\max{{\rm max}}

\def\spin{{\rm spin}}

\newcommand{\biota}{\mbox{\boldmath $\iota$}}

\begin{document}

\setcounter{section}{0}
\setcounter{equation}{0}

\begin{center}
{\Large\bf Representations of the Complex Classical Cayley-Klein Categories}

\bigskip

Stepan S. Moskaliuk

\medskip
{\it Bogolyubov Institute for Theoretical Physics\\
Metrolohichna Str., 14-b, Kyiv-143, Ukraine, UA-03143\\
e-mail: mss@bitp.kiev.ua}

\end{center}

\bigskip
\centerline{Abstract}

\medskip

Complex classical Cayley-Klein categories ${\bf A({\bj})}$, ${\bf B({\bj})}$, 
${\bf C({\bj})}$ and  ${\bf D({\bj})}$ are constructed by
the method of categorical extension of the complex classical
Cayley-Klein groups $SL(2n;{\bj};{\Bbb C}),\; SO(2n+1;{\bj};{\Bbb C}),\; 
Sp(2n;{\bj};{\Bbb C})$ and $SO(2n;{\bj};{\Bbb C})$, respectively. The
explicit
construction of the irreducible representations of the complex classical
Cayley-Klein categories is received. Completeness of lists of
 the irreducible representations of the complex classical Cayley-Klein
categories and the classification theorems are proved.

\newpage

\setcounter{equation}{0}
\begin{center}
\noindent

\section{Pimenov algebra and Complex Cayley-Klein spaces}
\end{center}
\noindent

Let us define {\it Pimenov algebra} $ {\cal D}_n( ${\mbox{\boldmath$ \iota$}}$; {\Bbb C}) $
as an associative algebra with unit over complex number field
and with nilpotent commutative generators
$ {\iota}_k, \ {\iota}_k^2=0,$
$ {\iota}_k{\iota}_m={\iota}_m{\iota}_k \not =0, \  k \neq m,
\  k,m=1, \ldots, n. $
The general element of ${\cal D}_n(${\mbox{\boldmath$ \iota$}}$; {\Bbb C}) $
is in the form
$$
d=d_0+\sum^{n}_{p=1}\sum_{k_1< \ldots < k_p}d_{k_1\ldots k_p}
{\iota}_{k_1} \ldots {\iota}_{k_p}, \quad  d_0,d_{k_1 \ldots k_p} \in
{\Bbb C}.
$$
For  $ n=1 $  we have  $ {\cal D}_1({\iota}_1;{\Bbb C})  \ni
 d=d_0+d_1{\iota}_1, $ i.e. the elements $d$ are dual 
numbers when $ d_0,d_1 \in {\Bbb R} $ \cite{Cliff, Rosenf}.
For $ n=2 $ the general element of  $ {\cal D}_2(\iota_1, \iota_2; \Bbb C)$ is
$ d=d_0+d_1{\iota}_1+d_2{\iota}_2+d_{12}{\iota}_1{\iota}_2. $
Two elements $ d, \tilde d \in {\cal D}_n(${\mbox{\boldmath$ \iota$}}$;{\Bbb C})$ are equal if and only if
$ d_0=\tilde d_0, \ d_{k_1 \ldots k_p}=\tilde d_{k_1\ldots k_p},\ p=1,
\ldots,n. $
If $ d=d_k{\iota}_k $ and  $ \tilde d=\tilde d_k{\iota}_k, $ then
the condition $ d=\tilde d, $ which is equivalent to
$ d_k{\iota}_k=\tilde d_k{\iota}_k, $
make possible the consistently definition of the division of nilpotent
generator  $ {\iota}_k $ by itself, namely:
$ {\iota}_k/{\iota}_k=1, \ k=1, \ldots, n. $
Let us stress that the division of different nilpotent generators
$ {\iota}_k/{\iota}_p, \ k \not = p, $ as well as the division of
complex number by nilpotent generators
$ a/{\iota}_k, \ a \in {\Bbb C} $  are not defined.
It is convenient to regard the algebras
$ {\cal D}_n(\bj; {\Bbb C}), $ where the parameters
$ j_k=1,\iota_k, \; k=1,\ldots,n. $
If $ m $ parameters are nilpotent $ j_{k_{s}}=\iota_{s}, \; s=1,\ldots,m $
and the other are equal to unit, then we have  Pimenov algebra
$ {\cal D}_m(${\mbox{\boldmath$ \iota$}}$; {\Bbb C}).$


 {\it Cayley-Klein spaces of ``complex''
type} $\Bbb C_{n+1}({\bj})$  come out from $(n+1)$-dimensional complex space $\Bbb 
C_{n+1}$ under the mapping 
\begin{eqnarray}
\label{1.91} 
\psi&:&\Bbb C_{n+1}\to \Bbb C_{n+1}({\bj}), \nonumber\\
\psi z_0^*=z_0,\quad \psi z_k^*&=&z_k\prod_{m=1}^k j_m,\quad
k=1,2,\dots,n,
\end{eqnarray}
  where
$z^*_0,z^*_k\in\Bbb C_{n+1},\ z_0,z_k\in\Bbb C_{n+1}(\bj)$
 are complex Cartesian coordinates,
$\bj=(j_1,\dots,j_n)$, each of
parameters $j_k$ takes three values:
$j_k=1,\iota_k,i$.

``Complex'' type means here that these spaces are coordinatized by
elements of a one-step extension of ${\Bbb R}$ through a labelled
Cayley-Dickson procedure ${\Bbb R} \Rightarrow {\Bbb R}(\mu_1)$ which adjoint 
a imaginary unit $i$ with $i^2 = - 1$ to ${\Bbb R}$, producing either the
{\it complex numbers} $\Bbb C$ ($\mu_1=1$), {\it double 
(split complex) numbers} ($\mu_1=-1$) and {\it dual numbers} 
($\mu_1=0$) \cite{Lohmus}. 

Within this family, only the spaces coordinatised by ordinary compex
numbers (where $i^2 = - 1$) are actually complex spaces; after this
restriction, there are only $3^N$ Cayley-Klein complex spaces in 
dimension N.

Quadratic form
$({\mbox{\boldmath$ z $}}^*, {\mbox{\boldmath$ 
z $}}^*)=\sum\limits_{m=0}^{n}|z_m^*|^2$ of the space $\Bbb C_{n+1}$ 
 turns into quadratic form
\begin{equation}
\label{1.92} 
({\mbox{\boldmath$ z $}}, {\mbox{\boldmath $z$}})=
|z_0|^2=\sum_{k=1}^{n}|z_k|^2 \prod_{m=1}^{k}j_m^2.  
\end{equation}
of the space
$\Bbb C_{n+1}( {\mbox{\boldmath$j $}})$
 under the mapping (\ref{1.91}). Here
$|z_k|=(x_k^2+y_k^2)^{1/2}$
 is absolute value (modulus)
of complex number
$z_k=x_k+iy_k$, and ${\mbox{\boldmath$ z $}}$ is complex vector:
${\mbox{\boldmath$ z $}}=(z_0,z_1,\dots,z_n)$.

A space
$\Bbb C_{n+1}({\mbox{\boldmath$ j $}})$
 is called non-fiber space, if no one of the parameters
$j_1,\dots,j_n$
 take dual
value. A space
$\Bbb C_{n+1}({\mbox{\boldmath$ j $}})$
 is called
$(k_1,k_2\dots,k_p)$
  -- fiber space, if
$1\leq k_1<k_2<\dots<k_p\leq n$
 and
$j_{k_1}=\iota_{k_1},\ \dots,\ j_{k_p}=\iota_{k_p}$, and the other
$j_k=1$.
These fiberings are trivial \cite{Bourb1} and can be characterised by the set of
consequently nested projections
$pr_1,pr_2,\dots,pr_p$; where for $pr_1$ the base is a subspace,
spanned over the basis vectors
$\{e_0,e_1,\dots,e_{k_1-1}\}$, and the fiber is a subspace,
spanned over
$\{e_{k_1},e_{k_1+1},\dots,e_n\}$;
for $pr_2$ the~base~is~a~subspace
$\{e_{k_1},e_{k_1+1},\dots,e_{k_2-1}\}$,
and the fiber is a subspace
$\{e_{k_2},e_{k_2+1},\dots,e_n\}$
 and so on.

\medskip

Let us study space $\Bbb C_n({\mbox{\boldmath$ j$}})$,
obtained from $n$-dimensional Euclidean space $\Bbb C_n$ using the
mapping (\ref{1.91}).

If all $j_k=1$, then $\psi$ is identical mapping.  For
$j_{k_1}=\iota_{k_1},\dots,j_{k_p}=\iota_{k_p}$ and other $j_k=1$ we
obtain fiber space with zero curvature $k_1,k_2,\dots,k_p$. This space
is characterised by the existence consequently nested projections
$pr_1,pr_2,\dots,pr_p$, where for $pr_1$ subspace
$\{z_1,z_2,\dots,z_{k_1-1}\}$ serves as a base, and space
$\{z_{k_1},z_{k_1+1},\dots,z_n\}$ -- as a
fiber; for $pr_2$ as a base subspace $\{z_{k_1},z_{k_1+1},\dots,
x_{k_2-1}\}$ and as a fiber - subspace $\{z_{k_2},z_{k_2+1},\dots,
x_n\}$ and so on.

From mathematical point of view, fibering in space $\Bbb
C_n({\mbox{\boldmath$ j $}})$ is trivial \cite{Bourb1}, i.e.  globally it has 
 the same structure as locally.  From physical point of view, fibering 
originates absolute physical quantities.  For example, complex Galilean space, 
which is realized on the sphere $S_4(1,\iota_2,1,1)$, can be 
characterised by the existence of absolute  time $t=z_1$ and absolute 
 space $\Bbb C_3= \{z_2,z_3,z_4\}$.

\medskip
\setcounter{equation}{0}
\section{Complex classical Cayley-Klein groups}
\subsection{Unitary Cayley-Klein groups}


\begin{definition}
\label{d1.6} 
 Group
$SU(n+1;{\mbox{\boldmath$ j $}}; \Bbb C)$
 consists of all transformations of
space
$\Bbb C_{n+1}( {\mbox{\boldmath$j $}})$
 with unit determinant, keeping invariant the quadratic form
(\ref{1.92}).
\end{definition}

In the
$(k_1,k_2,\dots,k_p)$
 -- fiber space
$\Bbb C_{n+1}({\mbox{\boldmath$ j $}})$
 we have $(p+1)$ -- quadratic form, which remains
invariant under the transformations of group
$SU(n+$\linebreak$+1;{\mbox{\boldmath$ j $}};\Bbb C)$. Under transformations of
group
$SU(n+1; {\mbox{\boldmath$j $}};\Bbb C)$, which do not affect coordinates
$z_0,z_1,\dots,z_{k_s-1}$, the form
\begin{equation}
\label{1.93} 
({\mbox{\boldmath$ z$}},
{\mbox{\boldmath$ z$}})_{s+1}=\sum_{\alpha=k_s}^{k_{s+1}-1}|z_\alpha|^2
\prod^{\alpha}_{l=k_s+1}j_l^2,
\end{equation}
where
$s+0,1,\dots,p,\ k_0=0$, remains invariant.

For $s=p$ the summation over $\alpha$ goes up to $n$.

The mapping (\ref{1.91}) induces the transition of classical group $SU(n+1;\Bbb C)$
into group
$SU(n+1;{\mbox{\boldmath$ j$}};\Bbb C)$. All $(n+1)^2-1$ generators of group
$SU(n+1;\Bbb C)$ are
Hermitean matrices. If group
$GL_{n+1}$ acts via left translations in the space of analytical
functions on $\Bbb C_{n+1}$, then its generators are
$X_{\alpha\beta}^*=z^{*\beta}\partial_\alpha{}^*$, where
$\partial_\alpha{}^*\equiv\partial/\partial z^{*\alpha}$. 
(The asterisk means that $A^*$ is a generator of a
classical group.)

 Commutators of generators $X^*$ satisfy the
following relation
\begin{equation}
\label{1.94} 
[X_{km}^*,X_{pq}^*]=\delta_{mp}X_{kq}^*-\delta_{kq}X_{pm}^*,
\end{equation}
where $\delta_{mp}$ is Kronecker symbol.

Independent Hermitean generators of group $SU(n+1;\Bbb C)$ are given by
equations
\begin{eqnarray}\nonumber
\label{1.95} 
Q_{\mu\nu}^*={i\over 2}(X_{\mu\nu}^*+X_{\nu\mu}^*),\quad
L_{\mu\nu}^*={1\over 2}(X_{\nu\mu}^*-X_{\mu\nu}^*),\\
P_k^*={i\over 2}(X_{k-1,k-1}^*-X_{kk}^*),
\end{eqnarray}
where
$\mu=0,1,\dots,n-1$; $\nu=\mu+1,\mu+2,\dots,n$; $k=1,2,\dots,n$.

Generators $X^*_{\alpha\beta}$ are transformed under the mapping (\ref{1.91}) as
follows:
\begin{equation}
\label{1.96} 
X_{kk}=z_k\partial_k,\quad
X_{\nu\mu}=z_\mu\partial_\nu,\quad
X_{\mu\nu}=\biggl(\prod^{\nu}_{l=\mu+1}j_l^2\biggr)z_\nu\partial_\mu,
\end{equation}
where
$k=1,2,\dots,n$; $\mu,\nu=0,1,\dots,n$, $\mu<\nu$.

Commutators of generators $X({\mbox{\boldmath$ j $}})$ can be easily 
found \cite{GMUniG}:  
\begin{equation}
\label{1.97} 
[X_{km},X_{pq}]=\prod^{l_2}_{l=l_1}j_l\prod^{l_4}_{l=l_3}j_l
\biggl(\delta_{mp}X_{kq}\prod^{l_6}_{l=l_5}j_l^{-1}-
\delta_{kq}X_{pm}\prod^{l_8}_{l=l_7}j_l^{-1}\biggr),
\end{equation}
where
$l_1=1+\min(k,m)$; $l_2=\max(k,m)$;
$l_3=1+\min(p,q)$; $l_4=$ \linebreak $=\max(p,q)$;
$l_5=1+\min(k,q)$; $l_6=\max(k,q)$;
$l_7=1+\min(m,p)$; $l_8=\max(m,p)$.

Hermitian generators (\ref{1.95}) are transformed in the same way under 
transition from group $SU(n+1;\Bbb C)$ to group $SU(n+1;{\mbox{\boldmath$ j 
$}};\Bbb C)$. This enables to find matrix generators of group $SU(n+1; 
{\mbox{\boldmath$ j $ }};\Bbb C)$:   
\begin{eqnarray}\nonumber 
\label{1.98} 
Q_{\mu\nu}({\mbox{\boldmath$ j 
$}})=\biggl(\prod^{\nu}_{m=\mu+1}j_m\biggr) Q_{\mu\nu}^*(\to)={i\over 
2}[X_{\nu\mu}({\mbox{\boldmath$ j$}})+X_{\mu\nu}({\mbox{\boldmath$ 
j$}})]= \nonumber\\ 
={i\over 2}\biggl(X_{\nu\mu}^*+X_{\mu\nu}^* 
\prod^{\nu}_{m=\mu+1}j_m^2\biggr),\nonumber\\
L_{\mu\nu}({\mbox{\boldmath$ j$}})=\biggl(\prod^{\nu}_{m=\mu+1}j_m\biggr)
L_{\mu\nu}^*(\to)={1\over 2}[X_{\nu\mu}({\mbox{\boldmath$
j$}})-X_{\mu\nu}({\mbox{\boldmath$ j$}})]=  \\
={1\over 2}\biggl(X_{\nu\mu}^*-X_{\mu\nu}^*
\prod^{\nu}_{m=\mu+1}j_m^2\biggr), \nonumber\\
P_k({\mbox{\boldmath$ j$}})=P_k^*={i\over 2}(X_{k-1\,k-1}^*-X_{kk}^*),
~~~~~~~k=1,2,\dots,n.\nonumber
\end{eqnarray}
We do not cite commutation relations for generators
$Q_{\mu\nu}({\mbox{\boldmath$ j$ }})$, $L_{\mu\nu}({\mbox{\boldmath$ 
j$}})$, $P_k({\mbox{\boldmath$ j $}})$ because they are cumbersome 
\cite{GMUniG}.  They can be found, using (\ref{1.97}).

Matrix generators (\ref{1.98}) make a basis of Lie algebra
$su\,(n+1;{\mbox{\boldmath$ j $}};\Bbb C)$. To the general element of the 
algebra 
\begin{equation} 
\label{1.101} 
Z({\mbox{\boldmath$ r$}},{\mbox{\boldmath$ 
s$}},{\mbox{\boldmath$ w$}}, {\mbox{\boldmath$ j$}})= 
\sum_{\lambda=1}^{n(n+1)/2}(r_\lambda Q_\lambda({\mbox{\boldmath$ j$}})+
s_\lambda L_\lambda({\mbox{\boldmath$ j$}}))+\sum_{k=1}^{n}w_kP_k,
\end{equation}
where index $\lambda$ is connected with the indices $\mu$, $\nu$,
$\mu<\nu$ by relation
\begin{equation}
\label{1.102} 
\lambda=\nu+\mu(n-1)-\mu(\mu-1)/2,
\end{equation}
and the group parameters
$r_\lambda$, $s_\lambda$, $w_k$
are complex, there
corresponds a finite group transformation of group 
$SU(n+1;{\mbox{\boldmath$j $}};\Bbb C)$
\begin{equation}
\label{1.103} 
\Xi({\mbox{\boldmath$ r$}}, {\mbox{\boldmath$s$}},{\mbox{\boldmath$ w$}},
{\mbox{\boldmath$j$}})=
\exp Z({\mbox{\boldmath$ r$}},{\mbox{\boldmath$ s$}},{\mbox{\boldmath$ w$}},  {\mbox{\boldmath$j$}})
\end{equation}

\subsection{Complex orthogonal Cayley-Klein group}

{\it Complex orthogonal Cayley-Klein group} $ SO(N;\bj;\Bbb C) $
is defined as the group of transformations
$ {\xi}'(\bj)=A(\bj)\xi(\bj)$ of complex vector space $O_N(\bj)$
with Cartesian coordinates
$ \xi^t (\bj) = (\xi_1,  (1,2)  \xi_2,  \ldots,(1,N)\xi_N)^t, \;$
which preserve the quadratic form
$$
inv(\bj) =\xi^t(\bj)\xi(\bj) = \xi^2_1 + \sum^N_{k=2}({1,k})^2{\xi}^2_k,
$$
where $ \bj=(j_1, \ldots, j_{N-1}) $, each parameter $ j_k $
takes {\it two} values:
$  j_r=1,{\iota}_r; \ r=1, \ldots, N-1, \ {\xi}_k \in \Bbb C$
and
$$
(\mu,\nu)=\prod^{max(\mu,\nu)-1}_{l=min(\mu,\nu)} j_l, \quad (\mu,\mu)=1.
$$
Let us stress, that Cartesian coordinates
of $  O_N(\bj) $ are special elements of Pimenov algebra
${\cal D}_{N-1}(\bj; \Bbb C). $
 It worth notice that the  orthogonal Cayley-Klein groups as well as 
the unitary and symplectic   Cayley-Klein   groups   have   been   
regarded  in \cite{GMOrtG,GMSymG} as the matrix groups with the real matrix 
elements.  Nevertheless there is a different approach, which gives the 
same results for ordinary groups \cite{GKKqortG,GKKqSymG}.   According with this approach,  
the Cayley-Klein group $ SO(N;\bj;{\Bbb C}) $ may be realised as the matrix 
group, whose elements are taken from algebra $ {\cal D}_{N-1}(\bj;{\Bbb C}) $ 
and in Cartesian basis consist of the $ N \times N $ matrices $A(\bj)$ 
with elements 
$$ (A(\bj))_{kp}=(k,p)a_{kp}, \ a_{kp} \in \Bbb C.  
$$ 
Matrices $A(\bj)$ are subject of the additional $\bj$-orthogonality relations 
\begin{equation}
A(\bj)A^t(\bj) = A^t(\bj)A(\bj)=I.
\label{3}
\end{equation}

Sometimes it is convenient to regard an orthogonal group in so-called
``symplectic'' basis. Transformation from Cartesian to symplectic basis
$ x(\bj)=D\xi (\bj) $ is made by unitary matrices $D$,
which are a solutions of equation
     \begin{equation}
 D^{t}C_{0}D=I,
     \label{1}
     \end{equation}
where $C_{0}\in M_{N},\; (C_{0})_{ik}= \delta_{ik'},\; k'=N+1-k$
To obtain all solutions of equation (\ref{1}),
 take one of them, namely
$$
D=\frac{1}{\sqrt{2}}
\left ( \begin{array}{cc}
      I &     -i{\tilde C_0} \\
      {\tilde C_0} &   iI
      \end{array} \right ),    \       N=2n,
$$
\begin{equation}
D=\frac{1}{\sqrt{2}}
\left ( \begin{array}{ccc}
      I & 0 &  -i{\tilde C_0} \\
      0 & \sqrt{2} &  0 \\
      {\tilde C_0} & 0  &  iI
      \end{array} \right ),    \       N=2n+1,
\label{2}
\end{equation}
where $n \times n $ matrix $ {\tilde C_0} $ is like $C_0,$
then  regard the matrix $ D_{\sigma}=DV_{\sigma},$
$V_{\sigma} \in M_{N},$ $(V_{\sigma})_{ik}= \delta_{\sigma_{i},k},$
and $ \sigma \in S(N) $ is a permutation of the $ N$th order.
It is easy to verify that $ D_{\sigma}$
is again a solution of equation (\ref{1}).
Then in symplectic basis the orthogonal Cayley Klein group
$SO(N;\bj;{\Bbb C})$ is described by the matrices
     \begin{equation}
 B_{\sigma}(\bj)=D_{\sigma}A(\bj)D^{-1}_{\sigma}
    \label{4}
    \end{equation}
with the additional relations of $ \bj$-orthogonality
$$
   B_{\sigma}(\bj)C_{0}B^{t}_{\sigma}(\bj)=
B^{t}_{\sigma}(\bj)C_{0}B_{\sigma}(\bj)=C_{0}.
$$

It should be noted that for  orthogonal groups $(\bj=1)$
the use of different matrices $D_{\sigma}$ makes no sense because  all
Cartesian coordinates of $O_N$ are equivalent up to a choice of
 its enumerations.
The different situation is for Cayley-Klein groups $(\bj \not =1).$
Cartesian coordinates
$(1,k)\xi_{k}, \ k=1,\dots,N$ for nilpotent values of some or all
parameters $j_k$ are different elements of the algebra $D_{N-1}(\bj; \Bbb C),$
therefore the same group $SO(N;\bj; \Bbb C)$ may be realized by  matrices
$ B_{\sigma} $  with a {\it different} disposition of nilpotent
generators among their elements.

Matrix elements of $ B_{\sigma}(\bj) $ are as follows
\begin{equation}
     \begin{array}{l}
(B_{\sigma})_{n+1,n+1}= b_{n+1,n+1},  \\
(B_{\sigma})_{kk}=b_{kk}+i\tilde{b}_{kk}(\sigma_{k},
\sigma_{k'}), \quad
(B_{\sigma})_{k'k'}=b_{kk}-i\tilde{b}_{kk}(\sigma_{k},
\sigma_{k'}), \\
(B_{\sigma})_{kk'}=b_{k'k}-i\tilde{b}_{k'k}(\sigma_{k},
\sigma_{k'}), \quad
(B_{\sigma})_{k'k}=b_{k'k}+i\tilde{b}_{k'k}(\sigma_{k},
\sigma_{k'}), \\
(B_{\sigma})_{k,n+1}=b_{k,n+1}(\sigma_k, \sigma_{n+1}) -
i\tilde{b}_{k,n+1}(\sigma_{n+1}, \sigma_{k'}), \\
(B_{\sigma})_{k',n+1}=b_{k,n+1}(\sigma_k, \sigma_{n+1}) +
i\tilde{b}_{k,n+1}(\sigma_{n+1}, \sigma_{k'}), \\
(B_{\sigma})_{n+1,k}=b_{n+1,k}(\sigma_k, \sigma_{n+1}) +
i\tilde{b}_{n+1,k}(\sigma_{n+1}, \sigma_{k'}), \\
(B_{\sigma})_{n+1,k'}=b_{n+1,k}(\sigma_k, \sigma_{n+1}) -
i\tilde{b}_{n+1,k}(\sigma_{n+1}, \sigma_{k'}), \; k \neq p, \\
(B_{\sigma})_{kp}=b_{kp}(\sigma_k, \sigma_{p})+
b_{kp}'(\sigma_{k'},\sigma_{p'})+
i\tilde{b}_{kp}(\sigma_{k}, \sigma_{p'})-
i\tilde{b}_{kp}'(\sigma_{k'}, \sigma_{p}),  \\
(B_{\sigma})_{kp'}=b_{kp}(\sigma_k, \sigma_{p})-
b_{kp}'(\sigma_{k'},\sigma_{p'})-
i\tilde{b}_{kp}(\sigma_{k}, \sigma_{p'})-
i\tilde{b}_{kp}'(\sigma_{k'}, \sigma_{p}),  \\
(B_{\sigma})_{k'p}=b_{kp}(\sigma_k, \sigma_{p})-
b_{kp}'(\sigma_{k'},\sigma_{p'})+
i\tilde{b}_{kp}(\sigma_{k}, \sigma_{p'})+
i\tilde{b}_{kp}'(\sigma_{k'}, \sigma_{p}),  \\
(B_{\sigma})_{k'p'}=b_{kp}(\sigma_k, \sigma_{p})+
b_{kp}'(\sigma_{k'},\sigma_{p'})-
i\tilde{b}_{kp}(\sigma_{k}, \sigma_{p'})+
i\tilde{b}_{kp}'(\sigma_{k'}, \sigma_{p}) . \\
\end{array}
   \label{5}
\end{equation}

Here $ b,  b', \tilde{b}, \tilde{b}' \in \Bbb C$
are expressed by the matrix elements of $A $ as
\begin{displaymath}
     \begin{array}{ll}
b_{n+1,n+1}=a_{\sigma_{n+1},\sigma_{n+1}}, &   \\
b_{n+1,k}=\displaystyle{\frac{1}{\sqrt{2}}a_{\sigma_{n+1},
\sigma_{k}}},
&
b_{k,n+1}=\displaystyle{\frac{1}{\sqrt{2}}a_{\sigma_k,\sigma_{n+1}}},
\\
\tilde{b}_{k,n+1}=
\displaystyle{\frac{1}{\sqrt{2}}a_{\sigma_{k'},\sigma_{n+1}}}, &
\tilde{b}_{n+1,k}=
\displaystyle{\frac{1}{\sqrt{2}}a_{\sigma_{n+1},\sigma_{k'}}}, \\
b_{kk}=\displaystyle{\frac{1}{2}
(a_{\sigma_k\sigma_k}+a_{\sigma_{k'}\sigma_{k'}})}, &
\tilde{b}_{kk}=\displaystyle{\frac{1}{2}
(a_{\sigma_{k}\sigma_{k'}}-a_{\sigma_{k'}\sigma_{k}})}, \\
b_{k'k}=\displaystyle{\frac{1}{2}
(a_{\sigma_k\sigma_k}-a_{\sigma_{k'}\sigma_{k'}})}, &
\tilde{b}_{k'k}=\displaystyle{\frac{1}{2}
(a_{\sigma_k\sigma_{k'}}+a_{\sigma_{k'}\sigma_{k}})}, \\
b_{kp}=\displaystyle{\frac{1}{2}a_{\sigma_k\sigma_p}}, \quad
b'_{kp}=\displaystyle{\frac{1}{2}a_{\sigma_{k'}\sigma_{p'}}}, &
\tilde{b}_{kp}=\displaystyle{\frac{1}{2}a_{\sigma_k\sigma_{p'}}},
\quad
\tilde{b}'_{kp}=\displaystyle{\frac{1}{2}a_{\sigma_{k'}\sigma_p}},
 \;  k \neq p. \\
\end{array}
\end{displaymath}

Let us observe that the elements $b$ of $ B_{\sigma}(\bj)$
are obtained from the elements $ b^* $ of $ B_{\sigma}(\bj=1)$
by multiplications on some products of  parameters $\bj,$ namely
\begin{equation}
     \begin{array}{ll}
b^*_{n+1,n+1}=b_{n+1,n+1}, &
b^*_{kk}=b_{kk}, \quad
b^*_{k'k}=b_{k'k}, \\
\tilde{b}^*_{kk}=(\sigma_k,\sigma_{k'})\tilde{b}_{kk}, &
\tilde{b}^*_{k'k}=(\sigma_k,\sigma_{k'})\tilde{b}_{k'k}, \\
b^*_{k,n+1}=(\sigma_k,\sigma_{n+1})b_{k,n+1}, &
b^*_{n+1,k}=(\sigma_k,\sigma_{n+1})b_{n+1,k},  \\
\tilde{b}^*_{k,n+1}=(\sigma_{k'},\sigma_{n+1})\tilde{b}_{k,n+1}, &
\tilde{b}^*_{n+1,k}=(\sigma_{k'},\sigma_{n+1})\tilde{b}_{n+1,k}, \\
b^*_{kp}=(\sigma_k,\sigma_p)b_{kp}, &
b^{*'}_{kp}=(\sigma_{k'},\sigma_{p'})b'_{kp}, \\
\tilde{b}^*_{kp}=(\sigma_k,\sigma_{p'})\tilde{b}_{kp}, &
\tilde{b}^{*'}_{kp}=(\sigma_{k'},\sigma_p)\tilde{b}'_{kp},
\;  k \neq p.  \\
\end{array}
   \label{6}
\end{equation}
A transformation of group by multiplications of some or all its
group parameters on zero tending parameter $\epsilon$ is  named
as group contraction \cite{GMClasG,GMSBook}, if a new group is obtained in the limit.
The formulas (\ref{6}) are just an example of such transformation,
where the nilpotent values $ j_{k}=\iota_{k} $ are used instead of the
limit $ \epsilon \rightarrow 0. $
In other words group contractions may be described mathematically
correctly by the replacement of real or complex group parameters
with a new one's which are elements of Pimenov algebra
$ {\cal D}_n(\iota; {\Bbb C}). $
In our case such replacement is made for matrix elements.

Let us consider as an example the group
 $ SO(3;\bj;\Bbb C).  $ For identical permutation  $ \sigma =(1,2,3) $
the matrix $ D_{\sigma} $ is given by equation (\ref{2}) for $ N=3 $
and in symplectic basis the group $ SO(3;\bj;\Bbb C) $
is described by the matrices
  \begin{displaymath}
        B_{\sigma}(j) =
     \left( \begin{array}{ccc}
 b_{11}+ij_{1}j_{2}\tilde{b}_{11} & j_{1}b_{12}-ij_{2}\tilde{b}_{12}
& b_{31}-ij_{1}j_{2}\tilde{b}_{31} \\
j_{1}b_{21}+ij_{2}\tilde{b}_{21}  & b_{22}
& j_{1}b_{21}-ij_{2}\tilde{b}_{21}       \\
 b_{31}+ij_{1}j_{2}\tilde{b}_{31} & j_{1}b_{12}+ij_{2}\tilde{b}_{12}
& b_{11}-ij_{1}j_{2}\tilde{b}_{11} \\
     \end{array} \right) .
 \end{displaymath}
For $ \sigma =(2,1,3)$ one obtain from equation (\ref{4})
 \begin{displaymath}
        B_{\sigma}(j) =
     \left( \begin{array}{ccc}
 b_{11}+ij_{2}\tilde{b}_{11} & j_{1}b_{12}-ij_{1}j_{2}\tilde{b}_{12}
& b_{31}-ij_{2}\tilde{b}_{31} \\
j_{1}b_{21}+ij_{1}j_{2}\tilde{b}_{21}  & b_{22}
& j_{1}b_{21}-ij_{1}j_{2}\tilde{b}_{21}       \\
 b_{31}+ij_{2}\tilde{b}_{31} & j_{1}b_{12}+ij_{1}j_{2}\tilde{b}_{12}
& b_{11}-ij_{2}\tilde{b}_{11} \\
     \end{array} \right) ,
 \end{displaymath}
finally the permutation
     $ \sigma =(1,3,2) $
leads to the matrices
 \begin{displaymath}
        B_{\sigma}(j) =
     \left( \begin{array}{ccc}
 b_{11}+ij_{1}\tilde{b}_{11} & j_{1}j_{2}b_{12}-ij_{2}\tilde{b}_{12}
& b_{31}-ij_{1}\tilde{b}_{31} \\
j_{1}j_{2}b_{21}+ij_{2}\tilde{b}_{21}  & b_{22}
& j_{1}j_{2}b_{21}-ij_{2}\tilde{b}_{21}       \\
 b_{31}+ij_{1}\tilde{b}_{31} & j_{1}j_{2}b_{12}+ij_{2}\tilde{b}_{12}
& b_{11}-ij_{1}\tilde{b}_{11} \\
     \end{array} \right) .
 \end{displaymath}
The same matrices are corresponded to three
remaining permutations from the group $S(3;\bj;\Bbb C)$.

For nilpotent values of both parameters
$ j_{1}=\iota_{1}, j_{2}=\iota_{2} $ we have the complex Galilei
group $ G(1+1;\Bbb C)=SO(3;\iota;\Bbb C), $ which is realized in Cartesian
basis by the matrices
 \begin{displaymath}
        A(\iota) =
     \left( \begin{array}{ccc}
 1  & \iota_{1}a_{12} & \iota_{1}\iota_{2}a_{13} \\
-\iota_{1}a_{12}  & 1 & \iota_{2}a_{23}       \\
 \iota_{1}\iota_{2}a_{31} & -\iota_{2}a_{23} & 1 \\
     \end{array} \right) ,
 \end{displaymath}
where $ a_{31}= -a_{13}+a_{12}a_{23}. $
The relations of $ \bj$-orthogonality (\ref{3})
have been taken into account. Three different realizations of
Galilei group in symplectic description are as follows
 \begin{displaymath}
B_{\sigma}(\iota) =   \left( \begin{array}{ccc}
 1+i\iota_{1}\iota_{2}\tilde{b}_{11}   &
\iota_{1}b_{12}-i\iota_{2}\tilde{b}_{12}
& -i\iota_{1}\iota_{2}\tilde{b}_{31} \\
-\iota_{1}b_{12}-i\iota_{2}\tilde{b}_{12}  & 1
& -\iota_{1}b_{12}+i\iota_{2}\tilde{b}_{12}       \\
 i\iota_{1}\iota_{2}\tilde{b}_{31} & \iota_{1}b_{12}+i\iota_{2}\tilde{b}_{12}
& 1-i\iota_{1}\iota_{2}\tilde{b}_{11} \\
     \end{array} \right) ,
 \end{displaymath}
where $ \tilde{b}_{31}=-b_{12}\tilde{b}_{12}, $
 \begin{displaymath}
        B_{\sigma}(\iota) =
     \left( \begin{array}{ccc}
 1+i\iota_{2}\tilde{b}_{11}   &
\iota_{1}b_{12}-i\iota_{1}\iota_{2}\tilde{b}_{12} & 0 \\
-\iota_{1}b_{12}+i\iota_{1}\iota_{2}\tilde{b}_{21}  & 1
& -\iota_{1}b_{12}-i\iota_{1}\iota_{2}\tilde{b}_{21}       \\
0 & \iota_{1}b_{12}+i\iota_{1}\iota_{2}\tilde{b}_{12}
& 1-i\iota_{2}\tilde{b}_{11} \\
     \end{array} \right) ,
 \end{displaymath}
where $ \tilde{b}_{21}=-\tilde{b}_{12} -b_{12}\tilde{b}_{11}, $
 \begin{displaymath}
        B_{\sigma}(\iota) =
     \left( \begin{array}{ccc}
 1+i\iota_{1}\tilde{b}_{11}   &
\iota_{1}\iota_{2}b_{12}-i\iota_{2}\tilde{b}_{12} & 0 \\
\iota_{1}\iota_{2}b_{21}-i\iota_{2}\tilde{b}_{12}  & 1
& \iota_{1}\iota_{2}b_{21}+i\iota_{2}\tilde{b}_{12}       \\
0 & \iota_{1}\iota_{2}b_{12}+i\iota_{2}\tilde{b}_{12}
& 1-i\iota_{1}\tilde{b}_{11} \\
     \end{array} \right) ,
 \end{displaymath}
where $ b_{21}=-b_{12} + \tilde{b}_{11}\tilde{b}_{12}. $

\subsection{ Complex symplectic Cayley-Klein groups   }

\noindent

Consider the space $ {\Bbb C^n}(\bj) \otimes {\Bbb C^n}(\bj) $, which is obtained
from $ 2n$-dimensional Euclidean space $ {\Bbb C^n} \otimes {\Bbb C^n} $
by mapping
\begin{equation}
\label{2.16} 
 \psi: \enskip {\Bbb C}^n \otimes {\Bbb C}^n \rightarrow {\Bbb C}^n(\bj) \otimes {\Bbb C^n }(\bj) , \\
\end{equation}
\beq
\psi  z_k = z_kJ_k, \quad  \psi z_{k^{\prime}} = z_{k^{\prime}}J_k, \quad
k^{\prime} = 2n + 1 - k, \enskip
J_k= \prod_{m=1}^{k-1} j_m ,
\enskip k=1,2 \ldots n . \nonumber
\eeq
Group $ Sp(n;\bj;\Bbb C) $ is defined as the group of transformation
(\ref{2.16}) of
 $ 2n-$dimensional space
$ {\Bbb C}^n(\bj) \otimes {\Bbb C}^n(\bj) $,
preserving the bilinear form
\beq
(z,y)= \sum_{k=1}^{n}(z_ky_{k^{\prime}}-z_{k^{\prime}}y_k)J_k^2.
\eeq
Cartesian  coordinates $ z_k, y_k, \enskip k=1,2 \ldots n, $
belongs to the first, and  $ z_{k^{\prime}},y_{k^{\prime}} $---
to the second factor in the direct product of spaces.
For $ \bj $-structure to be preserved, parameter $ \bj $ must be introduced
in matrix of transformation $ T(\bj) = (T_{ij})^n_{ i,j = 1} $ as follows
\beq
 T_{ij} =  J^k_m t_{ij} ,\  {\rm where} \enskip  k = i \enskip if \enskip i \leq n;
 \quad k = i^{\prime}, \enskip i > n;
 \quad m = j,  \enskip j \leq n;
\eeq
\beq
 m = j^{\prime} , \enskip j > n, \qquad
 J^k_m = \prod_{n=min(k,m)}^{max(k,m)-1} j_n .
\eeq
The matrix elements of $ T $ satisfies the additional relation
  \beq
     T^tC_0T=C_0.
  \eeq
The elements of matrix $ C_0 $ are as follows
$ (C_0)_{ij} = \varepsilon_i \delta_{ij} , \enskip
 \varepsilon_i = 1 $ if $ i=1,2 \ldots n $ and
$ \varepsilon_{i} = -1 $ for $ i = n-1, \ldots , 2n $.
The general element of the Shevalley basis of symplectic algebra
$ sp(n;\bj;\Bbb C) $ appears as
\beq
 H= \left (\begin{array}{cccccc}
   h_{1}    & j_1 z_1^+ &  \cdot      & \cdots    & \cdots  & \cdot      \cr
j_1 z_1^-   & h_2       &  j_2 z_2^+  & \cdots    & \cdots  & \cdot      \cr
   \cdot    & j_1 z_2^+ & \cdots      & \cdots    & \cdots  & \cdot      \cr
   \cdot    & \cdots    & \cdots      & \cdots    & \cdots  & \cdots     \cr
   \cdot    & \cdots    & \cdots      &-j_2 z_2^- & - h_1   & -j_1 z_1^+  \cr
   \cdot    & \cdots    & \cdots      & \cdots    & - j_1 z_1^- & - h_2
            \end{array} \right ) .
\eeq
Choosing the generators of algebra $ sp(n;\bj;\Bbb C) $ in the following form
 \begin{eqnarray}
 H_i &  = &  e_{ii} - e_{ i^{\prime}i^{\prime}} - e_{i+1,i+1}
 +  e_{(i+1)^{\prime},(i+1)^{\prime}} ,\quad
 H_n = e_{nn} - e_{n^{\prime},n^{\prime}} ,    \nonumber \\
 X_i^{+} &  = & j_i(e_{i,i+1} - e_{(i+1)^{\prime},i^{\prime}})  ,\quad
 X_i^{-} = j_i(e_{i+1,i} - e_{(i)^{\prime},(i+1)^{\prime}} )  ,   \nonumber \\
 X_n^{+} & = & j_n  e_{n,n+1} ,\quad
 X_n^{-} = j_n e_{n+1,n} ,  \quad i = 1, \ldots , n-1 ,
 \end{eqnarray}
where $ (e_{ij})_{km}=\delta_{ik}\delta_{jm}, $
the commutation relation for the Shevalley basis $ sp(n;\bj;\Bbb C) $
looks as
 \begin{eqnarray}
\label{2.18}
 \left [H_i,X_j^{ \pm} \right ]   =  \pm A_{ij} X_j^{ \pm} ,\quad 
 \left [X_i^{+},X_j^{-} \right ]  =  \delta_{ij} j_i^2 H_j , \quad
 \end{eqnarray}
where $ A_{ij} $ is a Cartan matrix ---
$ A_{ii} = 2, \enskip  A_{i,i-1} = A_{i-2,i-1} = -1, \enskip
A_{n-1,n} = - 2 $.

\setcounter{equation}{0}

\section{Structure of   classical Cayley-Klein algebras \\under 
contractions}

\subsection{Structure of Cayley-Klein orthogonal algebras}
\medskip

 In the theory
of Lie algebras, Cartan-Weyl commutation relations, connected with root
expansions of classical Lie algebras, are also often used.
To start with,
we shall give a brief description of the main stages of a root expansion,
following the monographs \cite{Knapp} and then construct analogues of such
expansion for algebras $so\,(n+1;\bj)$, which enables us to find
the structure of
contracted algebras.

Under Cartan subalgebra $h$ of a simple algebra $l$ we mean maximal Abelian
subalgebra in $l$. Let $\alpha$ be a linear function on $h\subset l$.
Let us denote
$E_\alpha$ linear subspace in $l$, defined by condition
\begin{equation}
\label{1.56}
E_\alpha=\{Y\in l\mid [X,Y]=\alpha(X)Y\quad\forall X\in h\}.
\end{equation}

If $E_\alpha\neq 0$ the function $\alpha$ is called a root,
and $E_\alpha$ is root vector. A set
of all roots is denoted $\Delta$. Algebra $l$ can be presented
as direct sum

\begin{equation}
\label{1.57}
l=h\oplus\sum_{\alpha\in\Delta}\!\!{}^\oplus\{E_\alpha\}.
\end{equation}

Let $H_k$ is a basis of Cartan subalgebra $h$.  Then for each root
$\alpha\in\Delta$ it is possible to choose vector $E_\alpha\in l$ such
that for all $\alpha,\beta\in\Delta$ the following commutation
relations are valid:
\beq
[H_k,E_\alpha]=\alpha(H_k)E_\alpha,
\eeq
\begin{equation}
\label{1.58}
[E_\alpha,E_\beta]=\left\{ \ba{cc}
 H_\alpha,&\alpha+\beta=0,\\
N_{\alpha\beta}E_{\alpha+\beta},&\alpha+\beta\in\Delta,\\
0,&\alpha+\beta\neq 0,\quad\alpha+\beta\notin\Delta,\ea
\right.
\end{equation}
They are called Cartan-Weyl commutators.

In general, root vectors
$E_\alpha$
 for $\alpha\in\Delta$
are linearly dependent. For this reason,
from the set $\Delta$ of all roots one singles out a subject
$\Pi$
 of
simple roots to which there correspond linearly independent root
vectors $E_\alpha$. Then to each root
$\alpha\in\Pi$
 one puts in correspondence an
element
$H_\alpha\in h$,
 defines the internal product of roots as
$\langle \alpha_k,\alpha_m\rangle=\Tr(H_{\alpha_k}H_{\alpha_m})$
  and takes into
consideration Cartan matrix

\begin{equation}
\label{1.59}
A_{km}={2\langle\alpha_k,\alpha_m\rangle\over
\langle\alpha_k,\alpha_k\rangle}.
\end{equation}

In the case of classical algebras Cartan
matrix has following properties:

\noindent
1) $A_{kk}=2$; 2) $A_{km}\leq0$,
$A_{km}=0,-1,-2,-3$ or $-4$, if $k\neq m$; 3) $A_{km}A_{mk}<4$,
$k\neq m$; 4) $A_{km}=0$ if and only if $A_{mk}=0$; 5) $\det(A_{km})$
is positive integer.  The structure of Cartan matrix can be
conveniently depicted graphically, using Dynkin diagrams. To each
 simple root $\alpha_k\in\Pi$ there corresponds a point in a plane
with a weight proportional to $\langle\alpha_k,\alpha_k\rangle$. Any
two  points $\alpha_k$ and $A_m$, $k\neq m$, are connected with lines.
It turns out that there are four infinite Dynkin diagrams,
corresponding to infinite series of classical Lie algebras:

\begin{equation}
\label{1.60} \noindent
(A_n),\quad u(n),\quad \Pi=\{e_{k-1}-e_k,\quad
k=1,2,\dots,n\}
\end{equation}

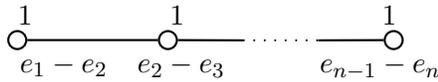
\begin{center}
\begin{picture}(1320,40)

\cnodeput(0,1){A}{}
 \cnodeput(2,1){B}{}
 \pnode(3,1){C}
 \cnodeput(5,1){D}{}
 \pnode(4,1){H}
 \ncline {-}{A}{B}\Aput{1 \ \qquad\qquad} \Bput{$e_1-e_2 \qquad$}
 \ncline {-}{B}{C}\Aput{1 \qquad}\Bput{$e_2-e_3 \qquad$}
 \ncline[linestyle=dotted] {-}{H}{C}
 \ncline {-}{H}{D}\Aput{\qquad \quad 1  }\Bput{$\qquad e_{n-1}-e_n$}

\end{picture}
\end{center}


\begin{equation}
\label{1.61} \noindent
(B_n),\quad so\,(2n+1),\quad\Pi=\{e_k-e_{k+1},\quad k=1,2,\dots,
n-1,e_n\},
\end{equation}

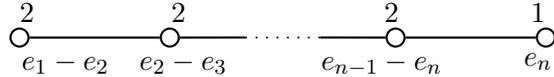
\begin{picture}(1320,40)

\cnodeput(0,1){A}{}
 \cnodeput(2,1){B}{}
 \pnode(3,1){C}
 \cnodeput(5,1){D}{}
 \cnodeput(7,1){E}{}
 \pnode(4,1){H}
 \ncline {-}{A}{B}\Aput{2 \ \qquad\qquad} \Bput{$e_1-e_2 \qquad$}
 \ncline {-}{B}{C}\Aput{2 \qquad}\Bput{$e_2-e_3 \qquad$}
  \ncline[linestyle=dotted] {-}{H}{C}
 \ncline {-}{H}{D}\Aput{\qquad \quad 2  }\Bput{$\qquad e_{n-1}-e_n$}
 \ncline {-}{D}{E}\Aput{\quad \qquad \qquad 1  }\Bput{$\ \ \qquad \qquad e_n$}

\end{picture}

\begin{equation}
\label{1.62} \noindent
(C_n),\quad sp\,(n),\quad\Pi=\{e_k-e_{k+1},\quad
k=1,2,\dots,n-1,2e_n\}
\end{equation}

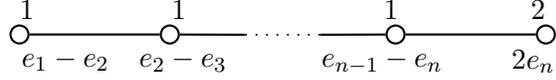
\begin{picture}(1320,40)

 \cnodeput(0,1){A}{}
 \cnodeput(2,1){B}{}
 \pnode(3,1){C}
 \cnodeput(5,1){D}{}
 \cnodeput(7,1){E}{}
 \pnode(4,1){H}
 \ncline {-}{A}{B}\Aput{1 \ \qquad\qquad} \Bput{$e_1-e_2 \qquad$}
 \ncline {-}{B}{C}\Aput{1 \qquad}\Bput{$e_2-e_3 \qquad$}
  \ncline[linestyle=dotted] {-}{H}{C}
 \ncline {-}{H}{D}\Aput{\qquad \quad 1  }\Bput{$\qquad e_{n-1}-e_n$}
 \ncline {-}{D}{E}\Aput{\quad \qquad \qquad 2  }\Bput{$\ \qquad \qquad 2e_n$}

\end{picture}


\begin{equation}
\label{1.63} \noindent
(D_n),\quad so\,(2n),\quad\Pi=\{e_k-e_{k+1},\quad
k=1,2,\dots,n-1,e_{n-1}+e_n\}
\end{equation}

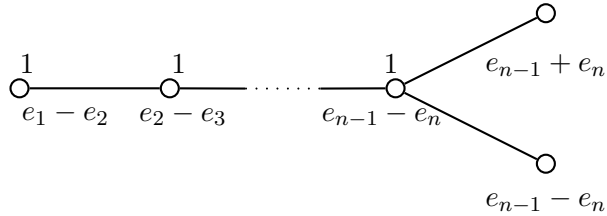
\begin{picture}(1320,90)

 \cnodeput(0,1){A}{}
 \cnodeput(2,1){B}{}
 \pnode(3,1){C}
 \cnodeput(5,1){D}{}
 \pnode(4,1){H}
 \cnodeput(7,0){Y}{}
 \cnodeput(7,2){Z}{}
 \ncline{-}{D}{Y} \rput(7,-0.5){$e_{n-1}-e_n$}
 \ncline{-}{D}{Z} \rput(7,1.3){$e_{n-1}+e_n$}
 \ncline {-}{A}{B}\Aput{1 \ \qquad\qquad} \Bput{$e_1-e_2 \qquad$}
 \ncline {-}{B}{C}\Aput{1 \qquad}\Bput{$e_2-e_3 \qquad$}
 \ncline[linestyle=dotted] {-}{H}{C}
 \ncline {-}{H}{D}\Aput{\qquad \quad 1  }\Bput{$\qquad e_{n-1}-e_n$}

\end{picture}
\vspace{1cm}

The number over circles at the diagrams (\ref{1.60}) -- (\ref{1.63})
indicate weights, Cartan notations are given in brackets, vectors
$e_k$ are orthogonal basis vectors of Euclidean space and have the
same (though arbitrary) length.  In addition to infinite sequences of
classical algebras, there are five exceptional Lie algebras. The
latter will not be considered here.

Let us first study algebras
$so\,(n+1;\bj)$
 of even dimension $n=2m$. A matrix $H$ from
Cartan subalgebra $h$ can be presented as
\begin{equation}
\label{1.64}
H=\left(\ba{cccccccc}
0&ij_1^2h_1&&&\ldots&&&0\\
-ih_1&0&&&&&&\vdots\\
&&\ddots&&&&&\\
\vdots&&&0&ij_{2k-1}^2h_k&&&\\
&&&-ih_k&0&&&\\
&&&&&\ddots&&\\
&&&&&&0&ij_{n-1}^2h_m\\
0&&\ldots&&&&-ih_m&0 
\ea\right)
\end{equation}

\noindent
Let matrix $A_{kr}$ has zero elements except for element
$(A_{kr})_{kr}=1$. As a
basis in Cartan subalgebra $h$ let us choose matrixes
 $H_k=-i
X_{2k-2,2k-1}= -i(A_{2k-1,2k-2} - j_{2k-1}^2 A_{2k-2,2k-1})$,
$k=1,2,\dots,m$ (we write out the constructions for $so\,(n+1;\bj)$;
the corresponding constructions for classical algebra $so\,(n+1)$ can
be obtained, putting all $j_k=1$). Linear functions on Cartan
 subalgebra are defined by $e_k(H)=h_k$.  Simple roots are given in
(\ref{1.61}). The system of all roots is as follows
 \bea
 \Delta &\!\!\!=\!\!\!& \{\pm
e_r\pm e_s,\quad r\neq s,\quad r,s=1,2,\dots,m\}\cup
\nonumber \\
  &\!\!\!\cup\!\!\!&
 \{\pm
e_k,\quad k=1,2,\dots,m\}.  
\label{1.65}
\eea
 To the roots
$\pm e_k$ there
correspond generators (root vectors)
 \begin{equation}
\label{1.66}
E_{\pm e_k}=-X_{2k-2,n}\pm ij_{2k-1}X_{2k-1,n}, 
 \end{equation}
to the roots
$e_r\pm e_s$ -- generators
\bea
\label{1.67}
E_{e_r\pm e_s} &\!\!\!=\!\!\!&
 -j_{2s-1}X_{2r-2,2s-2}\pm iX_{2r-2,2s-1}\pm
j_{2r-1}X_{2r-1,2s-1}+\nonumber \\
  &\!\!\!+\!\!\!&
ij_{2r-1}j_{2s-1}X_{2r-1,2s-2},
 \eea
to the roots
$-e_r\pm e_s$
 -- generators
 \bea
\label{1.68}
E_{-e_r\pm e_s} &\!\!\!=\!\!\!&
 -j_{2s-1}X_{2r-2,2s-2}\pm iX_{2r-2,2s-1}\mp\nonumber \\
  &\!\!\!\mp\!\!\!&
j_{2r-1}X_{2r-1,2s-1}-ij_{2r-1}j_{2s-1}X_{2r-1,2s-2}.
 \eea

Matrix generators $X_{\mu\nu}$ are described as follows \cite{GMOrtG}:
\be
(X_{\mu\nu})_{\nu\mu}=1,\quad (X_{\mu\nu})_{\mu\nu}=
-\prod_{m=\mu+1}^\nu j_m^2.
\ee

Cartan-Weyl commutators (\ref{1.58}) can be found using
commutation relations for Lie algebra $so\,(n+1;\bj)$, which can be of
the most simply derived from commutators of algebra $so\,(n+1)$, as it
has been done in \cite{GMOrtG}. The non-zero commutators are
\bea
&&[X_{\mu_1\nu_1}, X_{\mu_2\nu_2}]=\\
 &\!\!\!=\!\!\!& 
\left\{ \ba{ll}
\biggl(\prod\limits_{m=\mu_1+1}^{\nu_1} j_m^2\biggr)
X_{\nu_1\nu_2},& \mu_1=\mu_2,\quad \nu_1<\nu_2,\\
\biggl(\prod\limits_{m=\mu_2+1}^{\nu_2} j_m^2\biggr)
X_{\mu_1\mu_2},& \mu_1<\mu_2,\quad \nu_1=\nu_2,\\
-X_{\mu_1\nu_2},& \mu_1<\mu_2=\nu_1<\nu_2.
\ea \right. \nonumber 
\eea

So, Cartan-Weyl commutators (\ref{1.58}) can be written  as
follows (only non-zero commutators are written out):

\bea
\label{1.69}
\left[H_k,E_{\pm e_k}\right] &\!\!\!=\!\!\!&
\pm j_{2k-1}E_{\pm e_k},\nonumber \\
\left[H_k,E_{\pm e_r\pm e_s}\right] &\!\!\!=\!\!\!&
\left\{ \ba{ll}
\pm j_{2r-1}E_{\pm e_r\pm e_s},& k=r,\\
\pm j_{2s-1}E_{\pm e_r\pm e_s},& k=s, \ea \right. \\
\left[E_{\pm e_k},E_{\mp e_k}\right]
 &\!\!\!=\!\!\!& \mp 2j_{2k-1}\biggl(\prod_{l=2k}^{n}
j_l^2\biggr)H_k;\nonumber
\eea
\bea
\label{1.70}
\left[E_{e_r\pm e_s},E_{-e_r\mp e_s}\right]
&\!\!\!=\!\!\!& -4j_{2r-1}j_{2s-1}
\biggl(\prod_{l=2r}^{2s-2}j_l^2\biggr)\times\nonumber \\
 &\!\!\!\times\!\!\!&
(j_{2s-1}H_r\pm j_{2r-1}H_s),\qquad r<s;
\eea
\beq
\left[E_{\pm r_s},E_{\pm e_s}\right] &\!\!\!=\!\!\!& -
j_{2s-1}\biggl(\prod_{l=2s}^n j^2_l\biggr)E_{\pm e_r\pm e_s},\qquad
r<s,\nonumber \\
\left[E_{\pm e_r\pm e_s},E_{\mp e_s}\right] &\!\!\!=\!\!\!&
2j_{2s-1}E_{\pm e_r},\nonumber \\
\left[E_{\pm e_r\pm e_s},E_{\mp e_r}\right]
&\!\!\!=\!\!\!& -2j_{2s-1}
\biggl(\prod_{l=2r-1}^{2s-2}j_l^2\biggr)E_{\pm e_s};
\eeq
\bea
\label{1.71}
&&
\left[E_{\pm e_r\pm e_s},E_{\mp e_s\pm e_m}\right]
=\nonumber \\
 &\!\!\!=\!\!\!& 
\left\{ \ba{cc}
2j_{2s-1}E_{\pm e_r\pm e_m},&r<s<m,\\
-2j_{2m-1}\biggl(\prod\limits_{l=2m}^{2s-1}j_l^2\biggr)E_{\pm e_r\pm
e_m},&r<m<s,\\
-2j_{2r-1}\biggl(\prod\limits_{l=2s-1}^{2r-2}j_l^2\biggr)E_{\pm e_r\pm
e_m},&s<r<m.
\ea\right.
\eea 

The root technique is very propriate for description of contractions of
algebra
$so\,(n+1;\bj)$
 over parameters with even numbers
$j_2=\iota_2$, $j_4=\iota_4$, $\dots$, $j_n=\iota_n$, i.e. over
parameters which are not involved in Cartan subalgebra. Consequently,
Cartan subalgebra is not changed under such contractions. The
structure of contracted algebras $so\,(n+1;\bj)$ is determined
by commutators (\ref{1.69}) -- (\ref{1.71})
and can be found by induction over
dimension of algebra. We shall skip the formal derivation and proceed
directly to the description of results. Doing so, it is convenient to
introduce an auxiliary table $\Gamma_{n+1}(\bj)$.

Let us begin with some examples. In the case of algebra
$so\,(3;\bj)$
 the commutators
(\ref{1.69}) -- (\ref{1.71}) give us
$[E_{e_1},E_{-e_1}]=-2j_1j_2^2H_1$.
Let us construct a table
$\Gamma_3(\bj)$, as follows:
$\Gamma_3(\bj)=\stackrel{{\overline{\bigm|H_1\bigm|j^2_2\bigm|}}}
{{}_{\overline{ \ \ \ss{1} \  \  \ \ss{2} \ }}}$,
i.e. to the element $(\Gamma_3)_{11}$
of the table there corresponds generators $H_1$, to the element
$(\Gamma_3)_{12}$  -- generators $E_{\pm e_1}$,and in the cell
$(\Gamma_3)_{12}$ we inscribe parameter $j_2^2$, arising when these two
generators commute (let us remind that we have considered contractions
 over parameters with even numbers, so that parameters with odd numbers
can be put equal to unit and omitted). For contraction $j_2=\iota_2$ we
get the table $\Gamma_3(\iota_2)$=\fbox{$H_1 \bigm| \, 0$\, },  and 
algebra $so\,(3;\iota_2)$ gets the structure of semidirect sum of 
commutative ideal $T_2$, spanned over generators $E_{\pm e_1}$, and 
one-dimensional subalgebra $H_1$, spanned over $H_1$, i.e.  
$so\,(3;\iota_2)=T_2\notni H_1$.  

\begin{definition} 
\label{d1.1}
Algebra $L$ is called semidirect sum of
subalgebras $T$ and
$M,L=T\notni M$, if
$[T,T]\subset T$, $[M,M]\subset M$, $[M,T]\subset T$.
\end{definition}

It is easy to check up that these
conditions are satisfied for
$T_2=\{E_{\pm e_1}\}$ and $H_1$.

\begin{definition} 
\label{d1.2}
A subalgebra $T$ of algebra $L$ is called
ideal, if $[L,T]\subset T$.
\end{definition}

In the case of algebra
$so\,(5;\bj)$
 the commutators (\ref{1.69}) -- (\ref{1.71}) give us
$(j_1=j_3=1)$
\bea
\label{1.72}
\left[E_{e_1},E_{-e_1}\right] &\!\!\!=\!\!\!& -2j_2^2j_4^2H_1,\quad
\left[E_{e_2},E_{-e_2}\right]=-2j_4^2H_2,\nonumber \\
\left[E_{e_1-e_2},E_{e_2-e_1}\right] &\!\!\!=\!\!\!& -4j_2^2(H_1-H_2),
\eea
and the table
$\Gamma_5(\bj)$
 is as follows

\begin{equation}
\vbox{\offinterlineskip
\halign{#\on&#&\on #\on&#&\on #\on&#&\,#\,&#&\on #\cr
&&\strut\on$\ss{1}$&&\on$\ss{2}$&&\on$\ss{3}$&&\cr
&\multispan6\hrulefill&\cr
\phantom{a}&\vrule&&\vrule&&\vrule&&\vrule&\cr
$\strut\Gamma_5(\bj)=$&\vrule&$
H_1$&\vrule&$\,j^2_2$&\vrule&$j^2_2j^2_4$& \vrule&$\ss{1}$\cr
\phantom{a}&\vrule&&\vrule&&\vrule&&\vrule&\cr
&\multispan6\hrulefill&\cr
\phantom{a}&&&\vrule&&\vrule&&\vrule&\cr
&&&\vrule&$\strut H_2$&\vrule&$\on j^2_4$&\vrule&$\ss{2}$\cr
\phantom{a}&&&\vrule&&\vrule&&\vrule&\cr
&&&\multispan4\hrulefill&\cr}}\ ,
\end{equation}
i.e. to elements
$(\Gamma_5)_{kk}$
 we put in correspondence generator $H_k$ of Cartan subalgebra,
 $k=1,2$; to the element
$(\Gamma_5)_{12}$
  -- four generators
$E_{\pm e_1\pm e_2}$
 and inscribe in this
cell the parameter $j_2^2$, arising under commutation of these
generators; to the element
$(\Gamma_5)_{13}$
 -- two generators
$E_{\pm e_1}$
 with parameters
$j_2^2j_4^2$, and to the element
$(\Gamma_5)_{23}$
 -- two generators
$E_{\pm e_2}$
 with parameters
$j_4^2$. Under contraction
$j_2=\iota_2$
 some of these elements vanish:


\begin{equation}
\label{1.74}
\vbox{\offinterlineskip
\halign {#\on&#&\on #\on&#&\on #\on&#&\on #\on&#\cr
&\multispan6\hrulefill\cr
\phantom{a}&\vrule&&\vrule&&\vrule&&\vrule\cr
$\strut\Gamma_5(\iota_2)=$&\vrule&$
H_1$&\vrule&$\on 0$&\vrule&$\, 0\,$& \vrule\cr
\phantom{a}&\vrule&&\vrule&&\vrule&&\vrule\cr
&\multispan6\hrulefill\cr
\phantom{a}&&&\vrule&&\vrule&&\vrule\cr
&&$\oplus$&\vrule&$\strut H_2$&\vrule&$\, j^2_4\,$&\vrule\cr
\phantom{a}&&&\vrule&&\vrule&&\vrule\cr
&&&\multispan4\hrulefill\cr}}\ ,
\end{equation}
 and
algebra
$so\,(5;\iota_2)$
 gets the structure of semidirect sum
\begin{equation}
\label{1.75}
so\,(5;\iota_2)=T_6\notni (H_1\oplus so\,(3;j_3,j_4)),
\end{equation}
 where Abelian
subalgebra $T_6$ is spanned over generators corresponding to zero
cells of the table, i.e.
 $T_6=\{E_{\pm e_1\pm e_2,}E_{\pm e_1}\}$, and subalgebra $M$, 
according to the definition \ref{d1.1}, is a semidirect sum of 
one-dimensional subalgebra, spanned over generator $H_1$, and 
subalgebra $so\,(3;j_3,j_4)=$ \linebreak $=\{H_2,E_{\pm e_2}\}$, 
corresponding to the cells $(\Gamma_5)_{22}$ and $(\Gamma_5)_{23}$ (in 
 (\ref{1.74}) the latter fact is marked with the sign $\oplus$).

For contraction
$j_4=\iota_4$
 the table and the structure of algebra
$so\,(5;\iota_4)$
  are as follows
 \begin{equation}
\label{1.76}
\vbox{\offinterlineskip
\halign {#\on&#&\on #\on&#&\on #\on&#&\on #\on&#\cr
&\multispan6\hrulefill\cr
\phantom{a}&\vrule&&\vrule&&\vrule&&\vrule\cr
$\strut\Gamma_5(\iota_4)=$&\vrule&$
H_1$&\vrule&$\on j_2^2$&\vrule&$\, 0\,$& \vrule\cr
\phantom{a}&\vrule&&\vrule&&\vrule&&\vrule\cr
&\multispan6\hrulefill\cr
\phantom{a}&&&\vrule&&\vrule&&\vrule\cr
&&&\vrule&$\strut H_2$&\vrule&$\, 0\,$&\vrule\cr
\phantom{a}&&&\vrule&&\vrule&&\vrule\cr
&&&\multispan4\hrulefill\cr}}\ ,\quad
so\,(5;\iota_4)=T_4\notni so\,(4;\bj),
 \end{equation}
where Abelian ideal
$T_4=\{E_{\pm e_1},E_{\pm e_2}\}$
 and to non-zero part of the table
$\Gamma_5(\iota_4)$
there corresponds to subalgebra
$so\,(4;j_1,j_2,j_3)$, as we shall see further (see (\ref{1.84}))  
considering orthogonal algebras of even dimension. It must be noted 
that for one-dimensional contractions (when one parameter $j_k$ takes 
dual value) the ideal $T$ is always communicative. But this is not 
valid for multidimensional contractions (when two or more parameters 
$j_k$ are equal to dual units). For dual contraction $j_2=\iota_2$, 
$j_4=\iota_4$ of algebra $so\,(5)$ we get

 \begin{equation}
\label{1.77}
\vbox{\offinterlineskip
\halign {#\on&#&\on #\on&#&\on #\on&#&\on #\on&#\cr
&\multispan6\hrulefill\cr
\phantom{a}&\vrule&&\vrule&&\vrule&&\vrule\cr
$\strut\Gamma_5(\iota_2,\iota_4)=$&\vrule&$
H_1$&\vrule&$\on 0$&\vrule&$\ 0\ $& \vrule\cr
\phantom{a}&\vrule&&\vrule&&\vrule&&\vrule\cr
&\multispan6\hrulefill\cr
\phantom{a}&&&\vrule&&\vrule&&\vrule\cr
&&$\oplus$&\vrule&$\strut H_2$&\vrule&$\ 0\ $&\vrule\cr
\phantom{a}&&&\vrule&&\vrule&&\vrule\cr
&&&\multispan4\hrulefill\cr}}\ ,\quad
so\,(5;\iota_2,\iota_4)=T_8\notni H_1\oplus H_2,
 \end{equation}
where ideal
$T_8=\{E_{\pm e_1\pm e_2},E_{\pm e_1},E_{\pm e_2}\}$
 is not commutative, because, for example,
$[E_{\pm e_1\pm e_2},E_{\mp e_2}]=2E_{\pm e_1}$, but it
can be easily seen that $T_8$ is a nilpotent subalgebra and, moreover,
$T_8$ is a radical.

Let $L$ is Lie algebra. Let us introduce sequences of radicals
 \bea
\label{1.78}
L^{(0)}  &\!\!\!=\!\!\!& L,\quad L^{(1)}=[L^{(0)},L^{(0)}],\dots ,
L^{(n+1)}=[L^{(n)},L^{(n)}],\nonumber \\
L_{(0)} &\!\!\!=\!\!\!& L,
 \quad L_{(1)}=[L_{(0)},L],\dots,L_{(n+1)}=[L_{(n)},L],\nonumber \\
n &\!\!\!=\!\!\!& 0,1,2,\dots
 \eea
 
\begin{definition}
\label{d1.3}
 Algebra $L$ is called solvable, if for some
positive integer $n$
$L^{(n)}=0$.
 \end{definition}

 \begin{definition} 
\label{d1.4}
Algebra $L$ is called nilpotent, if for some
positive integer $n$
$L_{(n)}=0$.
 \end{definition}

A nilpotent algebra is solvable. The inverse is not true.

\begin{definition}  
\label{d1.5}
Maximal solvable ideal $T$, involving any other
solvable ideal of algebra $L$, is called a radical.
 \end{definition}

The structure of the radical $T_8$ is not fixed uniquely. It is
determined by partitions of zero elements of the table
$\Gamma_5(\iota_2,\iota_4)\equiv\Gamma(\biota)$
 into blocks,
generated by different consequent one-dimensional contractions:
contraction ``first
$j_2=\iota_2$, then $j_4=\iota_4$''gives

\begin{equation}
\label{1.79}
\ba{cclrr}
& \ss{2}  \ \ \  \ss{3} && & \cr 
& \frf{\phantom{H}\bigm|\phantom{H}} & \ss{1}& &\cr
&\quad \oplus &&& T_8=T_6\notni T_2,\cr
& \quad\, \frf{\phantom{\bigm|a }} & 
\ss{2}&& \cr
\ea
 \end{equation}
 where Abelian subalgebra
$T_6=\{E_{\pm e_1\pm e_2},E_{\pm e_1}\}$
  is spanned
over generators corresponding to the cells (1, 2) and (1, 3) of the
table (\ref{1.79}),
$T_2=\{E_{\pm e_2}\}$
 - over generators corresponding to the cell (2, 3).
The sequence of contractions ``first $j_4=\iota_4$, then
$j_2=\iota_2$'' gives
\begin{equation}
\label{1.80}
\vbox{\offinterlineskip
\halign {#&\on #\on&#&\,#\,&#&\on #&#&\on #\on&\on #\cr
\phantom{$\bigm| $}&$\ss{2}$&&&& $\ss{3}$&&& \cr 
\multispan2\hrulefill&&&&&&\cr
&&&&\multispan2\hrulefill&&\cr
\vrule&&\vrule&&\vrule&&\vrule&&\cr
\vrule&\phantom{$\bigm|a $}&\vrule&$\oplus$&\vrule&
\phantom{$\bigm|a $}&\vrule&$\ss{1}$&$T_8=T_4\oplus
\widetilde{T}_4$,\cr
\vrule&&\vrule&&\vrule&&\vrule&&\cr
\multispan2\hrulefill&&&&&&\cr
&&&&\multispan2\hrulefill&&\cr
&&&&\vrule&&\vrule&&\cr
&&&&\vrule&\phantom{$\int $}&\vrule&$\ss{2}$&\cr
&&&&\vrule&&\vrule&&\cr
&&&&\multispan2\hrulefill&&\cr}}
\end{equation}
where
$\widetilde{T}_4=\{E_{\pm e_1\pm e_2}\}$
 is spanned over the generators corresponding to the cell (1, 2), and
$T_4=\{E_{\pm e_1},E_{\pm e_2}\}$
 -- over generators corresponding to the cells (1, 3) and (2, 3).

In general, algebras
$so\,(n+1;\bj)$, $n=2m$, and the table
$\Gamma_{n+1}(\bj)$
 can be constructed as follows :
to the cell
$(\Gamma_{n+1}(\bj))_{kk}$, $k=1,2,\dots,m$, we put in correspondence
generator $H_k$ from Cartan subalgebra; to the cell
$(\Gamma_{n+1}(\bj))_{k,m+1}$
 we put in correspondence generators
$E_{\pm e_k}$ and inscribe
in the cell
$\prod\limits_{l=2k}^n j_l^2$; to the cell
$(\Gamma_{n+1}(\bj))_{rs}$, $1\leq r<s\leq m$,
we put in correspondence generators
$E_{\pm e_r\pm e_s}$ and inscribe
$\prod\limits_{l=2r}^{2s-2} j_l^2$
in the cell. The structure of algebra
$so\,(n+1;\bj)$
 under contraction
over even parameters is described by the following theorem \cite{GMOrtG}.

\begin{theorem}  
\label{t1.2}
Let integer $k_r$ satisfy inequalities
$k_0=0<k_1<\dots<k_s<\dots<k_p<k_{p+1}=(n+1)/2$ if
$p<m=n/2$, and
parameters
$j_{2k_1}=\iota_{2k_1}$,
$j_{2k_2}=\iota_{2k_2}$, $\dots$, $j_{2k_p}=\iota_{2k_p}$.
Then the following expansion is valid:
\begin{equation}
\label{1.81}
so\,(n+1;\bj)=T\notni M;
\end{equation}
where $T$ is nilpotent radical, and $M$ is semisimple algebra
\begin{equation}
\label{1.82}
M=\sum_{s=0}^{p}\!\!{}^\oplus so\,(2(k_{s+1}-k_s);j_{2k_s+1},
\dots,j_{2k_{s+1}-1}).
\end{equation}

\end{theorem}

If
$k_{s+1}=k_s+1$, then
$so\,(2(k_{s+1}-k_s);j_{2k_s+1},\dots,j_{2k_{s+1}-1})=$\linebreak
$=so\,(2;j_{2k_s+1})\equiv H_{s+1}(j_{2k_s+1})$. For
$k_p=n/2$ we can wright as $so\,(2(k_{p+1}-k_p);j_{2k_p+1})=H_m(j_{n-1})$. The
expansion (\ref{1.81}) is Levi-Maltsev expansion of contracted algebra
$so\,(n+1;\bj)$.
In the case of maximal contraction, when all even parameters of
algebra
$so\,(n+1;\bj)$
 are equal to dual units, semi-simple algebra $M$ can be
expanded as follows:
$M=\sum\limits_{s=0}^{m}\!\!{}^\oplus h_s$.

Let us consider algebras
$so\,(n+1;\bj)$
 of even dimensions for $n=2m-1$, Cartan subalgebra $h$
consists of matrices $H$ of the type (\ref{1.64}) with omitted last row and
column. Basis matrices $H_k$ are the same as for algebra
$so\,(2m+1;\bj)$. Linear
functions on $h$ are defined, like in the odd-dimensional case, as
$e_k(H)=h_k$.
Simple roots are given in (\ref{1.63}). The system of all roots is as
\begin{equation}
\label{1.83}
\Delta=\{\pm e_r\pm e_s,\ r\neq s,r,s=1,2,\dots,m\}.
\end{equation}
 To the roots
$\pm e_r\pm e_s$
 there correspond generators
$E_{\pm e_r\pm e_s}$
 of the type
of (\ref{1.67}), (\ref{1.68})
 with commutation relations described by (\ref{1.69}) --
(\ref{1.71}). We shall discuss contractions of algebra
$so\,(2;\bj)$
 over parameters $j_k$
with even numbers, and for this reason we shall not write out
parameters $j_k$ with odd numbers when constructing the table
$\Gamma_{2m}(\bj)$.

The system of roots for the algebra
$so\,(4;\bj)$
 is
$\Delta=\{\pm e_1\pm e_2\}$, and commutator of generators
$E_{e_1-e_2}$ and $E_{-e_1+e_2}$
 is given by (\ref{1.72}). To the elements
$(\Gamma_4(\bj))_{kk}$
 we put in correspondence
generators
$H_k,\ k=1,2$, of Cartan subalgebra, to the element
$(\Gamma_4(\bj))_{12}$
 -- four generators
$E_{\pm e_1\pm e_2}$
and inscribe in this cell parameter $j_2^2$. As a result, we obtain
the table

\begin{equation}
\label{1.84}
\vbox{\offinterlineskip
\halign {#\on&#&\on #\on&#&\on #\on&#&\on #\cr
&&\strut\on $\ss{1}$&&\on $\ss{2}$&&\cr
&\multispan4\hrulefill&\cr
\phantom{a}&\vrule&&\vrule&&\vrule&\cr
$\strut\Gamma_4(\bj)=$&\vrule&$
H_1$&\vrule&$\,j^2_2$&\vrule&$\ss{1}$\cr
\phantom{a}&\vrule&&\vrule&&\vrule&\cr
&\multispan4\hrulefill&\cr
\phantom{a}&&&\vrule&&\vrule&\cr
&&&\vrule&$\strut H_2$&\vrule & $\ss{2}$\cr
\phantom{a}&&&\vrule&&\vrule&\cr
&&&\multispan2\hrulefill&\cr}}
\end{equation}
 For contraction
$j_2=\iota_2$
 we get
\begin{equation}
\label{1.85}
\vbox{\offinterlineskip
\halign {#\on&#&\on #\on&#&\on #\on&#&\on #&\on #\cr
&&\strut\on $\ss{1}$&&\on $\ss{2}$&&&\cr
&\multispan4\hrulefill&&\cr
\phantom{a}&\vrule&&\vrule&&\vrule&&\cr
$\strut\Gamma_4(\iota_2)=$&\vrule&$
H_1$&\vrule&$\,j^2_2$&\vrule&$\ss{1}$&,\quad
$so\,(4;\iota_2)=T_4\notni(H_1\oplus H_2),$\cr
\phantom{a}&\vrule&&\vrule&&\vrule&&\cr
&\multispan4\hrulefill&&\cr
\phantom{a}&&&\vrule&&\vrule&&\cr
&&$\oplus$&\vrule&$\strut H_2$&\vrule&$\ss{2}$&\cr
\phantom{a}&&&\vrule&&\vrule&&\cr
&&&\multispan2\hrulefill&\cr}}
\end{equation}
 where
$T_4=\{E_{\pm e_1\pm e_2}\}$
 is Abelian subalgebra.

For algebra
$so\,(6;\bj)$
 we get the system of roots
$\Delta=\{\pm e_1\pm e_2,\pm e_1\pm$\linebreak $\pm e_3,\pm e_2\pm
e_3\}$; commutators required for inserting in the cells of the table
$\Gamma_6(\bj)$, satisfy the following relations:
\bea
\label{1.86}
\left[E_{e_1-e_2},E_{-e_1+e_2}\right] 
&\!\!\!=\!\!\!& -4j_2^2(H_1+H_2),\quad\nonumber \\ 
\left[E_{e_2-e_3},E_{-e_2+e_3}\right] 
 &\!\!\!=\!\!\!& -4j_4^2(H_2+H_3),\\
\left[E_{e_1-e_3},E_{-e_1+e_3}\right] 
&\!\!\!=\!\!\!& -4j_2^2j_4^2(H_1+H_3).\nonumber   
\eea
To the element $(\Gamma_6(\bj))_{kk}$
 we put in correspondence generator $H_k,\quad k=$\linebreak$=1,2,3$;
to the element $(\Gamma_6(\bj))_{rs}$ -- four generators
$E_{\pm e_r\pm e_s}$ and inscribe in the cell parameter
$\prod\limits_{e=2r}^{2s-2}j_e^2,\quad r,s=1,2,3,\quad r<s$.  As a
result we get the table

\begin{equation}
\label{1.87}
\vbox{\offinterlineskip \halign {#\on&#&\on #\on&#&\on
#\on&#&\,#\,&#&\on #\cr
&&\strut\on$\ss{1}$&&\on$\ss{2}$&&\on$\ss{3}$&&\cr
&\multispan6\hrulefill&\cr
\phantom{a}&\vrule&&\vrule&&\vrule&&\vrule&\cr
$\strut\Gamma_6(\bj)=$&\vrule&$
H_1$&\vrule&$\,j^2_2$&\vrule&$j^2_2j^2_4$& \vrule&$\ss{1}$\cr
\phantom{a}&\vrule&&\vrule&&\vrule&&\vrule&\cr
&\multispan6\hrulefill&\cr
\phantom{a}&&&\vrule&&\vrule&&\vrule&\cr
&&&\vrule&$\strut H_2$&\vrule&$\on j^2_4$&\vrule&$\ss{2}$\cr
\phantom{a}&&&\vrule&&\vrule&&\vrule&\cr
&&&\multispan4\hrulefill&\cr
\phantom{a}&&&&&\vrule&&\vrule&\cr
&&&&&\vrule&$\strut\, H_3$&\vrule&$\ss{3}$\cr
&&&&&\multispan2\hrulefill&\cr}}\ ,
\end{equation}
which enables us easily to find the structure of contracted algebras
\bea
\label{1.88}
so\,(6;\iota_2) &\!\!\!=\!\!\!& T_8\notni so\,(4;j_3,j_4,j_5),
\nonumber \\
so\,(6;\iota_4) &\!\!\!=\!\!\!& \widetilde{T}_8\notni 
so\,(4;j_1,j_2,j_3),\\ 
 so\,(6;\iota_2,\iota_4) &\!\!\!=\!\!\!& T_{12}\notni (H_1\oplus 
H_2\oplus H_3),\nonumber  
\eea
where $T_8=\{E_{\pm e_1\pm 
e_2},E_{\pm e_1\pm e_3}\}$ and $\widetilde{T}_8=\{E_{\pm e_1\pm 
e_3},E_{\pm e_2\pm e_3}\}$ are commutative ideals; $T_{12}=\{E_{\pm 
e_1\pm e_2},E_{\pm e_1\pm e_3},E_{\pm e_2\pm e_3}\}$ is nilpotent 
 radical.

In general case of algebra
$so\,(2m,\bj)$
 the table
$\Gamma_{2m}(\bj)$
  can be constructed in the same
way as for
$so\,(6,\bj)$, but with
$r,s=1,2,\dots,m$, $r<s$, and the structure of contracted algebras is
described by the following theorem.

\begin{theorem} 
\label{t1.3}
Let integer $k_r$ satisfy inequalities
$k_0=0<k_1<$\linebreak$<\dots<k_s<\dots<k_p<k_{p+1}=m$, and
parameters
$j_{2k_1}=\iota_{2k_1},\dots,j_{2k_s}=$\linebreak 
 $=\iota_{2k_s},\dots,j_{2k_p}= \iota_{2k_p}$. Then Levi-Maltsev 
expansion is valid: 
\begin{equation} 
\label{1.89}
so\,(2m;\bj)=T\notni M, 
\end{equation} 
where $T$ is nilpotent radical; $M$ is semisimple algebra of
the type 
\begin{equation}
\label{1.90}
M=\sum_{s=0}^{p}\!\!{}^\oplus so\,(2(k_{s+1}-k_s);j_{2k_s+1},\dots,
j_{2k_{s+1}-1}),
\end{equation}
and
$so\,(2;j_{2k_s+1})\equiv H_{s+1}(j_{2k_s+1})$.
\end{theorem}

\subsection{Structure of unitary Cayley-Klein algebras}
\medskip

 The structure of unitary algebras under contractions we shall also
describe, using the root technique. In contrast to orthogonal algebras, 
this method in the case of algebras $u\,(n+1;\bj)$, 
$su\,(n+1;\bj)$ enables us to consider the contractions over all 
parameters $\bj$, because Cartan subalgebra consists of the 
diagonal matrices which do not depend on $\bj$.

A matrix of the Cartan subalgebra $h$ can
be written as
\begin{equation}
\label{1.125}
H=\left( \ba{cccc}
h_0&&&0\\
&h_1&&\\
&&\ddots&\\
0&&&h_n.
\ea \right)
\end{equation}
As a basis in $h$, we choose matrices
$H_k=A_{kk}$, $k=0,1,\dots,n$. The linear
functions on Cartan subalgebra are defined by relation $e_k(H)=h_k$,
where $H$ is matrix (\ref{1.125}). The simple roots are given by 
(\ref{1.60}). The
system of all roots is
$\Delta=\{e_k-e_m,\ k\neq m,\ k,m=0,1,\dots,n\}$.
To the root $e_k-e_m$ there corresponds generators
$E_{e_k-e_m}=A_{km}(\bj)$
\bea
\label{1.126}
E_{e_\mu-e_\nu} 
&\!\!\!=\!\!\!& A^*_{\mu\nu}\prod^{\nu}_{l=\mu+1}j_l^2,\quad 
E_{e_\nu-e_\mu}=A^*_{\nu\mu},\quad \mu<\nu,\nonumber \\
\mu,\nu &\!\!\!=\!\!\!& 0,1,\dots,n.
\eea
 Weyl-Cartan commutation relations (\ref{1.58}) can be found, using 
(\ref{1.97}), and are as follows (we have written only non-zero commutators):  
\beq
\left[H_k,E_{e_k-e_p}\right] &\!\!\!=\!\!\!& E_{e_k-e_p},\quad
\left[H_k,E_{-e_k+e_p}\right]=-E_{-e_k+e_p},\\
\left[E_{e_k-e_m},E_{-e_k+e_m}\right] &\!\!\!=\!\!\!& 
\biggl(\prod_{l=k+1}^{m}j_l^2\biggr)
(H_k-H_m),
\eeq
\bea
\label{1.127}
&&\left[E_{e_k-e_m},E_{e_m-e_q}\right]=\nonumber \\
 &\!\!\!=\!\!\!& \left\{ \ba{cccc}
E_{e_k-e_q},& k<m<q &  {\rm or}& q<m<k,\\
\biggl(\prod\limits^{m}_{l=1+\max(k,q)}j_l^2\biggr)E_{e_k-e_q},&
m>k,&& m>q,\\
\biggl(\prod\limits^{\min(k,q)}_{l=m+1}j_l^2\biggr)E_{e_k-e_q},&
m<k, && m<q.\ea \right.
\eea

The structure of contracted algebras $u\,(n+1;\bj)$ is
completely determined by the commutators (\ref{1.127}) and described by the
auxiliary table $\Gamma_{n+1}(\bj)$ which can be compiled as
follows. To the diagonal elements of the table there correspond basis
generators of the Cartan subalgebra
$(\Gamma_{n+1}(\bj))_{kk}=H_k$, $k=0,1,\dots,n$. To the elements
$(\Gamma_{n+1}(\bj))_{\mu\nu}$, $\mu<\nu$,
$\mu,\nu=0,1,\dots,n$, of the
table there correspond two generators
$E_{e_\mu-e_\nu}$, $E_{e_\nu-e_\mu}$, and we inscribe in the cell
the product
$\prod\limits^{\nu}_{e=\mu+1}j_e^2$,
arising under commutation of these two generators.

To algebra $u\,(2;j_1)$ there corresponds the table
\begin{equation}
\label{1.128}
\vbox{\offinterlineskip
\halign {#\on&#&\on #\on&#&\on #\on&#&\on #\cr
&&\strut\on$\ss{0}$&&\on$\ss{1}$&&\cr
&\multispan4\hrulefill&\cr
\phantom{a}&\vrule&&\vrule&&\vrule&\cr
$\strut\Gamma_2(j_1)=$&\vrule&$
H_0$&\vrule&$\,j^2_1$&\vrule&${\ss{0}}\, ,$\cr
\phantom{a}&\vrule&&\vrule&&\vrule&\cr
&\multispan4\hrulefill&\cr
\phantom{a}&&&\vrule&&\vrule&\cr
&&&\vrule&$\strut H_1$&\vrule&$\ss{1}$\cr
\phantom{a}&&&\vrule&&\vrule&\cr
&&&\multispan2\hrulefill&\cr}}
\qquad
\vbox{\offinterlineskip
\halign {#\on&#&\on #\on&#&\on #\on&#&\on #&\on #\cr
&&&&&&&\cr
&\multispan4\hrulefill&&\cr
\phantom{a}&\vrule&&\vrule&&\vrule&&\cr
$\strut\Gamma_2(\iota_1)=$&\vrule&$
H_0$&\vrule&$\on 0$&\vrule&&\cr
\phantom{a}&\vrule&&\vrule&&\vrule&&\cr
&\multispan4\hrulefill&&\cr
\phantom{a}&&&\vrule&&\vrule&&\cr
&&$\oplus$&\vrule&$\strut H_1$&\vrule&&\cr
\phantom{a}&&&\vrule&&\vrule&&\cr
&&&\multispan2\hrulefill&\cr}}\;.
\end{equation}
For contraction $j_1=\iota_1$ the table is the same as (\ref{1.128}), and
the structure of contracted algebra is as follows:

\begin{equation}
\label{1.129}
u\,(2;\iota_1)=T_2\notni(H_1\oplus H_2),\quad T_2=\{E_{\pm e_1\pm
e_2}\},
\end{equation}
where it is
necessary to take either the upper, or the lower signs, i.e.  either
$(+,-)$, or $(-,+)$. To algebra $u\,(3;j_1,j_2)$ and to algebras
obtained from the latter using all possible contractions, there
correspond the tables
\beq
\vbox{\offinterlineskip \halign {#\on&#&\on #\on&#&\on
#\on&#&\,#\,&#&\on #\cr
&&&&&&&&\cr
&\multispan6\hrulefill&\cr
\phantom{a}&\vrule&&\vrule&&\vrule&&\vrule&\cr
$\strut\Gamma_3(\bj)=$&\vrule&$
H_0$&\vrule&$\,j^2_1$&\vrule&$j^2_1j^2_2$& \vrule&\cr
\phantom{a}&\vrule&&\vrule&&\vrule&&\vrule&\cr
&\multispan6\hrulefill&\cr
\phantom{a}&&&\vrule&&\vrule&&\vrule&\cr
&&&\vrule&$\strut H_1$&\vrule&$\on j^2_2$&\vrule&\cr
\phantom{a}&&&\vrule&&\vrule&&\vrule&\cr
&&&\multispan4\hrulefill&\cr
\phantom{a}&&&&&\vrule&&\vrule&\cr
&&&&&\vrule&$\strut\, H_{2_{\phantom{9}}}$&\vrule&\cr
&&&&&\multispan2\hrulefill&\cr}}\ ,
\quad
\vbox{\offinterlineskip \halign {#\on&#&\on #\on&#&\on
#\on&#&\,#\,&#&\on #\cr
&&&&&&&&\cr
&\multispan6\hrulefill&\cr
\phantom{a}&\vrule&&\vrule&&\vrule&&\vrule&\cr
$\strut\Gamma_3(\iota_1,j_2)=$&\vrule&$
H_0$&\vrule&$\on 0$&\vrule&$\on 0\on$&\vrule&\cr
\phantom{a}&\vrule&&\vrule&&\vrule&&\vrule&\cr
&\multispan6\hrulefill&\cr
\phantom{a}&&&\vrule&&\vrule&&\vrule&\cr
&&$\oplus$&\vrule&$\strut H_1$&\vrule&$\on
j^2_2$&\vrule&\cr
\phantom{a}&&&\vrule&&\vrule&&\vrule&\cr &&&\multispan4\hrulefill&\cr
\phantom{a}&&&&&\vrule&&\vrule&\cr
&&&&&\vrule&$\strut\, H_{2_{\phantom{9}}}$&\vrule&\cr
&&&&&\multispan2\hrulefill&\cr}}\ ,
\eeq
\begin{equation}
\label{1.130}
\end{equation}
\beq
\vbox{\offinterlineskip \halign {#\on&#&\on #\on&#&\on
#\on&#&\,#\,&#&\on #\cr
&&&&&&&&\cr
&\multispan6\hrulefill&\cr
\phantom{a}&\vrule&&\vrule&&\vrule&&\vrule&\cr
$\strut\Gamma_3(j_1,\iota_2)=$&\vrule&$
H_0$&\vrule&$\,j_1^2$&\vrule&$\on 0\phantom{0}$& \vrule&\cr
\phantom{a}&\vrule&&\vrule&&\vrule&&\vrule&\cr
&\multispan6\hrulefill&\cr
\phantom{a}&&&\vrule&&\vrule&&\vrule&\cr
&&&\vrule&$\strut H_1$&\vrule&$\on 0$&\vrule&\cr
\phantom{a}&&&\vrule&&\vrule&&\vrule&\cr
&&&\multispan4\hrulefill&\cr
\phantom{a}&&&&&\vrule&&\vrule&\cr
&&&&$\oplus$&\vrule&$\strut\, H_{2_{\phantom{9}}}$&\vrule&\cr
&&&&&\multispan2\hrulefill&\cr}}\ ,
\quad
\vbox{\offinterlineskip \halign {#\on&#&\on #\on&#&\on
#\on&#&\,#\,&#&\on #\cr
&&&&&&&&\cr
&\multispan6\hrulefill&\cr
\phantom{a}&\vrule&&\vrule&&\vrule&&\vrule&\cr
$\strut\Gamma_3(\iota_1,\iota_2)=$&\vrule&$
H_0$&\vrule&$\on 0$&\vrule&$\on 0\phantom{0}$& \vrule&\cr
\phantom{a}&\vrule&&\vrule&&\vrule&&\vrule&\cr
&\multispan6\hrulefill&\cr
\phantom{a}&&&\vrule&&\vrule&&\vrule&\cr
&&$\oplus$&\vrule&$\strut H_1$&\vrule&$\on 0$&\vrule&\cr
\phantom{a}&&&\vrule&&\vrule&&\vrule&\cr
&&&\multispan4\hrulefill&\cr
\phantom{a}&&&&&\vrule&&\vrule&\cr
&&&&$\oplus$&\vrule&$\strut\, H_{2_{\phantom{9}}}$&\vrule&\cr
&&&&&\multispan2\hrulefill&\cr}}\ ,
\eeq
from which we can find the structure of
contracted algebras
\bea
\label{1.131}
u\,(3;\iota_1,j_2) &\!\!\!=\!\!\!& T_4\notni (H_0\oplus u_2(j_2)),
\quad T_4=\{E_{\pm e_1\mp e_2},E_{\pm e_1\mp e_3}\},\nonumber \\
u\,(3;j_1,\iota_2) &\!\!\!=\!\!\!& 
\widetilde{T}_4\notni (u_2(j_1)\oplus H_2),
\quad \widetilde{T}_4=\{E_{\pm e_1\mp e_3},E_{\pm e_2\mp e_3}\},
\nonumber \\
u\,(3;\iota_1) &\!\!\!=\!\!\!& T_6\notni (H_0\oplus H_1\oplus H_2),
\quad T_6=T_4\cup\widetilde{T}_4.
\eea

If only one parameter $j_{k_1}=\iota_{k_1}$ of algebra
$u\,(n+1;\bj)$ takes a dual value (i.e. under one-dimensional
contraction), then its structure is as follows:
\bea
\label{1.132}
u\,(n &\!\!\!+\!\!\!& 1;\dots,\iota_{k_1},\dots)=T_{k_1}\notni
(u(k_1;j_1,\dots,j_{k_1-1})\oplus\nonumber \\
 &\!\!\!\oplus\!\!\!& 
 u\,(n-k_1+1;j_{k_1+1,\dots,j_n})),
\eea
where Abelian ideal
$T_{k_1}$ is spanned over generators, corresponding to the zero cells
of the table $\Gamma_{n+1}(\bj)$
\bea
\label{1.133}
T_{k_1} &\!\!\!=\!\!\!& \left\{E_{\pm e_\mu\mp e_\nu},\quad 
\mu=0,1,\dots,k_1-1,  \right. \nonumber \\ 
\nu &\!\!\!=\!\!\!&  \left. k_1, \quad k_1+1,\dots,n\right\}, 
\eea
and its dimension is $2k_1(n-k_1+1)$.  Under two-dimensional 
contraction $j_{k_1}=\iota_{k_1}$, $j_{k_2}=\iota_{k_2}$, $k_1<k_2$ 
ideal $T_2$ can be obtained as a
union of ideals $T_{k_1}$ and $T_{k_2}$, i.e.
\bea
\label{1.134}
T_2 &\!\!\!=\!\!\!& T_{k_1}\cup T_{k_2}=\left\{E_{\pm e_\mu\mp
e_\nu},\quad\mu=0,1,\dots,k_1-1, \right.\nonumber \\
\nu &\!\!\!=\!\!\!& k_1,\, k_1+1,\dots,n\, \mu=k_1,k_1+1,\dots,k_2-1,
\nonumber \\
\nu &\!\!\!=\!\!\!& \left.  k_2,k_2+1,\dots,n\right\}
\eea
and is nilpotent. Under multidimensional contraction
$j_{k_1}=\iota_{k_1},\dots,j_{k_m}=$\linebreak$=\iota_{k_m}$
nilpotent $T$
is as follows
\bea
\label{1.135}
T &\!\!\!=\!\!\!& 
 \bcup\limits^{m}_{p=1} T_{k_p}=\left\{E_{\pm e_\mu\mp
e_\nu},\quad\mu= k_{p-1},k_{p-1}+1,\dots,k_p-1,\right. \nonumber \\
\nu &\!\!\!=\!\!\!& \left.  k_p,k_p+1,\dots,n,\quad k_0=0,\quad 
p=1,2,\dots,m\right\}.  \eea 

As in the case of orthogonal algebras, the structure of nilpotent ideal $T$
under multidimensional contractions is not fixed uniquely and determined by
the partitions of zero parts of the table $\Gamma(\bj)$ into
blocks, arising under different sequences of one-dimensional
contractions, which lead to a given multidimensional contraction. In
general the structure of contracted unitary algebra is described by
following theorem \cite{GMUniG}.

\begin{theorem} 
\label{t1.5}
Let integer $k_p$ satisfy inequalities
$k_0=0<k_1<k_2<$\linebreak $<\dots<k_m<k_{m+1}=n+1$, and let
$j_{k_1}=\iota_{k_1},\dots,j_{k_m}=\iota_{k_m}$.
Then the following Levi-Maltsev expansion takes place:
\begin{equation}
\label{1.136}
u\,(n+1;\bj)=T\notni M,
\end{equation}
where
$T$ is nilpotent radical (\ref{1.135}); $M$ is semisimple algebra
\begin{equation}
\label{1.137}
M=\sum_{p=0}^{m}\!\!{}^\oplus u\,(k_{p+1}-k_p;\,
j_{k_p+1},\dots,j_{k_{p+1}-1}),
\end{equation}
where  for
$k_{p+1}=k_p+1$   algebra
$u\,(k_{p+1}-k_p)=u\,(1)=H_{k_p}$.
\end{theorem}

In the case of maximal contraction, when all $j_k=\iota_k$, algebra
$M$ coincides with Cartan subalgebra
$h=\sum\limits_{p=0}^{m}\!\!{}^\oplus H_p=M$.

The matrices of special unitary algebra $su\,(n+1;\bj)$ have
zero trace; for this reason its Cartan subalgebra $h$ consists of
matrices $\widetilde{H}$ of the type (\ref{1.125})  with additional 
condition $\sum\limits_{k=0}^{n}h_k=0$.  As a basis in $h$ we choose 
matrices $\widetilde{H}_k=A_{k-1,k-1}-A_{kk}$, $k=1,2,\dots,n$, and the 
linear functions on $h$ are $e_k(\widetilde{H})=h_k$. Algebra 
$su\,(n+1;\bj)$ has the same system of roots $\Delta$ and generators 
are put in correspondence to the roots as in the case of algebra 
$u\,(n+1;\bj)$, but the commutators, involving generators 
$\widetilde{H}_k$ are somewhat different:  
\bea 
\label{1.138}
\left[\widetilde{H}_k,E_{e_k-e_p}\right] &\!\!\!=\!\!\!& 
-E_{e_k-e_p},\quad p<k-1,\nonumber \\
\left[\widetilde{H}_{k+1},E_{e_k-e_p}\right] &\!\!\!=\!\!\!& 
E_{e_{n_k}-e_p},\nonumber \\
\left[\widetilde{H}_{k},E_{e_k-e_{k-1}}\right] &\!\!\!=\!\!\!& 
-2E_{e_k-e_{k-1}},\\
\left[E_{e_k-e_m},E_{-e_k+e_m}\right] &\!\!\!=\!\!\!& 
\biggl(\prod^{m}_{l=k+1}j_l^2\biggr)\sum_{r=k+1}^{m}\widetilde{H}_r,
\quad k<m. \nonumber 
\eea
The table
$\widetilde{\Gamma}_{n+1}(\bj)$
differs from the table
$\Gamma_{n+1}(\bj)$ for unitary
algebras \linebreak in that to diagonal elements there correspond
generators $\widetilde{H}_k$, i.e. \linebreak
$(\widetilde{\Gamma}_{n+1}(\bj))_{kk}=\widetilde{H}_k$,
$k=1,2,\dots,n$. To algebra $su\,(2;j_1)$ there corresponds the table
$\widetilde{\Gamma}_2(j_1)=\frf{\,\widetilde{H}_1\bigm|j^2_1\,}$\,,
and to contracted algebra $su\,(2;\iota_1)=T_2\notni \widetilde{H}_0$,
where $T_2$ is given by (\ref{1.129}), -- the table
$\widetilde{\Gamma}_2(\iota_1)=\frf{\,\widetilde{H}_1\bigm|0\,}\,$.
To algebra $su\,(3;j_1,j_2)$ and
algebras, derived from it by contractions, there correspond the tables
\beq
\vbox{\offinterlineskip
\halign {#\on&#&\on #\on&#&\on #\on&#&\,#\,&#&\on #\cr
&&\strut\on$\ss{0}$&&\on$\ss{1}$&&\on$\ss{2}$&&\cr
&\multispan6\hrulefill&\cr
\phantom{a}&\vrule&&\vrule&&\vrule&&\vrule&\cr
$\strut\widetilde{\Gamma}_3(j)=$&\vrule&$
\widetilde{H}_1$&\vrule&$\,j^2_1$&\vrule&$j^2_1j^2_2$&
\vrule&$\ss{0}$\cr
\phantom{a}&\vrule&&\vrule&&\vrule&&\vrule&\cr
&\multispan6\hrulefill&\cr
\phantom{a}&&&\vrule&&\vrule&&\vrule&\cr
&&&\vrule&$\strut \widetilde{H}_2$&\vrule&$\on
j^2_2$&\vrule&$\ss{1}$\cr
\phantom{a}&&&\vrule&&\vrule&&\vrule&\cr
&&&\multispan4\hrulefill&\cr}}\ ,
\quad
\vbox{\offinterlineskip
\halign {#\on&#&\on #\on&#&\on #\on&#&\on #\on&#\cr
&&&&&&&\cr
&\multispan6\hrulefill\cr
\phantom{a}&\vrule&&\vrule&&\vrule&&\vrule\cr
$\strut\widetilde{\Gamma}_3(\iota_1,j_2)=$&\vrule&$
\widetilde{H}_1$&\vrule&$\on 0$&\vrule&$\, 0\,$& \vrule\cr
\phantom{a}&\vrule&&\vrule&&\vrule&&\vrule\cr
&\multispan6\hrulefill\cr
\phantom{a}&&&\vrule&&\vrule&&\vrule\cr
&&$\oplus$&\vrule&$\strut \widetilde{H}_2$&\vrule&$\,
j^2_2\,$&\vrule\cr
\phantom{a}&&&\vrule&&\vrule&&\vrule\cr
&&&\multispan4\hrulefill\cr}}\ ,
\eeq
\begin{equation}
\label{1.139}
\end{equation}
\beq
\vbox{\offinterlineskip
\halign {#\on&#&\on #\on&#&\on #\on&#&\on #\on&#\cr
&\multispan6\hrulefill\cr
\phantom{a}&\vrule&&\vrule&&\vrule&&\vrule\cr
$\strut\widetilde{\Gamma}_3(j_1,\iota_2)=$&\vrule&$
\widetilde{H}_1$&\vrule&$\on j_1^2$&\vrule&$\, 0\,$& \vrule\cr
\phantom{a}&\vrule&&\vrule&&\vrule&&\vrule\cr
&\multispan6\hrulefill\cr
\phantom{a}&&&\vrule&&\vrule&&\vrule\cr
&&&\vrule&$\strut \widetilde{H}_2$&\vrule&$\,
0\,$&\vrule\cr
\phantom{a}&&&\vrule&&\vrule&&\vrule\cr
&&&\multispan4\hrulefill\cr}}\ ,
\quad
\vbox{\offinterlineskip
\halign {#\on&#&\on #\on&#&\on #\on&#&\on #\on&#\cr
&\multispan6\hrulefill\cr
\phantom{a}&\vrule&&\vrule&&\vrule&&\vrule\cr
$\strut\widetilde{\Gamma}_3(\iota_1,\iota_2)=$&\vrule&$
\widetilde{H}_1$&\vrule&$\on 0$&\vrule&$\, 0\,$& \vrule\cr
\phantom{a}&\vrule&&\vrule&&\vrule&&\vrule\cr
&\multispan6\hrulefill\cr
\phantom{a}&&&\vrule&&\vrule&&\vrule\cr
&&$\oplus$&\vrule&$\strut \widetilde{H}_2$&\vrule&$\,
0\,$&\vrule\cr
\phantom{a}&&&\vrule&&\vrule&&\vrule\cr
&&&\multispan4\hrulefill\cr}}\ ,
\eeq 
using which it is possible to find the structure of algebras:
\bea
\label{1.140}
su\,(3;\iota_1,j_2) &\!\!\!=\!\!\!& 
T_4\notni (\widetilde{H}_0\oplus su\,(2;j_2))=
T_4\notni u\,(2;j_2),\nonumber \\
su\,(3;j_1,\iota_2) &\!\!\!=\!\!\!& 
\widetilde{T}_4\notni u\,(2;j_1),\nonumber \\
su\,(3;\iota) &\!\!\!=\!\!\!& T_6\notni(\widetilde{H}_1\oplus 
\widetilde{H}_2), 
\eea 
where ideals $T$ are the 
same as in (\ref{1.131}).

The first of the formulas (\ref{1.140}) gives the same structure of algebra
$su\,(3;\iota_1,j_2)$
which is obtained from the commutators (see the formulas (23), page 548
in [5]) for $j_1=\iota_1$, when
taking into account the relation between generators
\bea
\label{1.141}
E_{e_\mu-e_\nu} &\!\!\!=\!\!\!& -{1\over2}(L_{\mu\nu}+iQ_{\mu\nu}),
\nonumber \\
E_{e_\nu-e_\mu} &\!\!\!=\!\!\!& {1\over2}(L_{\mu\nu}-iQ_{\mu\nu}),\\
\widetilde{H}_k &\!\!\!=\!\!\!& -iP_k,\quad k=1,2,\dots,n.
\eea 

In general the structure of contracted algebras $su\,(n+1;\bj)$
is described by the following theorem \cite{GMUniG}.

\begin{theorem} Let integer $k_p$ satisfy the inequalities
$k_0= 0<k_1<$\linebreak$<k_2<\dots<k_m<n+1$, and parameters take the
values $j_{k_1}=$\linebreak $=\iota_{k_1},\dots,j_{k_m}=\iota_{k_m}$.
Then the following Levi-Maltsev expansion is valid:  
\begin{equation}
\label{1.142}
su\,(n+1;\bj)=T\notni M,
\end{equation}
where $T$ is nilpotent
radical (\ref{1.135}); $M$ is semisimple algebra 
\bea
\label{1.143}
M &\!\!\!=\!\!\!& \sum_{p=0}^{m-1}\!\!{}^\oplus
u\,(k_{p+1}-k_p;j_{k_p+1},\dots,j_{k_{p+1}-1})\oplus
\nonumber \\ 
 &\!\!\!\oplus\!\!\!& 
su\,(n-k_m+1;j_{k_m+1},\dots,j_n), 
\eea
where for $k_{p+1}=k_p+1$ algebra
$u\,(k_{p+1}-k_p)=u\,(1)=\widetilde{H}_{k_p+1}$, and for $k_m=n$
$su\,(1)=\widetilde{H}_n$.
\end{theorem}

In the case of maximal contraction $\bj=\biota$
algebra $M$ coincides with Cartan subalgebra
$h=\sum\limits_{p=1}^{n}\!\!{}^\oplus \widetilde{H}_p=M$.

\subsection{Structure of symplectic Cayley-Klein algebras}
\medskip

For symplectic algebras $sp\,(n;\bj)$ Cartan subalgebra $h$
consists of the diagonal matrices
\begin{equation}
\label{1.158}
H=\left[ \ba{cccccc}
h_1&&&&&\\
&\ddots&&&&0\\
&&h_n&&&\\
&&&-h_1&&\\
0&&&&\ddots&\\
&&&&&-h_n
\ea \right],
\end{equation}
which do not depend on parameters
$\bj$. Therefore, the root technique makes possible to consider
contractions over all parameters $\bj$ too. The linear functions on
Cartan subalgebra are as follows: $e_k(H)=h_k$, $k=1,2,\dots,n$. As a
basis in $h$ one chooses matrices $H_k=A_{kk}-A_{k+n,k+n}$. The simple
roots are given by (\ref{1.62}). The system of all roots is given by
$\Delta=\{\pm e_k\pm  \pm e_m\;,\quad k\neq m;\quad
k,m=1,2\dots,n\}\cup \{\pm2e_k\}$, where the signs plus and minus are
chosen independently.  To the roots there correspond generators
\bea
\label{1.159}
E_{e_k-e_m} &\!\!\!=\!\!\!& 
\biggl(\prod_{l=k+1}^{m} j_l^2\biggr)A_{km}-A_{m+n,k+n},
\quad E_{2e_k}=A_{k,k+n},\nonumber \\
E_{e_m-e_k} &\!\!\!=\!\!\!& 
A_{mk}-\biggl(\prod_{l=k+1}^{m} j_l^2\biggr)A_{k+n,m+n},
\quad E_{-2e_k}=A_{k+n,k},\nonumber \\
E_{-e_k-e_m} &\!\!\!=\!\!\!& \biggl(\prod_{l=k+1}^{m}
j_l^2\biggr)A_{k+n,m}+A_{m+n,k},\nonumber \\
E_{e_k+e_m} &\!\!\!=\!\!\!& \biggl(\prod_{l=k+1}^{m}
j_l^2\biggr)A_{k,m+n}+A_{m,k+n},\\
k<m,  &\!\!\!\!\!\!&  k,m=1,2,\dots,n. \nonumber 
\eea
 Cartan-Weyl commutators (\ref{1.58}) can be found using (\ref{2.18})
 and can be presented as follows (only non-zero 
commutators are written out):  
\bea
\label{1.160}
\left [H_k,E_{\pm e_m\pm e_p}\right] &\!\!\!=\!\!\!& 
\left\{ \ba{ll} \pm E_{\pm e_k\pm 
e_p},& k=m,\\ 
\pm E_{\pm e_m\pm e_k},&k=p,\\ \ea \right. 
\nonumber \\ 
\left[H_k,E_{\pm2e_k}\right] &\!\!\!=\!\!\!& 
\pm2E_{\pm2e_k},\nonumber \\ 
\left[E_{e_k\pm e_m},E_{-e_k\mp 
e_m}\right] &\!\!\!=\!\!\!& 
\biggl(\prod_{l=1+\min(k,m)}^{\max(k,m)} j_l^2\biggr) 
(H_k\pm H_m),\nonumber \\ 
\left[E_{2e_k},E_{-2e_k}\right] &\!\!\!=\!\!\!& H_k,\\
\left[E_{e_k\pm e_m},E_{-e_k\pm e_m}\right]  &\!\!\!=\!\!\!& 
-2\biggl(\prod_{l=1+\min(k,m)}^{\max(k,m)} j_l^2\biggr) 
E_{\pm2e_m},\nonumber \\
\left[E_{-e_k\pm e_m},E_{2e_k}\right] &\!\!\!=\!\!\!& 
\pm E_{e_k\pm e_m}, \quad
\left[E_{e_k\pm e_m},E_{-2e_k}\right]=\pm E_{-e_k\pm e_m},
\nonumber \\
\left[E_{e_k\pm e_m},E_{-e_k\pm e_p}\right] 
&\!\!\!=\!\!\!& 
\left\{ \ba{llll}
-\biggl(\prod\limits_{l=k+1}^{\min(m,p)}j_l^2\biggr)E_{\pm 
e_m\pm e_p},\;\;  k<m,\;\;  k<p,\\ 
-\biggl(\prod\limits_{l=1+\max(m,p)}^{k}j_l^2\biggr)E_{\pm e_m\pm
e_p},\;\; k>m,\;\; k>p,\\ 
-E_{\pm e_m\pm e_p},\;\; m<k<p \;\; {\rm or}\;\;
p<k<m. \\ \ea \right. \nonumber 
\eea
The signs $\pm$ of the generator in the first commutator are chosen as
in $\pm e_m$ for $k=m$ or in $\pm e_p$ for $k=p$. In the last
commutator the sign minus of the combination $+e_m,-e_p$ is
changed for the sign plus.

The commutation  relations (\ref{1.160}) completely determine the 
structure of symplectic algebras $sp\,(n;\bj)$ under contractions, 
which can be conveniently described, using the table $\Gamma_n(\bj)$. 
To the diagonal elements of the table $(\Gamma_n(\bj))_{kk}$ we put in 
correspondence the subalgebras $M_k=\{H_k,E_{\pm2e_k}\}$, 
$k=1,2,\dots,n$, which are isomorphic to the algebra $sp\,(1)$ and do 
not depend on the parameters $\bj$, because both generators $H_k,\quad 
E_{\pm2e_k}$ and their commutators $[H_k,E_{\pm2e_k}]=\pm2E_{\pm2e_k}$, 
$[E_{2e_k}$, $E_{-2e_k}]=H_k$ do not involve the parameters
$\bj$. To each element $(\Gamma_n(\bj))_{km}$, $k<m$,
$k,m=1,2,\dots,n$ we put in correspondence four generators $E_{\pm
e_k\pm e_m}$ and inscribe in this cell the product
$\prod\limits_{l=k+1}^{m}j_l^2$.

For the one-dimensional contraction $j_{k_1}=\iota_{k_1}$, $2\leq
k_1\leq n$ (the other parameters $j_k=1,i$) we get
\bea
\label{1.161}
sp\,(n; &\!\!\!\dots,\!\!\!& 
\iota_{k_1},\dots) =
T_{k_1}\notni(sp\,(k_1-1;j_2,\dots,j_{k_1-1})
\oplus\nonumber \\
 &\!\!\!\oplus \!\!\!& 
sp\,(n-k_1+1;j_{k_1+1},\dots,j_n)),
\eea
where Abelian
ideal $T_{k_1}$ is spanned over generators $E_{\pm e_k\pm e_m}$,
corresponding to zero cells of the table $\Gamma_n(\bj)$,
and subalgebras $sp\,(k_1-1)$ and $sp\,(n-k_1+1)$ are spanned over
generators, corresponding to non-zero cells i.e.
\bea
\label{1.162}
T_{k_1} &\!\!\!=\!\!\!& 
\{E_{\pm e_k\pm e_m},\quad1\leq k<k_1,\quad k_1\leq m\leq
n\},\nonumber \\
sp\,(k_1-1) &\!\!\!=\!\!\!& \left\{M_r,\quad1\leq r<k_1, \right.
\nonumber \\
E_{\pm e_k\pm e_m}, &\!\!\!\!\!\!& \left.
1\leq k<k_1-1,\quad 1<m<k_1,\quad k<m\right\},\nonumber \\
sp\,(n-k_1+1) &\!\!\!=\!\!\!& \left\{M_r,\quad k_1\leq r\leq n,
\right.\nonumber \\
E_{\pm e_p\pm e_s}, &\!\!\!\!\!\!&
\left.  k_1\leq p<n,\quad k_1<s\leq n,\quad
p<s\right\}.
\eea

For two-dimensional contraction $j_{k_1}=\iota_{k_1}$,
$j_{k_2}=\iota_{k_2}$, $2\leq k_1<k_2\leq n$ of symplectic algebra we
get
\bea
\label{1.163}
sp\,(n;\iota_{k_1},\iota_{k_2}) &\!\!\!=\!\!\!& 
T_2\notni(sp\,(k_1-1;j_2,\dots,j_{k_1-1})
\oplus\nonumber \\
 &\!\!\!\oplus\!\!\!& 
sp\,(k_2-k_1;j_{k_1+1},\dots,j_{k_2-1})\oplus\nonumber \\
 &\!\!\!\oplus\!\!\!& 
sp\,(n-k_2+1;j_{k_2+1},\dots,j_n)),
\eea
where $T_2=T_{k_1}\cup T_{k_2}$ is nilpotent ideal,
\bea
\label{1.164}
T_2 &\!\!\!=\!\!\!& 
\left\{E_{\pm e_k\pm e_m},\quad(1\leq k<k_1,\quad k_1\leq m\leq n),
\right. \nonumber \\
 &\!\!\!\!\!\!& \left.
(k_1\leq k<k_2,\quad k_2\leq m\leq n)\right\},
\eea
subalgebra $sp\,(k_1-1)$
and $sp\,(n-k_2+1)$ are described in (\ref{1.162}) (in the latter case $k_1$
should be changed for $k_2$) and rest subalgebra is spanned over the
following generators:
\bea
\label{1.165}
sp\,(k_2-k_1) &\!\!\!=\!\!\!& \left\{M_r,\quad k_1\leq r<k_2,\right.
\nonumber \\
E_{\pm e_p\pm e_s},\quad k_1\leq p<k_2-1,
 &\!\!\!\!\!\!& \left.
 k_1<s<k_2,\quad p<s\right. \}.
\eea

 In general case of multidimensional contraction the structure of 
symplectic algebras is described by the following theorem \cite{GMSymG}:
\begin{theorem} Let parameters be $j_{k_1}=\iota_{k_1},\
\dots$, $j_{k_m}=\iota_{k_m}$, where the integers $k_p$ satisfy the
inequalities $k_0=1<k_1<\dots< k_m<k_{m+1}=n+1$.
Then the following Levi-Maltsev expansion is valid:  
\begin{equation}
\label{1.166}
sp\,(n;\bj)=T\notni M, 
\end{equation}
where $T$ is
nilpotent radical, 
\bea
\label{1.167}
T &\!\!\!=\!\!\!& \bigcup\limits_{p=1}^m T_{k_p}=\left\{E_{\pm e_\mu\pm
e_\nu},\quad k_{p-1}\leq\mu<k_p\leq\nu\leq n,\quad\mu<\nu,\right.
\nonumber \\
p &\!\!\!=\!\!\!& 1,2,\dots,m\left. \right\}, 
\eea
and $M$ is semisimple algebra
\begin{equation}
\label{1.168}
M=\sum_{p=0}^{m}\!\!{}^\oplus sp\,(k_{p+1}-k_p; j_{k_p+1},\dots,
j_{k_{p+1}-1}), 
\end{equation}
where for $k_{p+1}=k_p+1$ the algebras 
$sp\,(k_{p+1}-k_p)=M_{k_p}$.
\end{theorem}
For contraction of maximal dimension $\bj=\biota$
we get algebra $M=\sum\limits_{k=1}^{n}\!\!{}^\oplus M_k$.
\setcounter{equation}{0}
\section{Representation of Cayley-Klein categories: \\terminology and examples}
In the next Sections of the work, we follow the program of
categorification and categorical extension developed in 
\cite{ Neretin, MAdStr} and apply it to the problem of
categorification of the representation theory of  Cayley-Klein groups 
\cite{GMROrt,MWRCKSymC}.

 To define a category $\rm{CK} \bf{K}$ we require the following 
data:

(a) a set $\Ob(\rm{CK} \bf{K})$ of elements called the objects of the category $\rm{CK} \bf{K}$;

(b) for any two objects $V(\bj), W(\bj) \in \Ob(\rm{CK} \bf{K})$ a set 
$\Mor_{\rm{CK} \bf{K}} (V(\bj), W(\bj))$ is 
defined, called  the {\it morphisms} from $V(\bj)$ to $W(\bj)$ (when it is clear 
what the category in question is, we omit the index $\rm{CK} \bf{K}$ and merely 
write $\Mor(V(\bj), W(\bj))$);

(c) for any $P$ in $\Mor(V(\bj), V'(\bj))$ and $Q$ in $\Mor(V'(\bj), V''(\bj))$ their 
product $QP$ is defined in $\Mor(V(\bj), V''(\bj))$. 

\noindent
The product must be 
associative:  for any $P$ in $\Mor(V(\bj), V'(\bj))$, $Q$ in $\Mor(V'(\bj), V''(\bj))$, and 
$R$  in $\Mor(V''(\bj), V'''(\bj))$, the formula
\beq
R(QP)=(RQ)P \ \ \ \ \ \ \ \ \ {\rm holds;}
\eeq

(d) it is usually assumed that the set $\Mor(V(\bj), V(\bj))$ contains an element 
$1_{V(\bj)}$ called  the {\it identity} such that, for any $P\in 
\Mor(Y,(\bj) V(\bj))$, 
we have $P\cdot 1_{V(\bj)} = P$ and, for any $Q \in \Mor(V(\bj), W(\bj))$, we have 
$1_{V(\bj)}\cdot Q=Q$.

{\it The Cayley-Klein category $\rm{CK} \bf{GD}$}.  The objects of the category $\rm{CK} \bf{GD}$ are the
even-dimensional complex linear spaces $V(\bj)$ endowed with a
non-degenerate symmetric bilinear form $\{\cdot,\cdot\}_{V(\bj)}$.

Let $V(\bj), W(\bj) \in \Ob(\rm{CK} \bf{GD})$. We introduce in 
$V(\bj) \oplus W(\bj)$ the
bilinear form
\be
\{(v,w),(v',w')\}_{V(\bj)\oplus W(\bj)} = \{v,v'\}_{V(\bj)} - \{w,w'\}_{W(\bj)}.
\ee
The morphisms from $V(\bj)$ to $W(\bj)$ in the category $\rm{CK} \bf{GD}$ are of two
types:
\begin{description}
\item{(a)} maximal isotropic subspaces of $V(\bj)\oplus W(\bj)$;
\item{(b)} the formal element $\nully_{V(\bj),W(\bj)}$, more usually denoted
simply by $\nully$. We emphasize that this element does not
correspond to any subspace of $V(\bj)\oplus W(\bj)$.
\end{description}

\noindent
Next we need to define the product $QP$ of morphisms $P\in
\Mor(V(\bj), W(\bj))$ and $Q\in\Mor(W(\bj),Y(\bj))$. First, the product of $\nully$ and any
other morphism is  equal $\nully$, that is,
\beq
\nully_{W(\bj),Y(\bj)}\cdot P = \nully_{V(\bj),Y(\bj)}, 
\quad Q \cdot \nully_{V(\bj),W(\bj)} =
\nully_{V(\bj),Y(\bj)}
\eeq
for any $P \in \Mor(V(\bj), W(\bj))$ and $Q\in \Mor(W(\bj), Y(\bj))$. Now let $P$ and $Q$
be linear relations. If
\be
\label{262}
\ker(Q)\cap \Indef P=0,
\ee
then $Q$ and $P$ are multiplied as linear relations. If, on the other
hand, (\ref{262}) does not hold, then $QP =\nully_{V(\bj),Y(\bj)}$.

\noindent
{\it The Cayley-Klein category $\rm{CK} \bf{GA}$}. The objects of $\rm{CK} \bf{GA}$ are the finite-dimensional
complex linear Cayley-Klein spaces. 

\noindent
The morphisms 
$\Mor_{\rm{CK} \bf{GA}}(V(\bj),W(\bj))$ consist of all
possible linear relations  $V(\bj)\rightrightarrows W(\bj)$, together with the
formal element $\nully = \nully_{V(\bj),W(\bj)}$,  which is not identified with any
linear relation. Let $P \in \Mor_{\rm{CK} \bf{GA}}(V(\bj),W(\bj))$ and
$Q\in \Mor_{\rm{CK} \bf{GA}}(W(\bj),Y(\bj))$. 

\noindent
We
define their product $QP\in \Mor_{\rm{CK} \bf{GA}}(V(\bj),Y(\bj))$ by the following rule:
\begin{description}
\item{(a)} the product of null with any morphism is null:
\beq
\nully_{W(\bj),Y(\bj)}\cdot P = \nully_{V(\bj),Y(\bj)},\quad  
Q\cdot\nully_{V(\bj),W(\bj)} =
\nully_{V(\bj),Y(\bj)};
\eeq
\item{(b)} if $P\neq \nully$ and $Q\neq \nully$, and
\be
\label{272}
\ker\, Q\cap \Indef\, P=0,
\ee
\be
\label{273}
\im \, P+D(Q)=W,
\ee
then $Q$ and $P$ are multiplied as linear relations. Otherwise,
$QP=\nully$.
\end{description}
We introduce a topology on $\Mor(V(\bj), W(\bj))$ as follows. 

The set
$\Mor(V(\bj),W(\bj))\backslash \nully$ is endowed with the topology of the
disconnected union of $(\dim V(\bj) + \dim W(\bj))$ Grassmannians and the morphism
$\nully_{V(\bj),W(\bj)}$ is contained in just one open set, namely
the whole of $\Mor(V(\bj),W(\bj))$.
\begin{proposition}
\label{p271}
\begin{description}

\item{(a)} Let $P\in \Mor_{\rm{CK} \bf{GA}}(V(\bj),W(\bj))$ and\\ 
$Q \in \Mor_{\rm{CK} \bf{GA}}(W(\bj),Y(\bj))$, and
suppose that $QP \neq \nully$. Then
\be
\label{274}
\dim\, QP = \dim \, P + \dim\,  Q - \dim\, W .
\ee
\item{(b)} Multiplication of morphisms is a jointly continuous
operation.
\end{description}
\end{proposition}

\noindent
{\bf Proof} \ \ \ 

(a) Let $H=V(\bj)\oplus W(\bj)\oplus W(\bj)\oplus Y(\bj)$, and let $Z$ be the
subspace of vectors of the form $(v,w, w, y)$. We define the subspace
$T = P\oplus Q$ as the set of all vectors of the form $(v,w,w',y)$,
where $(v,w)\in Q$ and $(w',y) \in P$. In view of (\ref{273}) we have
$T + Z = H$.  Thus $T\cap Z$ has dimension
\beq
\dim \, Z + \dim\,  T - \dim\,  H = \dim \, P + \dim \, Q - \dim\,  W.
\eeq

Next we denote the projection of $H$ onto $V(\bj) \oplus Y(\bj)$ along $W(\bj)\oplus
W(\bj)$ by $\pi$.  Then, as is easily seen, $\pi(T\cap Z)$ is the product
$QP$.  Furthermore, $\pi$ is injective on $T\cap Z$. In fact,
$\pi(v,w,w',y) = 0$ implies that $v = 0$ and $y = 0$, while
$(v,w,w',y)\in Z$ implies that $w = w'$.  Finally, $(0,w,w,0)\in T$
implies that $w$ lies in $\ker\, \cap\, \Indef\,  P$, which, by virtue of
(\ref{272}), consists merely of the origin.  This completes the proof of
assertion (a).  Assertion (b) is proved in similar fashion.
\hfill{$\Box$}

\vspace{0.2cm}
\noindent
{\it Duality}. We denote by $V'(\bj)$ the dual of $V(\bj)$, that is, the space of
linear functionals on $V(\bj)$.
Let $P\in\Mor_{\rm{CK} \bf{GA}}(V(\bj),W(\bj))$. The dual morphism 
$P'\in \Mor_{\rm{CK} \bf{GA}}(V'(\bj),W'(\bj))$
is constructed in terms of $P$ as follows. If $P$ is a non-null
element, then by definition, $P'$ consists of all pairs $(f',g')\in
V'(\bj)\oplus W'(\bj)$ such that $f(v) = f(w)$ for all $(v, w) \in P$.  Also, we
set $\nully_{V(\bj),W(\bj)} = \nully_{V'(\bj),W'(\bj)}$.

The dual morphisms are described slightly differently. Thus, consider the
subspace $P^0$ of $V(\bj)\oplus W(\bj)$ consisting of all vectors of the form
$(v,-w)$, where $(v, w) \in P$. Then $P'\subset(V(\bj)\oplus W(\bj))' = 
V'(\bj)\oplus
W'(\bj)$ is the same as the annihilator of the subspace $P^0$ of $V(\bj)\oplus W(\bj)$.

We denote by $\Ann\,  Q \subset H'$ the annihilator of the subspace
$Q\subset H$.

\begin{lemma}
\label{l272}
 Let $P : V(\bj) \rightrightarrows W(\bj)$ and $Q : W(\bj) \rightrightarrows 
 Y(\bj)$ be
morphisms of the Cayley-Klein category $\rm{CK} \bf{GA}$. Then
\begin{description}
\item{(a)} $(P')'=P$;

\item{(b)} $
\dim \, P' + \dim\,  P = \dim \, V(\bj) + \dim\,  W(\bj);
$

\item{(c)} $\ker\,  P' = \Ann \, D(P)$, $D(P') = \Ann\,  \ker\,  P$, $\Indef
\, P'=\Ann\, \im \, P$, and $\im \, P' =\Ann\, \Indef\,  P$;

\item{(d)} $Q'P' = \nully$ in the category $\rm{CK} \bf{GA}$ if and only if
$QP=\nully$ in the category $\rm{CK} \bf{GA}$;

\item{(e)} $(QP)'= Q'P'$.
\end{description}
\end{lemma}

\noindent
{\bf Proof} 
(a) $\Ann\,  \Ann \, T = T$.

(b) and (c) are simple exercises in linear algebra.

(d) Conditions (\ref{272}) and (\ref{273})  change places on passing to
the dual morphism.

(e) Let $Q,P$, and $QP$ be non-null morphisms. Let $(f'',f)\in Q'P'$.
Then there exists $f'\in W'(\bj)$ such that $(f'',f')\in Q'$ and $(f',f)\in
P'$.  Let $(y,v)$ be an  arbitrary element of $QP$. Then there exists
$w\in W$ such that $(y,w) \in Q$ and
$(v,w) \in P$. Then $f''(y) = f'(w) = f(v)$ by definition of the dual
morphism.  Thus $(f'',f')$ belongs to $(QP)'$, and hence $Q'P'\subset
(QP)'$.  On the other hand, the dimensions of $Q'P'$ and $(QP)'$ are
the same (see Proposition~\ref{p271}). This completes the proof.
\hfill{$\Box$}

\begin{lemma}
\label{l273}
 Multiplication of morphisms of the  category $CK{\bf GA}$ is
associative.
\end{lemma}
\noindent
{\bf Proof} \\  Let $P\in \, \Mor(V(\bj),W(\bj))$, $Q\in \Mor(W(\bj),Y(\bj))$, $R\in \Mor(Y(\bj),Z(\bj))$,
and suppose that $R(QP)=\nully$. This means that one of the following
four conditions holds:

(1) $\ker\,  Q \cap \Indef\,  P \neq O$,

(2) $\ker\,  R \cap \Indef \, (QP) = 0$,

(3) $\im\,  P+D(Q) \neq W$,

(4) $D(R)+\im\,  QP \neq Y$.

Suppose first that (1) holds. Clearly $\ker(RQ)\supset \ker Q$, and
therefore we have $\Indef\,  P\cap \ker(RQ)\neq 0$. Hence $(RQ)R = \nully$.

Suppose next that (2) holds. Then $y \in \Indef(QP)\cap \ker R$
contains a non-zero vector. Take $w\in W$ such that $(0,w) \in P$ and
$(w,y) \in Q$.  Then $(y,0)\in R$, which implies that $w\in \Indef
P\cap \ker(RQ)$.  Hence $(RQ)P = \nully$.

The cases (3) and (4) reduce to (1) and (2) by going over to the dual
morphisms. \hfill{$\Box$}

\begin{proposition}
\label{p274}
 The definition of the Cayley-Klein category $CK{\bf GD}$ 
is correct.
\end{proposition}

\noindent
{\bf Proof}  We have just verified associativity. Suppose  that $P\in
\Mor(V(\bj),W(\bj))$, $Q\in \Mor(W(\bj),Y(\bj))$ and $QP\neq \nully$. Then, by Proposition
 \ref{p271},
\beq
\dim \, QP = {1\over2}\left(\dim \, V(\bj) + \dim \, W(\bj)\right).
\eeq
Hence the isotropic subspace $QP$ has maximum possible dimension.
\hfill{$\Box$}

Let $\rm{CK} \bf{Op}$ be the category of {\it linear Cayley-Klein spaces and operators}. 
The objects of the category $\rm{CK} \bf{Op}$ are the finite-dimensional linear 
Cayley--Klein spaces, and the morphisms are the linear operators.

 A {\it covariant functor} $(F,\phi)$ 
from a category $CK{\bf K}$ to a category $\rm{CK} \bf{Op}$ is determined by the following 
data:

(i) a map $F:\Ob(CK{\bf K}) \rightarrow \Ob(\rm{CK} \bf{Op})$;

(ii) the collection of maps 
\beq
\phi_{V(\bj),W(\bj)}:\Mor(V(\bj), W(\bj))\rightarrow 
\Mor(F(V(\bj)),  F(W(\bj))) 
\eeq
is defined for all $V(\bj), W(\bj) \in \Ob(CK{\bf K})$, and these maps 
must satisfy the condition that
\beq
\phi_{V(\bj),Y(\bj)}(PQ) =\phi_{W(\bj),Y(\bj)}(P) \phi_{V(\bj),W(\bj)}(Q), \qquad \phi(1_V(\bj))=1_{F(V(\bj))}.
\eeq
Similarly, by a {\it contravariant functor} from a category $CK{\bf K}$ to a 
category $\rm{CK} \bf{Op}$, we mean a map $F:\Ob(CK{\bf K}) \rightarrow \Ob(\rm{CK} \bf{Op})$ and a 
collection of maps
\beq
\phi_{V(\bj),W(\bj)}:\Mor(V(\bj),W(\bj)) \rightarrow \Mor(F(W(\bj)),F(V(\bj))),
\eeq
defined for all $V(\bj), W(\bj) \in\Ob(\Ob(CK{\bf K})$, such that
\beq
\phi_{V(\bj),Y(\bj)}(PQ) =\phi_{V(\bj),W(\bj)}(Q) \phi_{W(\bj),Y(\bj)}(P).
\eeq
By a representation of a Cayley-Klein category $\rm{CK} \bf{K}$ we mean a covariant functor 
$(T,\tau)$ from $\rm{CK} \bf{K}$ to the category $\rm{CK} \bf{Op}$. In other words, we associate 
with each $V(\bj)
\in\Ob(\rm{CK} \bf{K})$ a linear space $T(V(\bj))$, and with each morphism 
$P:V(\bj) \rightarrow W(\bj)$ an operator $\tau(P)$ from $T(V(\bj))$ to $T(W(\bj))$, such 
that for each triple $(V(\bj),W(\bj),Y(\bj))\in\Ob(\rm{CK} \bf{K})$, and for each 
$P\in\Mor(V(\bj),W(\bj))$ 
and $Q\in\Mor(W(\bj),Y(\bj))$ we have 

\begin{equation}
\tau(QP) =\tau(Q) \tau(P).
\end{equation}

\noindent
Contravariant functor from $\rm{CK} \bf{K}$ to $\rm{CK} \bf{Op}$ gives antirepresentation.

{\it Topological Cayley-Klein categories}. Let $\rm{CK} \bf{K}$  be a category, and suppose that 
a topology is introduced on each set $\Mor(V(\bj),W(\bj))$. Then 
$\rm{CK} \bf{K}$ 
is called a {\it topological category} if for any objects
$V(\bj),W(\bj),Y(\bj)$ and  any 
$P_0\in\Mor(V(\bj),W(\bj))$, and  any $Q_0 \in \Mor(W(\bj),Y(\bj))$ the maps $Q\mapsto QP_0$ from 
$\Mor(W(\bj),Y(\bj))$ to $\Mor(Y(\bj),Y(\bj))$ and $P\mapsto Q_0P$ from 
$\Mor(V(\bj),W(\bj))$ to 
$\Mor(V(\bj), Y(\bj))$ are continuous. We emphasize that we only require separate 
continuity for multiplication of morphisms.

Representations  of topological categories must, of course, satisfy 
certain continuity requirements. Let $(R,\rho)$ be a representation of 
the category $\rm{CK} \bf{K}$. We require that all the spaces $R(V(\bj))$ be complete 
locally convex spaces  and that all the operators $\rho(L)$ be bounded. We denote 
by $B(H_1,H_2)$ the set of all bounded operators from $H_1$ to $H_2$. 
The functions $Q \mapsto\rho(Q)$ from $\Mor(V(\bj),W(\bj))$ to 
$B(R(V(\bj)),R(W(\bj)))$ 
must be continuous.  An ambiguity can arise here: there exist several 
natural topologies on the set $B(R(V(\bj)),R(W(\bj)))$.  We shall always assume 
that $B(R(V(\bj)),R(W(\bj)))$ is endowed with the weak topology (that is, the 
topology in which for each $h\in R(V(\bj))$ and each continuous linear 
functional $l$ on $R(W(\bj))$ the function
\beq
\phi_{h,l}(Q)=l(Qh)
\eeq
is continuous).  In other words, we require continuity of all the 
``matrix elements'' 
\beq
m_{h,l}(P)=l(\rho(P)h)
\eeq
of the representation $(R,\rho)$.

Next we define continuity of the projective representation $(R,\rho)$. 
We consider the quotient spaces $B(R(V(\bj)),R(W(\bj)))/\Cb^*$ (we identify the 
operators $A$ and $\lambda A$, where $\lambda\subset\Cb^*$) endowed 
with the natural quotient topology which, incidentally, is not 
Hausdorff.  All the functions $Q\mapsto\rho(Q)$ from $\Mor(V(\bj), W(\bj))$ to 
$B(R(V(\bj)), R(W(\bj)))$ must be continuous.

\vspace{2mm}
\noindent
{\it Subrepresentations}. Let $(T,\tau)$ be a representation of the 
category $\rm{CK} \bf{K}$.  Suppose that a closed subspace $A(V(\bj))$
is 
chosen in each 
space $T(V(\bj))$ such that $\tau(P)A(V(\bj))$ lies in $A(W(\bj))$ for any 
$V(\bj),W(\bj)\in\Ob(\rm{CK} \bf{K})$ and $P\in \Mor(V(\bj),W(\bj))$.  Then we obtain the
following representation of the category $\rm{CK} \bf{K}$: associated with each 
object $V(\bj)$ is the space $A(V(\bj))$, and with each morphism 
$P:V(\bj)\rightarrow  W(\bj)$ the restriction $a(P)$ of the operator 
$\tau(P)$ to $A(V(\bj))$. In this 
case we say that $(A,a)$ is a {\it subrepresentation} of the 
representation $(T,\tau)$.

A representation $(T,\tau)$ is said to be {\it irreducible} if it has 
no subrepresentations other than itself and the zero representation.

Let $(T,\tau)$ be a representation of $\rm{CK} \bf{K}$, and let $S$ be a subset of 
some space $T(V(\bj))$. Let $A(W(\bj))$ be the closed linear span of the set of 
all vectors of the form $\tau(P)h$, where $P\in \Mor(V(\bj), W(\bj))$
and 
$h\in 
S$.  Then it is not difficult to see that the collection of subspaces 
$A(W(\bj))$ gives a subrepresentation of $(T,\tau)$.

This subrepresentation is called the {\it cyclic} span of the set $S$.

Let $(T,\tau)$ be a representation of the category $\rm{CK} \bf{K}$. Then the 
operators $\tau(P)$ of the semigroup $\End(V(\bj))$ and also of the group 
$\Aut(V(\bj))$ act in each space $T(V(\bj))$.  These representations of 
$\End(V(\bj))$ 
and $\Aut(V(\bj))$ are said to be representations {\it subordinated} to the 
given representation.

\begin{lemma} Let $(T,\tau)$ be an irreducible representation. Then 
all the subordinated representations of the semigroups $\End(V(\bj))$ are 
irreducible.
\end{lemma}
\noindent
{\bf Proof} 

 Suppose that $T(V(\bj))$ contains an $\End(V(\bj))$-invariant 
subrepresentation $A$.  Then the cyclic span of $A$ is a non-trivial 
subrepresentation of $T$. \hfill{$\Box$}


We note that  the converse is false in general. It is, however, true 
for ordered Cayley-Klein categories (see Subsection 5.3).

\vspace{0.2cm}
\noindent
{\it Direct sums}. Let $(T,\tau)$ be a representation of the category 
$\rm{CK} \bf{K}$, let $A_1,A_2$ be subrepresentations of $(T,\tau)$ and suppose that
\beq
T(V(\bj))=A_1(V(\bj))\oplus A_2(V(\bj))
\eeq
for any $V(\bj)\in \Ob(\rm{CK} \bf{K})$. Then we say that $(T,\tau)$ is the {\it direct 
sum} of the subrepresentations $A_1, A_2$.

We say that $(T,\tau)$ is {\it decomposible} if $(T,\tau)$ decomposes 
into a direct sum of subrepresentations $A_1,A_2,\ldots$\,. 
Incidentally, even in the finite-dimensional case, (that is, in the 
case when all the $T(V(\bj))$ are finite-dimensional) the number of these 
$A_j$ can be infinite (if for any $V(\bj)$ only finitely many spaces 
$A_j(V(\bj))$ differ from 0).

Finally, if we are given a family of linear representations 
$(T_i,\tau_i)$ of a family $\rm{CK} \bf{K}$, then their 
{\it exterior direct sum} 
$(S,\sigma)$ is defined: the spaces $S(V(\bj))$ of the representation are 
$\bigoplus T(V_i(\bj))$, and the operators $\sigma(S(P))$ are 
$\bigoplus\tau_i(P)$.

A representation is said to be {\it completely reducible} if it 
decomposes into a direct sum of irreducible representations.

\vspace{0.2cm}
\noindent
{\it Intertwining transformations}. Let $(T,\tau)$ and $(T',\tau')$ be 
two representations of a category $\rm{CK} \bf{K}$. Then by an 
{\it intertwining 
transformation} $A:(T,\tau)\rightarrow  (T',\tau')$ we mean a 
collection of bounded operators $A(V(\bj)): T(V(\bj))\rightarrow  
T'(V(\bj))$ such 
that for all $V_1(\bj),V_2(\bj)\in\Ob(\rm{CK} \bf{K})$ and 
$P\in\Mor(V_1(\bj),V_2(\bj))$ the diagram
\makeatletter
$$
\begin{CD}
T({V_1(\bj)}) @>\tau(P)>>  T({V_2(\bj)})\\
@V{A}({V_1(\bj)})VV @VV{A}({V_2(\bj)})V\\
T'({V_1(\bj)}) @>\tau'(P)>> T'({V_2(\bj)})
\end{CD}
$$

\makeatother

\noindent
commutes (that is, $\tau'(P)A(V_1(\bj)) = A(V_2(\bj))\tau(P))$. If we are dealing 
with projective representations then we require that both $\tau'(P)A(V_1(\bj))$ 
and $A(V_2(\bj))\tau(P)$ should coincide to within multiplication by a 
constant.

Let $(T,\tau)$ and $(T',\tau')$ be representations of the category 
$\rm{CK} \bf{K}$. 
Suppose further that $A:(T,\tau)\rightarrow  (T',\tau')$ and 
$B:(T',\tau')\rightarrow  (T,\tau)$ are intertwining transformations 
such that $A(V(\bj))B(V(\bj)) = E$, and $B(V(\bj))A(V(\bj)) = E$ for all $V(\bj)$. Then we say 
that $(T,\tau)$ and $(T',\tau')$ are {\it equivalent representations}.

It should be noted that our definition of equivalence is valid only  in 
the case of finite-dimensional representations and $*$-representations 
(see below); it does not work in the general case. 

\vspace{0.2cm}
\noindent
{\it $*$-representations}. Suppose that for any 
$V(\bj),W(\bj) \in\Ob(\rm{CK} \bf{K})$ we are 
given a map $s:P \mapsto P^*$ from $\Mor(V(\bj),W(\bj))$ to 
$\Mor(W(\bj), V(\bj))$. Then we 
call $s$ an {\it involution} if the following identities hold:
\be
P^{**}=P, \qquad (PQ)^* =Q^*P^*.
\ee
Let $(T,\tau)$ be a representation of the category $\rm{CK} \bf{K}$. 
Then we say 
that $(T,\tau)$ is a $*$-representation if the $T(V(\bj))$ are Hilbert 
spaces and for each $P\in\Mor(V(\bj), W(\bj))$
\beq
\tau(P)^* =\tau(P^*).
\eeq
\begin{lemma} Let $(T,\tau)$ be a $*$-representation of the category 
$\rm{CK} \bf{K}$, and let $A$ be a subrepresentation of it. We choose in
 each 
$T(V(\bj))$ 
the orthocomplement of $A(V(\bj))$ which we denote by $B(V(\bj))$. 
Then $B(V(\bj))$ is 
also a subrepresentation.
\end{lemma}
\noindent
{\bf Proof} 
Let $v\in B(V(\bj))$. For each $P\in\Mor(V(\bj), W(\bj))$ and each 
$w\in A(W(\bj))$ we have 
\beq
\langle Pv,w\rangle_{W(\bj)} =\langle v,P^* w\rangle_{V(\bj)}=0,
\eeq
since $P^*w \in A(V(\bj))$. Hence $Pv\in B(W(\bj))$. \hfill{$\Box$}

We call an element $P\in \Aut_{\rm{CK} \bf{K}}(V(\bj))$ {\it unitary} if $P$ is invertible 
and $P^* = P^{-1}$ (in other words, $P^*P = PP^* = E$). We denote the 
set of all unitary elements of $\Aut_{\rm{CK} \bf{K}}(V(\bj))$ by $\Aut^*_{\rm{CK} \bf{K}}(V(\bj))$. Clearly 
$\Aut^*_{\rm{CK} \bf{K}}(V(\bj))$ is a subgroup of $\Aut_{\rm{CK} \bf{K}}(V(\bj))$. Let $\tau$ be a linear 
$*$-representation of the category $\rm{CK} \bf{K}$.  Then the subordinated 
representation of the group $\Aut^*_{\rm{CK} \bf{K}}(V(\bj))$ is clearly unitary. 
\begin{lemma} Let $\tau$ be a projective $*$-representation of the 
category $\rm{CK} \bf{K}$.  Then, for any 
$P\in \Aut^*_{\rm{CK} \bf{K}}(V(\bj))$, the operator $\tau(P)$ 
is unitary to within multiplication by a constant.
\end{lemma}
\noindent
{\bf Proof} 
In fact,
\beq
\tau(P)^*\tau(P)  &\!\!\!=\!\!\!& \lambda\tau(P^*P) =\lambda\tau(1_V) 
=\lambda E,\\
\tau(P)\tau(P^*)  &\!\!\!=\!\!\!& \mu\tau(PP^*) =\mu\tau(1_V)=\mu E
\eeq
for some $\lambda,\mu\in\Cb$. We then observe that $\tau(P)^*\tau(P)$ 
and $\tau(P)\tau(P^*)$ are positive self-adjoint operators, from which 
it follows that $\lambda> 0$ and $\mu > 0$. But
\beq
\lambda\tau(P)=\tau(P)(\tau(P)^*\tau(P)) =(\tau(P)\tau(P^*)) \tau(P) 
=\mu\tau(P).
\eeq
Consequently $\lambda=\mu$, and therefore the operator 
$\lambda^{-1/2}\tau(P)$ is unitary. \hfill{$\Box$}

\vspace{0.2cm}
\noindent
{\it Tensor products}. Let $(T,\tau)$ and $(T',\tau'$) be 
finite-dimensional representations or $*$-representations of the 
category $\rm{CK} \bf{K}$. We define their {\it tensor product}
\beq
(T\otimes T', \tau\otimes\tau')
\eeq
as follows. The spaces 
$(T\otimes T')(V(\bj))$ are $T(V(\bj))\otimes T'(V(\bj))$, and 
\beq
(\tau\otimes \tau')(P): T(V(\bj))\otimes T'(V(\bj)) 
\rightarrow T(W(\bj))\otimes T'(W(\bj))
\eeq
are $\tau(P)\otimes \tau'(P)$. We then define the $k$th {\it exterior 
power} of the representation $(T, \tau)$ by setting
\beq
(\Lambda^kT)(V(\bj)):=\Lambda^k(T(V(\bj))), \qquad 
(\Lambda^k\tau)(P):=\Lambda^k\tau(P).
\eeq
The $k$th {\it symmetric power} of a representation is defined 
similarly.

\vspace{0.2cm}
\noindent
{\it Equivalence of categories}. Let $\rm{CK} \bf{K}$ and $\rm{CK} \bf{L}$ be 
two categories. It 
might appear that we should define an {\it isomorphism} by requiring 
that there exist a bijection $H : \Ob(\rm{CK} \bf{K})\leftrightarrow  
\Ob(\rm{CK} \bf{L})$ 
and a family of bijections
\beq
h_{V(\bj),W(\bj)}: \Mor_{\rm{CK} \bf{K}}(V(\bj),W(\bj)) \leftrightarrow 
\Mor_{\rm{CK} \bf{L}}(H(L(\bj)), H(W(\bj)))
\eeq
such that
\beq
h_{W(\bj),Y(\bj)}(Q) h_{V(\bj),W(\bj)}(P) =h_{V(\bj),Y(\bj)}(QP)
\eeq
for all $Q \in\Mor_{\rm{CK} \bf{K}}(W(\bj),Y(\bj))$ and $P\in \Mor_{\rm{CK} \bf{K}}(V(\bj),W(\bj))$. Unfortunately, this 
definition is inadequate, as the following example shows.

\begin{example} Let $\rm{CK} \bf{K}$ be a category of finite-dimensional 
complex Cayley-Klein spaces and linear operators. Let $\bf{L}$ be the category whose 
objects are the natural numbers and whose morphisms $m\rightarrow  n$ 
are $m \times n$ matrices. Then $\rm{CK} \bf{K}$ and $\bf{L}$ are not isomorphic because 
there is no bijection $\Ob(\rm{CK} \bf{K})\leftrightarrow \Ob(\bf{L})$.
\end{example}

Let $\rm{CK} \bf{K}$ be a category. We say that $V(\bj),W(\bj) \in\Ob(\rm{CK} \bf{K})$ are {\it isomorphic}
if there exists a pair of isomorphisms $P:V(\bj) \rightarrow  W(\bj)$ and 
$Q:W(\bj)\rightarrow V(\bj)$ such that $PQ = 1$ and $QP = 1$. We form a 
{\it reduced category} $\rm{CK} \bf{K'}$ as follows.  We choose one object from each 
isomorphism class of objects of $\rm{CK} \bf{K}$. This will be our set 
$\Ob(\rm{CK} \bf{K'})$. 
The sets $\Mor_{\rm{CK} \bf{K'}}(V(\bj),Y(\bj))$ are the same as 
$\Mor_{\rm{CK} \bf{K}}(V(\bj),Y(\bj))$.

Two categories $\rm{CK} \bf{K}$ and $\rm{CK} \bf{L}$ are called {\it equivalent} if their reduced 
categories are isomorphic.                        

\begin{example} The categories $\rm{CK} \bf{K}$ and $\bf{L}$ in Example 4.1 are 
equivalent.
\end{example}

\noindent
{\it Subcategories}. We say that a category $\rm{CK} \bf{L}$ is a {\it subcategory} 
of $\rm{CK} \bf{K}$ if $\Ob(\rm{CK} \bf{L})\subset \Ob(\rm{CK} \bf{K})$
and 
$\Mor_{\rm{CK} \bf{L}}(V(\bj), W(\bj))\subset 
\Mor_{\rm{CK} \bf{K}}(V(\bj), W(\bj))$ 
for any $V(\bj), W(\bj) \in\Ob(\rm{CK} \bf{L})$.

\vspace{0.2cm}
\noindent
{\it Quotient representations}. Let $(T,\tau)$ be a representation of 
the category $\rm{CK} \bf{K}$. Let $(A,\alpha)$ be a subrepresentation of it. Then 
we define the {\it quotient representation}
\beq
(S,\sigma)=(T,\tau)/(A,\alpha)
\eeq
as follows. The spaces $S(V(\bj))$ are $T(V(\bj))/A(V(\bj))$, and the operators
\beq
\sigma(P) : T(V(\bj))/A(V(\bj)) \rightarrow  T(W(\bj))/A(W(\bj))
\eeq
are the natural quotient maps induced by the map $\tau(P): 
T(V(\bj))\rightarrow T(W(\bj))$.

\begin{example} {\rm{(See \cite{SRCKOrtC})}} 

The \underline{functor Spin}. 
Let $V(\bj)\in\Ob(\rm{CK} \bf{GD})$ with $\dim V(\bj)=2n$. We 
decompose $V(\bj)$ into a direct sum $V(\bj)=V_+(\bj)\oplus V_-(\bj)$ of maximal isotropic 
subspaces. 

\noindent
Next, we choose bases $\{e_1^+(V(\bj)),\ldots, e_n^+(V(\bj))\}$ and 
$\{e_1^-(V(\bj)),\ldots, e_n^-(V(\bj))\}$ in $V_+(\bj)\oplus V_-(\bj)$ such that, for 
$i,j\in\{1,\ldots,n\}$, $\{e_i^+(V(\bj)), e_j^+(V(\bj))\}=
 \{e_i^-(V(\bj)), 
e_j^-(V(\bj))\}=0$, $\{e_i^+(V(\bj)), e_j^-(V(\bj))\}=\delta_{ij}$. We denote the 
coordinates of the vector $v\in V$ in this basis by $(v_1^+,\ldots, 
v_n^+; v_1^-,\ldots,v_n^-)$.

Thus $V(\bj)$ is identified with the coordinate space $V_{2n}(\bj)$. Similarly, 
the bilinear forms in $V(\bj)$ and $V_{2n}(\bj)$ are identified. We associate 
with each object $V(\bj)$ in $\rm{CK} \bf{GD}$ the subspace ${\rm 
Spin}(V(\bj)):=\Lambda(V_+(\bj))$, that is, the exterior algebra on $V_+(\bj)$ (or the 
Grassmann algebra on $\Lambda_n$).

\end{example}
\begin{example} {\rm{(See \cite{SRCKOrtC})}} \label{e4.4}

\underline{Projective 
representation of the category $\rm{CK} \bf{GD}$.} \\
 (a) {\it Let $P\in\Mor(V(\bj),W(\bj))$ be a non-null morphism. 
Then there exists a non-zero operator ${\rm Spin}(P):\Lambda (V_+(\bj)) 
\rightarrow \Lambda (W_+(\bj))$ such that }
\begin{equation}
\hat{a}(\omega){\rm Spin}(P) ={\rm Spin}(P)\hat{a}(v)
\end{equation}
{\it for all $(v,\omega)\in P$. Further, ${\rm Spin}(P)$ is unique to 
within proportionality.}\\ 
(b) {\it Let $P\in\Mor(V(\bj),W(\bj))$ and $Q\in\Mor(W(\bj),Y(\bj))$ 
be non-null morphisms.  Then, in the case where $QP\neq\nully$, we have}
\begin{equation}
{\rm Spin}(Q){\rm Spin}(P) =\lambda(P,Q){\rm Spin}(QP),
\end{equation}
{\it where $\lambda(P,Q)\in\Cb^*$. \\
If, on the other hand, $QP=\nully$, 
then ${\rm Spin}(Q){\rm Spin}(P)=0$. Thus the map \\
$P\rightarrow {\rm 
Spin}(P)$, $\nully\rightarrow 0$ 
is a projective 
representation of the category $\rm{CK} \bf{GD}$.}
\end{example}
\begin{example}  {\rm{(See \cite{SRCKOrtC})}} \label{e4.5}

\underline{Fundamental representation of the category $\rm{CK} \bf{GA}$}. 
First we 
construct a functor $(M,\mu)$ from the category $\rm{CK} \bf{GA}$ to the 
category 
$\rm{CK} \bf{GD}$. Let $V(\bj)\in\Ob(\rm{CK} \bf{GA})$. We set 
$M(V(\bj))=V(\bj)\oplus V'(\bj)$ and endow $M(V(\bj))$ with 
the symmetric bilinear form $\{(v_1,f_1), 
(v_2,f_2)\}=f_1(v_2)+f_2(v_1)$, so that $M(V(\bj))$ is an object of the 
category $\rm{CK} \bf{GD}$.

Let $P\in\Mor_{\rm{CK} \bf{GA}}(V(\bj),W(\bj))|\nully$. We set
\begin{equation}
\mu(P): =P\oplus P' \subset (V(\bj)\oplus W(\bj))\oplus(V'(\bj)\oplus W'(\bj))''.
\end{equation}
Finally we set $\mu(\nully)=\nully$.

The composition of the functor Spin and the functor $M$ gives us a 
projective representation $(\Lambda,\lambda)$ of the category $\rm{CK} \bf{GA}$. 
This representation is called the fundamental representation of the 
category $\rm{CK} \bf{GA}$.

\end{example}

\setcounter{equation}{0}


\section{Representations of complex classical \\Cayley-Klein categories}
\subsection{ The complex classical Cayley-Klein categories \\
$\rm{CK} \bf{A}, 
\rm{CK} \bf{B}, 
\rm{CK} \bf{C}, \rm{CK} \bf{D}$}

The category $\rm{CK} \bf{A}$ has the linear spaces as its objects, and 
the linear
operators defined to within multiplication by a non-zero constants as
its morphisms.

An object of category $\rm{CK} \bf{B}$ is  an odd-dimensional complex 
linear space
$V(\bj)$ endowed with a non-degenerate symmetric bilinear form
$\{\cdot,\cdot\}_{V(\bj)}$. Let $V(\bj), W(\bj)$ be objects of  the
category 
$\rm{CK}
\bf{B}$. 
We
introduce in $V(\bj)\oplus W(\bj)$ the symmetric bilinear form
\be \label{5.1}
\{(v,w), (v', w')\}_{V(\bj)\oplus W(\bj)} = \{v, v'\}_{V(\bj)} -
 \{w,w'\}_{W(\bj)}.
\ee
The set $\Mor_{\rm{CK} \bf{B}}(V(\bj),W(\bj))$ consists of null and 
the maximal isotropic
subspaces of $V(\bj)\oplus W(\bj)$. The morphisms are multiplied in the same
way as in the category $\rm{CK} \bf{GA}$. The verification that this
multiplication is well defined is done in the same way as for
 the category $\rm{CK} \bf{GD}$; it is based on Proposition \ref{p271}.

An object of the category $\rm{CK} \bf{C}$ is a finite-dimensional 
complex linear
space $V(\bj)$ endowed with a non-degenerate skew-symmetric bilinear form
$\{\cdot, \cdot\}_{V(\bj)}$.  If $V(\bj), W(\bj)$ are objects in the 
category $\rm{CK} \bf{C}$,
then we introduce in $V(\bj)\oplus W(\bj)$ the skew-symmetric bilinear form
defined by (\ref{5.1}). The set $\Mor_{\rm{CK} \bf{C}}(V(\bj),W(\bj))$ 
consists of null and
all maximal isotropic subspaces of $V(\bj)\oplus W(\bj)$. The morphisms are
multiplied in the  same way as  in the category $\rm{CK} \bf{GA}$.

The automorphism groups of an object $V(\bj)$ for the above categories $\rm{CK} \bf{A}, 
\rm{CK} \bf{B}, 
\rm{CK} \bf{C},$ are $GL(V(\bj))$, $O(V(\bj))$, and $Sp(V(\bj))$,
respectively. 

The next
category in the list might have been the category $\rm{CK} \bf{GD}$ related to  the
series of groups $O(2n;\bj,\Cb)$.
However we prefer the category $\rm{CK} \bf{D}$, which is close to it 
and is related
to the groups $SO(2n;\bj,\Cb)$.

{\it The category $\rm{CK} \bf{D}$}. Let $Y(\bj)$ be a complex even-dimensional space
endowed with a non-degenerate symmetric bilinear form.  We denote by
$ {\mathcal Gr}(Y(\bj))$ the Grassmannian of all maximal isotropic subspaces of $Y(\bj)$.
\begin{lemma} \label{l5.1}
The set ${\mathcal Gr}(Y(\bj))$ has two connected components. Two
subspaces $H_1, H_2\in {\mathcal Gr}(Y(\bj))$ lie in the same connected component if
and only if the codimension  of $H_1\cap  H_2$ in $H_1$ and $H_2$ is
even.
\end{lemma}
\noindent
{\bf Proof} 
By Witt's theorem \cite{Burbaki} the space ${\mathcal Gr}(Y(\bj))$ is
$O(2n;\bj, \Cb)$-homogeneous. The stabiliser of a point is
isomorphic to $GL(n;\bj,  \Cb)$ and is therefore connected. Hence ${\mathcal Gr}(Y(\bj))$
has two connected components. Let $V(\bj)$ be a fixed member of 
${\mathcal Gr}(Y(\bj))$,
and let $W(\bj)\in {\mathcal Gr}(Y(\bj))$. We choose a basis $e_1,\ldots, e_s$,
$f_1,\ldots,f_k$, $f'_1,\ldots, f'_k$, $g_1,\ldots, g_s$ in $Y(\bj)$ such
that
\beq
\{e_i,g_j\}=\delta_{ij}, \qquad \{ f_i,f'_j\} =\delta_{ij},
\eeq
\beq
\{e_i,e_j\}=\{e_i,f_j\}=\{e_i,f'_j\} =\{g_i,f_j\} = \{g_i,f'_j\}
=\{g_i,g_j\}=0,
\eeq
\beq
e_i\in V(\bj)\cap W(\bj),\quad f_j\in V(\bj), \, f_j\in W(\bj), \, 
g_i\notin V(\bj), \,
g_i\notin W(\bj).
\eeq
Suppose that $k$ is even. We construct a curve $V(t;\bj)$ in 
${\mathcal Gr}(Y(\bj))$
joining $V(\bj)$ and $W(\bj)
$ ($0\leq t \leq 1$). By definition, the space $V(t;\bj)$
is spanned by vectors of the form
\be \label{n.5.2}
e_i,t\left( f_{2j-1} \pm f_{2j}\right) +(1-t) \left(f'_{2j-1} \mp
f'_{2j}\right).
\ee
Suppose that $k$ is  odd. We construct a curve $\tilde{V}(t;\bj)$ in
${\mathcal Gr}(Y(\bj))$ such that the space $V(t;\bj)$ is spanned by vectors of the form
(\ref{n.5.2}) and also $f'_k$.  It is clear that $V(0) =W$  and $V(1)$ is
spanned by the vectors $e_1,\ldots,e_s$, $f_1,\ldots,f_{k-1}$, $f'_k$,
so that $V(1)\cap V$ has codimension 1 in $V$.  Let $\gtL$ be the set
of all subspaces of $V(\bj)$ of codimension 1, and let $M$ be the set of all
$W(\bj)\in {\mathcal Gr}(Y(\bj))$ such that $W(\bj)\cap V(\bj) \in\gtL$.   The map $\phi:W(\bj)\mapsto W(\bj)
\cap V(\bj)$ from $M$ to $\gtL$  is a continuous bijection
($\phi^{-1}(H)\simeq H^\bot$). But $\gtL$  is connected, which means
that $M$ is connected and hence the set of al $Q\in {\mathcal
Gr}(Y(\bj))$ 
such that
$Q\cap V$ has odd codimension in $V(\bj)$ is also connected. \hfill{$\Box$}


An object ot the category $\rm{CK} \bf{D}$ is an object $V(\bj)$ of the category $\rm{CK} \bf{GD}$
in which  one of the connected components ${\mathcal Gr}_+(V(\bj))$ of the set
${\mathcal Gr}(V(\bj)) = Mor_{\rm{CK} \bf{GD}}(0,V(\bj))\backslash\nully$
is fixed. Let $V(\bj),W(\bj)$ be objects in $\rm{CK} \bf{D}$. 
Let $W_+(\bj)\in{{\mathcal Gr}}_+(W(\bj))$,
$V_+(\bj)\in{{\mathcal Gr}}_+(V(\bj))$, and let $W_-(\bj)\in {{\mathcal Gr}}(W(\bj))$ be a
complement of $W_+(\bj)$ with respect to $W(\bj)$. (Note that $W_-(\bj)$ is either
contained or is not contained in ${\mathcal Gr}_+(W(\bj))$ depending on whether the
number ${1\over2}\dim W(\bj)$ is even or odd.) The set 
$\Mor_{\rm{CK} \bf{D}}(V(\bj), W(\bj))$
consists of null and the connected component of
$\Mor_{\rm{CK} \bf{GD}}(V(\bj),W(\bj))\backslash\nully =
{\mathcal Gr}(V(\bj)\oplus W(\bj))$  containing
$V_+(\bj)\oplus W_-(\bj)$. The morphisms are multiplied in the same way as in
the category ${\rm{CK} \bf{GD}}$.
\begin{lemma} Multiplication of morphisms is well defined in the
sense that if $P\in\Mor_{\rm{CK} \bf{D}}(V(\bj),W(\bj))$ and $Q \in 
\Mor_{\rm{CK} \bf{D}}(W(\bj),Y(\bj))$,
then $QP\in\Mor_{\rm{CK} \bf{D}}(V(\bj),Y(\bj))$.
\end{lemma}

\noindent
{\bf Proof} 
From continuity considerations it suffices to verify that
for some $P_0$ in $\Mor_{\rm{CK} \bf{D}}(V(\bj),W(\bj))$ and $Q_0$ 
in $\Mor_{\rm{CK} \bf{D}}(W(\bj),Y(\bj))$ such that
$Q_0P_0\neq\nully$, the product $Q_0P_0$ belongs to 
$\Mor_{\rm{CK} \bf{D}}(V(\bj),Y(\bj))$. For
$P_0$ and $Q_0$ we can take $V_+(\bj)\oplus W_-(\bj)$ and $W_+(\bj)\oplus Y_-(\bj)$
($V_\pm(\bj)$ and $Y_\pm(\bj)$ are defined in the same way as $W_\pm(\bj)$).
\hfill{$\Box$}


It it easy to see that the group  $\Aut_{\rm{CK} \bf{D}}(V(\bj))$ is 
the same as $SO(V(\bj))$
 (in fact, $SO(V(\bj))$ must be entirely contained in one of the two 
components of the manifold $\Mor_{\rm{CK} \bf{GD}}(V(\bj), V(\bj))\backslash\nully$; 
furthermore, it is easy to see that identity element of the group lies 
in $\Mor_{\rm{CK} \bf{D}}(V(\bj),V(\bj))\backslash\nully$).


\subsection{ Classification theorems}

\begin{theorem} \label{t.5.1} (a) The holomorphic projective representations of the
categories $\rm{CK} \bf{A}$, $\rm{CK} \bf{B}$, $\rm{CK} \bf{C}$,  $\rm{CK} \bf{D}$, $\rm{CK} \bf{GD}$ are completely reducible.

(b) The irreducible holomorphic representations of the categories
$\rm{CK} \bf{A}$, $\rm{CK} \bf{B}$, $\rm{CK} \bf{C}$,  $\rm{CK} \bf{D}$ are numbered in accordance with the following
diagram 
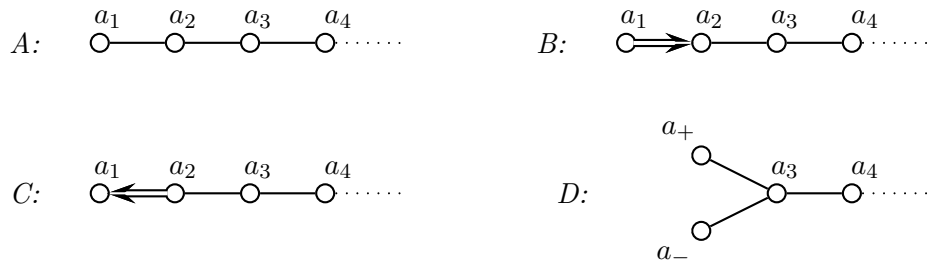
\begin{picture}(1320,90)
 \cnodeput(1,2){A}{}
 \cnodeput(2,2){B}{}
 \cnodeput(3,2){C}{}
 \cnodeput(4,2){D}{}
 \pnode(5,2){E}
 \rput(0,2){\rnode{H}{A:}}
 \ncline {-}{A}{B}\Aput{$a_1 \qquad$}
 \ncline {-}{B}{C}\Aput{$a_2 \qquad$}
 \ncline {-}{C}{D}\Aput{$a_3 \qquad$}
 \ncline[linestyle=dotted]{-}{D}{E}\Aput{$a_4 \qquad$}
  \cnodeput(8,2){A}{}
 \cnodeput(9,2){B}{}
 \cnodeput(10,2){C}{}
 \cnodeput(11,2){D}{}
 \pnode(12,2){E}
 \rput(7,2){\rnode{H}{B:}}
\psline [doubleline=true,
linearc=.5,doublesep=1.5pt]{->}(8.1,2)(8.9,2)
 \ncline[linecolor=white] {-}{A}{B}\Aput{$a_1 \qquad$}
 \ncline {-}{B}{C}\Aput{$a_2 \qquad$}
 \ncline {-}{C}{D}\Aput{$a_3 \qquad$}
 \ncline[linestyle=dotted]{-}{D}{E}\Aput{$a_4 \qquad$}
  \cnodeput(1,0){A}{}
 \cnodeput(2,0){B}{}
 \cnodeput(3,0){C}{}
 \cnodeput(4,0){D}{}
 \pnode(5,0){E}
 \rput(0,0){\rnode{H}{C:}}
\psline [doubleline=true,
linearc=.5,doublesep=1.5pt]{->}(1.9,0)(1.1,0)
 \ncline[linecolor=white] {-}{A}{B}\Aput{$a_1 \qquad$}
 \ncline {-}{B}{C}\Aput{$a_2 \qquad$}
 \ncline {-}{C}{D}\Aput{$a_3 \qquad$}
 \ncline[linestyle=dotted]{-}{D}{E}\Aput{$a_4 \qquad$}
  \cnodeput(9,-0.5){A}{}
 \cnodeput(9,0.5){B}{}
 \cnodeput(10,0){C}{}
 \cnodeput(11,0){D}{}
 \pnode(12,0){E}
 \pnode(8,-0.5){a}{}
 \pnode(8,0.5){b}{}
 \rput(7,0){\rnode{H}{D:}}
 \ncline[linecolor=white] {-}{A}{a}\Aput{$\quad a_- $}
  \ncline[linecolor=white] {-}{b}{B}\Aput{$\quad a_+ $}
 \ncline {-}{A}{C}
 \ncline {-}{B}{C}
 \ncline {-}{C}{D}\Aput{$a_3 \qquad$}
 \ncline[linestyle=dotted]{-}{D}{E}\Aput{$a_4 \qquad$}
\end{picture}

\vspace{1cm}

\noindent
where the $a_j$ are non-negative numbers of which only finitely many
ore non-zero. Let $a_\alpha$ be the rightmost non-zero label. If
$n<\alpha-1$, then the corresponding subordinate
representation of the group $G_n(\bj)=A_n(\bj), B_n(\bj), C_n(\bj),$ 
$D_n(\bj)$ is
zero-dimensional. If, on the other hand, $n\geq\alpha-1$, then the
subordinate representation $G_n(\bj)$ is irreducible and has the numerical
labels ($\ldots, a_{n-1},A_n$) on the corresponding Dynkin diagram.
\end{theorem}
Below we shall give an explicit construction for all the
representations of the categories
$\rm{CK} \bf{A}$, $\rm{CK} \bf{B}$, $\rm{CK} \bf{C}$,  $\rm{CK} \bf{D}$.
 The verification of the
correctness of these constructions and the proof of Theorem \ref{t.5.1} is
deferred until Subsection 5.3.

It should be noted that the holomorphic projective representations of
the category $\rm{CK} \bf{A}$ are linearisable. This will be clear from the
explicit construction.

 We note also that it can be shown that the irreducible holomorphic
projective representations of the category $\rm{CK} \bf{GD}$ are numbered by diagrams
of the form 
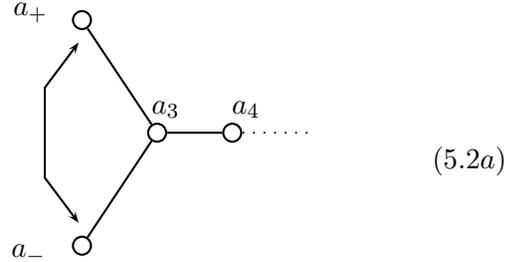
\vspace{-1,8cm}
\noindent
$$
\begin{picture}(1320,90) 
  \cnodeput(7,-1.5){A}{}
 \cnodeput(7,1.5){B}{}
 \cnodeput(8,0){C}{}
 \cnodeput(9,0){D}{}
 \pnode(10,0){E}
 \pnode(6,-1.){a}{}
 \pnode(6,1.){b}{}
\psline []{<->}(6.7,1.2)(6.25,0.6)(6.25,-0.6)(6.7,-1.2)
 \ncline[linecolor=white] {-}{A}{a}\Aput{$\quad a_- $}
  \ncline[linecolor=white] {-}{b}{B}\Aput{$\quad a_+ $}
 \ncline {-}{A}{C}
 \ncline {-}{B}{C}
 \ncline {-}{C}{D}\Aput{$a_3 \qquad$}
 \ncline[linestyle=dotted]{-}{D}{E}\Aput{$a_4 \qquad$}
\end{picture}
\eqno(5.2a)
$$
\vspace{9mm}

\noindent
where the $a_\alpha$ are non-negative integers of which only finitely
many are non-zero and interchanging $a_+$ and $a_-$ does not change the
representation. If $a_+=a_-$, then the restriction of this
representation to the category $\rm{CK} \bf{D}$ is irreducible and has the numerical
labels (5.2a). If $a_+\neq a_-$, then the restriction to $\rm{CK} \bf{D}$
decomposes into a sum of two representations with numerical labels 

\noindent
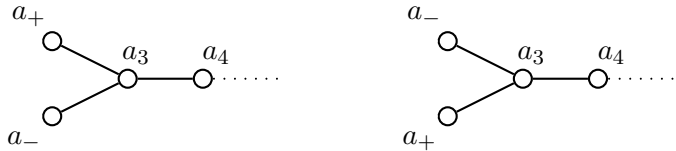
\begin{picture}(1320,50)
  \cnodeput(4,-0.5){A}{}
 \cnodeput(4,0.5){B}{}
 \cnodeput(5,0){C}{}
 \cnodeput(6,0){D}{}
 \pnode(7,0){E}
 \pnode(3,-0.5){a}{}
 \pnode(3,0.5){b}{}
 \ncline[linecolor=white] {-}{A}{a}\Aput{$\quad a_- $}
  \ncline[linecolor=white] {-}{b}{B}\Aput{$\quad a_+ $}
 \ncline {-}{A}{C}
 \ncline {-}{B}{C}
 \ncline {-}{C}{D}\Aput{$a_3 \qquad$}
 \ncline[linestyle=dotted]{-}{D}{E}\Aput{$a_4 \qquad$}
  \cnodeput(9,-0.5){A}{}
 \cnodeput(9,0.5){B}{}
 \cnodeput(10,0){C}{}
 \cnodeput(11,0){D}{}
 \pnode(12,0){E}
 \pnode(8,-0.5){a}{}
 \pnode(8,0.5){b}{}
 \ncline[linecolor=white] {-}{A}{a}\Aput{$\quad a_+ $}
  \ncline[linecolor=white] {-}{b}{B}\Aput{$\quad a_- $}
 \ncline {-}{A}{C}
 \ncline {-}{B}{C}
 \ncline {-}{C}{D}\Aput{$a_3 \qquad$}
 \ncline[linestyle=dotted]{-}{D}{E}\Aput{$a_4 \qquad$}
\end{picture}
\vspace{1cm}

\noindent
{\it The fundamental representations}. Let $\rm{CK} \bf{K}$ be one of the
categories $\rm{CK} \bf{A}$, $\rm{CK} \bf{B}$, $\rm{CK} \bf{C}$,  
$\rm{CK} \bf{D}$. We denote by $(\Pi_\alpha^{\rm{CK} \bf{K}},\pi_\alpha^{\rm{CK} \bf{K}})$ the
irreducible representation of $\rm{CK} \bf{K}$ for which the numerical label
$a_\alpha$ is equal to 1 and the other labels are all zero.

 Strictly speaking, while the classification theorem remains unproved
us do not have the right to talk about numerical labels of irreducible
representations.  Therefore for the moment, by a fundamental
representation $(\Pi_\alpha^{\rm{CK} \bf{K}},\pi_\alpha^{\rm{CK} \bf{K}})$  we mean a representation
of the category $\rm{CK} \bf{K}$ whose subordinate representations are fundamental
representations $\pi_\alpha$ of the groups $K_n(\bj)= A_n(\bj), B_n(\bj), C_n(\bj), D_n(\bj)$
(for $n > a -1$).  We now construct these representations.

Thus we consider the spinor representation (Spin, spin) of the
category $\rm{CK} \bf{GD}$ (see Example \ref{e4.4}) and restrict it to
the 
category $\rm{CK} \bf{D}$. If $P\in
\Mor_{\rm{CK} \bf{D}}(V(\bj),W(\bj))$, then the kernel operator $\spin(P)$ is an
even 
function
(see  Example \ref{e4.4}), therefore $\spin(P)$ takes
even functions in $\Lambda(V_+(\bj))$ to even functions in $\Lambda(W_+(\bj))$,
and odd functions to odd functions.  The representations
$(\Pi_+^{\rm{CK} \bf{K}},\pi_+^{\rm{CK} \bf{K}})$ and 
$(\Pi_-^{\rm{CK} \bf{K}},\pi_-^{\rm{CK} \bf{K}})$
 are subrepresentations in the restriction of (Spin, spin) to
$\rm{CK} \bf{D}$ consisting of the even and odd functions respectively.

Next, we embed the category $\rm{CK} \bf{B}$ in the category 
$\rm{CK} \bf{GD}$. Let $V(\bj)\in
\Ob(\rm{CK} \bf{GD})$, and  let $L(\bj)$ be an one-dimensional complex
space endowed with a non-zero bilinear form.
Then $V(\bj)\oplus L(\bj) \in \Ob(\rm{CK} \bf{GD})$. Suppose further
that 
$P$ is a 
member of
$\Mor_{\rm{CK} \bf{GD}}(V(\bj),W(\bj))\backslash\nully$.  We define the subspace $Q
\subset(V\oplus L)\oplus(W\oplus L)$ as the set of all vectors of the
form $((v,l),(w,l))$, where $(v,w) \in P$ and $l\in L$. Then $Q= Q(P)
\in \Mor_{\rm{CK} \bf{GD}}(V(\bj),W(\bj))$.  The restriction of (Spin,
spin) 
to $\rm{CK} \bf{B}$
decomposes into a sum of two equivalent representations of the  form
$(\Pi_1^{\rm{CK} \bf{B}}, \pi_1^{\rm{CK} \bf{B}})$, one of which is
realized 
in even functions,
and the other in odd functions.  (A more explicit and more
interesting construction of the representation
$(\Pi_1^{\rm{CK} \bf{B}}, \pi_1^{\rm{CK} \bf{B}})$ is given below.)

 Let $V(\bj), W(\bj) \in\Ob(\rm{CK} \bf{GA})$ and let $P\subset
V(\bj)\oplus W(\bj)$ be a non-zero morphism of the category $\rm{CK}
\bf{GA}$
 with $\dim P= s$.
We see that the operator $\Lambda(P)$ takes the subspace
$\Lambda^kV(\bj)\subset \Lambda(V(\bj))$ to $\Lambda^{k-\dim(V(\bj))+s}(W(\bj))$. We
restrict the representation $(\Lambda,\lambda)$ to the category
$\rm{CK} \bf{K}=\rm{CK} \bf{B}, \rm{CK} \bf{C}, \rm{CK} \bf{D}$.  
Let $V(\bj), W(\bj) \in \Ob(\rm{CK} \bf{K})$, and let $P \in 
\Mor_{\rm{CK} \bf{K}}(V(\bj),
W(\bj))\backslash\nully$.  Then
\beq
\dim P = {1\over2}(\dim V(\bj) + \dim W(\bj)).
\eeq
Thus the restriction of $(\Lambda,\lambda)$ to $\rm{CK} \bf{K}$
decomposes 
into a
countable direct sum
\beq
(\Lambda,\lambda)|_{\rm{CK} \bf{K}} =\bigoplus\limits_{j=-\infty}^{+\infty}
\left(L_j^{\rm{CK} \bf{K}}, l_j\right),
\eeq
where
\beq
L_j^{\rm{CK} \bf{K}}(V(\bj)) =\Lambda^{[{1\over2}\dim V(\bj)]-j+1}(V(\bj)).
\eeq
It is easy to see that:

1. $(\Pi_j^{\rm{CK} \bf{B}}, \pi_j^{\rm{CK} \bf{B}})= (L_j^{\rm{CK} \bf{B}},l_j)$ for $j\geq2$;

2. $(\Pi_j^{\rm{CK} \bf{D}}, \pi_j^{\rm{CK} \bf{D}})=(L_j^{\rm{CK} \bf{D}},l_j)$ for $j\geq 3$;

3. $(\Pi_j^{\rm{CK} \bf{C}}, \pi_j^{\rm{CK} \bf{C}})$ is the factor representation $(L_j^{\rm{CK} \bf{C}},l_j)/Q
(L_{j+2}^{\rm{CK} \bf{C}},l_{j+2})$, where $Q$ is the embedding of the
representation $(L_{j+2}^{\rm{CK} \bf{C}},l_{j+2}$ in $(L_j^{\rm{CK} \bf{C}},l_j)$, described as
follows. We denote by $q$ an element of $\Lambda^2V(\bj)$ which is invariant
with respect to $Sp(V(\bj))$. Then $Qf :=qf$.

Finally, the representation $(\Pi_1^{\rm{CK} \bf{A}}, \pi_1^{\rm{CK} \bf{A}})$ is the
self-representation of the category $\rm{CK} \bf{A}$  and
\beq
\left(\Pi_j^{\rm{CK} \bf{A}}, \pi_j^{\rm{CK} \bf{A}}\right) = \Lambda^j
\left(\Pi_1^{\rm{CK} \bf{A}}, \pi_1^{\rm{CK} \bf{A}}\right).
\eeq

\noindent
{\it Construction of the remaining representations}. 

Let 
$\rm{CK} \bf{K}=\rm{CK} \bf{A}, \rm{CK} \bf{B}, \rm{CK} \bf{C},$ or  $\rm{CK} \bf{D}$.  We now construct an irreducible representation $T = (T, r)$ of
the category $\rm{CK} \bf{K}$ with set of numerical labels $\{a_\alpha\}$. For
this we consider the tensor product

\beq
(S,\sigma):=\bigotimes\limits_\alpha
\left(\Pi_\alpha^{\rm{CK} \bf{K}}, \pi_\alpha^{\rm{CK} \bf{K}}\right)^{\otimes a_\alpha}.
\eeq
In each space $\Pi_\alpha^{\rm{CK} \bf{K}}(V(\bj))$ we consider the
vector 
$h_\alpha(V(\bj))$ of
highest weight with respect to $\Aut(V(\bj))$. 

Let $T(V(\bj))$ be the 
cyclic span
of the vector $\otimes_\alpha h_\alpha(V(\bj))^{\otimes a_\alpha} \in
S(V(\bj))$ under the action of the group $\Aut(V(\bj))$.  The set of subspaces
$T(V(\bj)) \subset S(V(\bj))$ defines a subrepresentation in $(S,\sigma)$ which
has the numerical labels $\{a_\alpha\}$.

\noindent
{\it   The category $\rm{CK} \bf{GA}$}.
\begin{theorem} {\rm (see [\cite{SRCKOrtC}])}

(a) The holomorphic projective representations of
the  category $\rm{CK} \bf{GA}$ are completely reducible.

(b) The irreducible holomorphic projective representations of the
category $\rm{CK} \bf{GA}$ are numbered by diagrams of the form

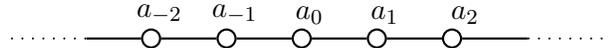
\begin{picture}(1320,30)
\cnodeput(4,0){A}{}
 \cnodeput(5,0){B}{}
 \cnodeput(6,0){C}{}
 \cnodeput(7,0){D}{}
 \cnodeput(8,0){E}{}
 \pnode(9,0){H}{}
 \pnode(2,0){a}
 \pnode(3,0){b}
 \pnode(10,0){c}
 \ncline []{-}{A}{b}
 \ncline [linestyle=dotted]{-}{b}{a}
 \ncline {-}{A}{B}\Aput{$a_{-2} \qquad$}
 \ncline {-}{B}{C}\Aput{$a_{-1} \qquad$}
 \ncline {-}{C}{D}\Aput{$a_{0} \qquad$}
 \ncline {-}{D}{E}\Aput{$a_{1} \qquad$}
 \ncline[] {-}{E}{H}\Aput{$a_{2} \qquad$}
 \ncline [linestyle=dotted]{-}{H}{c}
\end{picture}
\vspace{5mm}

\noindent
where the $a_j$ are non-negative integers of which only finitely many
are  nonzero.  Furthermore, diagrams differing from each other by a
shift correspond to the same representation. Let $a_\alpha$ be the
leftmost non-zero numerical label, and $a_\beta$ the rightmost. Then
the subordinate representation of the group $A_n(\bj)\simeq SL(n+1;\bj,\Cb)$ is
the simple direct sum of the irreducible representations $A_n(\bj)$ with
numerical labels
$$
\begin{picture}(1320,30)

 \cnodeput(3.5,0){A}{}
 \cnodeput(4.5,0){B}{}
 \cnodeput(5.5,0){C}{}
 \pnode(6.5,0){D}{}
 \pnode(7.5,0){E}{}
 \cnodeput(8.5,0){H}{}
 \pnode(1.5,0){a}
 \pnode(2.5,0){b}
 \pnode(9.5,0){c}
 \ncline [linecolor=white]{-}{A}{b}
 \ncline [linestyle=dotted,linecolor=white]{-}{b}{a}
 \ncline {-}{A}{B}\Aput{$a_{\gamma} \qquad$}
 \ncline {-}{B}{C}\Aput{$a_{\gamma+1} \qquad$}
 \ncline {-}{C}{D}
 \ncline {-}{H}{E}
 \ncline[linestyle=dotted] {-}{E}{D}%
 \ncline
[linestyle=dotted,linecolor=white]{-}{H}{c}\Aput{$a_{\gamma+n-1}
\qquad$}
\end{picture}
\eqno(2.2b)
$$

\vspace{5mm}
\noindent
where $\gamma\leq \alpha+1$ and $\gamma+n-1\geq \beta-1$.
\end{theorem}

\vspace{0.5cm}

Note that the fundamental representation in Example \ref{e4.5} has the
numerical labels

\begin{picture}(1320,40)

 \cnodeput(4.5,0){B}{}
 \cnodeput(5.5,0){C}{}
 \cnodeput(6.5,0){D}{}
 \cnodeput(7.5,0){E}{}
 \pnode(8.5,0){H}{}
 \pnode(2.5,0){a}
 \pnode(3.5,0){b}
 \pnode(9.5,0){c}
 \ncline []{-}{B}{b}
 \ncline [linestyle=dotted]{-}{b}{a}
 \ncline {-}{B}{C}\Aput{$0 \qquad$}
 \ncline {-}{C}{D}\Aput{$0 \qquad$}
 \ncline {-}{D}{E}\Aput{$1 \qquad$}
 \ncline[] {-}{E}{H}\Aput{$0 \qquad$}
 \ncline [linestyle=dotted]{-}{H}{c}

\end{picture}

\vspace{8mm}

\noindent
all the irreducible representations of $\rm{CK} \bf{GA}$ being realized in the tensor
powers of the fundamental representation.

Note also that if $\mu_\gamma$ is the representation with the numerical
labels (5.2b) and if $R\in\End(\Cb^{n+1}(\bj))\backslash\nully$, then
$R\mu_\gamma\in \mu_{\gamma+\dim R-(n+1)}$.

\vspace{0.3cm}
\noindent
{\it Spinor representation of the category $\rm{CK} \bf{B}$}. 
Let $V(\bj) \in \Ob(\rm{CK} \bf{GD})$.
Let $L(\bj)\simeq \Cb(\bj)$ be a one-dimensional complex linear space endowed
with the bilinear form $\{z, u\} = zu$. Let $\tilde{V}(\bj) = V(\bj) \oplus
L(\bj)$. 
Clearly,
$\tilde{V}(\bj) \in \Ob(\rm{CK} \bf{B})$.  
For each $\tilde{v}=(v,s)\in\tilde{V}(\bj)$ we
define the operator $\hat{a}(\tilde{v})$ in $\Lambda(V_+(\bj))$ (here, as
before, $V_+(\bj)$ is a maximal isotropic subspace of $(V(\bj))$ by the formula
\beq
\hat{a}(\tilde{v}) f=\hat{a}(v) f+sf.
\eeq
\begin{theorem} (a) For each $P\in \Mor_{\rm{CK} \bf{B}}(\tilde{V}(\bj),
\tilde{W}(\bj))$ there
exists a unique (to within proportionality) non-zero operator
\beq
\spin_{\rm{CK} \bf{B}}(P): \Lambda(V_+(\bj)) \rightarrow \Lambda(W_+(\bj))
\eeq
such that
\beq
\spin_{\rm{CK} \bf{B}}(P)\hat{a}(\tilde{v}) = 
\hat{a}(\tilde{w})\spin_{\rm{CK} \bf{B}}(P)
\eeq
for all $(\tilde{v}, \tilde{w})\in P$.

(b) The map $P \mapsto \spin_{\rm{CK} \bf{B}}(P)$, $\null\mapsto 0$ is a projective
representation of the category $\rm{CK} \bf{B}$.
\end{theorem}
For an arbitrary function $f\in\Lambda_m$ we define the expression
$\exp(f)=\sum f^j/j!$.  If at least one of the functions $f_1$, $f_2$
is even, then $f_1f_2=f_2f_1$ and therefore
\beq
exp(f_1+f_2)=\exp(f_1)\exp(f_2)=\exp(f_2)\exp(f_1).
\eeq
Suppose further that $\nu_1,\nu_2,\ldots$ are linear expressions in
the variables $\xi$. Then it is easily seen that
\be \label{5.3}
(1+\nu_1)(1+\nu_2)=\exp(\nu_1+\nu_2+\nu_1\nu_2)
=\exp(\nu_1\nu_2)\exp(\nu_1+\nu_2).
\ee
Hence we can prove the following equality by induction:
\be \label{5.4}
(1+\nu_1)(1+\nu_2) \ldots (1+\nu_k) = \exp \left(\sum_{i<j}
\nu_i\nu_j\right) \exp \left(\sum \nu_i\right)
\ee
It also follows easily from (\ref{5.3}) that
\be \label{5.5}
\exp(\nu_1\nu_2) =(1+\nu_1)(1+\nu_2) (1-\nu_1-\nu_2).
\ee
\noindent
{\it The generalised Beresin operator} $\Lambda_n\rightarrow
\Lambda_m$.  This is the operator with kernel of the form
\beq
K(\xi,\bar{\eta}) =\lambda \prod\limits_{i=1}^s l_i(\xi,\bar{\eta})
\exp\left\{ {1\over2} (\xi\,\bar{\eta}\,1)
\left(\ba{ccc}
K& L& p^t\\
-L^t& M& q^t\\
p& q& 0\\ \ea \right)
\left(\ba{c} \xi\\
\bar{\eta}\\
1\\ \ea \right) \right\},
\eeq
where the $l_i(\xi,\bar{\eta})$ are linear expressions in $\xi$,
$\bar{\eta}$, $K=-K^t$, $M=M^t$, and $p$  and $q$ are row
matrices.

An equivalent definition is as follows. A generalised Berezin operator
is an  operator with kernel of the form
\beq
\prod\limits_{j=1}^N \left(a_j +\sum \alpha_{jk}\xi_k +\sum
\beta_{js}\bar{\eta}_s\right).
\eeq
The equivalence of the two definitions can be seen from the equalities
(\ref{5.4}) and (\ref{5.5}).

It can be shown that the operators $\spin_{\rm{CK} \bf{B}}(P)$ are generalised Berezin
operators and, conversely, any generalised Berezin operator has the
form $\spin_{\rm{CK} \bf{B}}(P)$.

\subsection{Ordered Cayley-Klein categories}

In this Subsection we prove the classification theorems of Subsection 5.2.

\noindent
Let $\Sigma$ be a partially ordered set such that for any
$\sigma_1,\sigma_2\subset\Sigma$ there exists $\sigma_3\in \Sigma$ such
that $\sigma_3>\sigma_1$ and $\sigma_3>\sigma_2$. in this article the
set $\Sigma$ is always $\Zb_+$.  Let $\rm{CK} \bf{K}$ be a category 
whose objects
$V_\sigma(\bj)$ are numbered by an element of the set $\Sigma$.  Suppose
that for any two elements $\sigma,\tau\in \Sigma$ such that
$\sigma<\tau$ there are fixed morphisms $\lambda_{\sigma\tau}:
V_\sigma(\bj)\rightarrow V_\tau(\bj)$ and $\mu_{\tau\sigma}:V_\tau(\bj)\rightarrow
V_\sigma(\bj)$ such that
\be \label{5.6}
\mu_{\tau\sigma} \lambda_{\tau\sigma}=1.
\ee
Suppose that for any $\sigma<\sigma'<\sigma''$ we have
\be \label{5.7}
\lambda_{\sigma'\sigma''} \lambda_{\sigma\sigma'}
=\lambda_{\sigma\sigma''}, \qquad \mu_{\sigma'\sigma}
\mu_{\sigma''\sigma'} =\mu_{\sigma''\sigma}.
\ee
 We call such categories {\it purely ordered}. Categories that are
equivalent to purely ordered categories are called {\it ordered Cayley-Klein 
categories}.

\vspace{0.3cm}
\noindent
{\it Examples of ordered Cayley-Klein categories}. We discuss the
categories 
$\rm{CK} \bf{A}, \rm{CK} \bf{B},
\rm{CK} \bf{C}, \rm{CK} \bf{D}$ in detail.  The simplest is the category 
$\rm{CK} \bf{K}=\rm{CK} \bf{A}$ and the most
complex is $\rm{CK} \bf{K}= \rm{CK} \bf{D}$.

\noindent
{\it The category $\rm{CK} \bf{A}$}. We consider the Cayley-Klein category whose objects are the
spaces $\Cb^n(\bj)$ and whose morphisms are linear operators. This category
is equivalent to the category $\rm{CK} \bf{A}$. Let $m < n$ and define
$\lambda_{mn}: \Cb^m(\bj)\rightarrow \Cb^n(\bj)$, $\mu_{mn}: \Cb^n(\bj)\rightarrow
\Cb^m(\bj)$ by the formulae:
\beq
\lambda_{mn}(x_1,\ldots,x_m) &\!\!\!=\!\!\!&
(x_1,\ldots,x_m,0,\ldots,0),\\
\mu_{mn}(x_1,\ldots,x_n)  &\!\!\!=\!\!\!&
(x_1,\ldots,x_m).
\eeq
\noindent
{\it The category $\rm{CK} \bf{D}$.} We denote a $2n$-dimensional object of the
category $\rm{CK} \bf{D}$ by $V_n(\bj)$. Let $V(\bj)$ be a two-dimensional complex space
endowed with a non-degenerate symmetric bilinear form. Then $V(\bj)$
contains two isotropic lines which we denote
by $l$ and $l'$. We identify $V_{n+1}(\bj)$ with $V_n(\bj)\oplus V(\bj)$. Let $S$ be
the graph of the natural embedding $V_n(\bj) \rightarrow  V_{n+1}(\bj)$. Consider
the following maximal isotropic subspaces $H(\bj)$ and $H'(\bj)$ of $V_n(\bj) \oplus
V_{n+1}(\bj)$:
\beq
H=S\oplus l, \qquad H'=S\oplus l'.
\eeq
It is easily seen that $H\cap H'=S$, therefore (see Lemma \ref{l5.1}) $H$
and $H'$ lie different connected components of the Grassmannian
${\mathcal Gr}(V_n(\bj)\oplus V_{n+1}(\bj)$). Note that the
sets $\Mor_{\rm{CK} \bf{D}}(V_n(\bj), V_{n+1}(\bj))\backslash\nully$ and 
$\Mor_{\rm{CK} \bf{D}}(V_{n+1}(\bj),
V_n(\bj))\backslash\nully$ are the connected components of the manifold
${\mathcal Gr}(V_n(\bj)\oplus V_{n+1}(\bj))$ and it is easily seen that these components
are distinct.  (In fact, $V_n^-(\bj)\oplus V_{n+1`}^+(\bj)\in 
\Mor_{\rm{CK} \bf{D}}(V_n(\bj),
V_{n+1}(\bj))$, while $V_{n+1}^-(\bj)\oplus V_n^+(\bj)\in 
\Mor_{\rm{CK} \bf{D}}(V_{n+1}(\bj), V_n(\bj))$,
the intersection of $V_n^-(\bj)\oplus V_{n+1}^+(\bj)$ and  $V_{n+1}^-(\bj)\oplus
V_n^+(\bj)$ is  zero, and the dimensions of the subspaces are odd.) But now
$\lambda_{n,n+1}$ is whichever of the two subspaces $H, H'$ that lies
in $\Mor_{\rm{CK} \bf{D}}(V_n(\bj),V_{n+1}(\bj))$, while $\mu_{n+1,n}$ is the one
which lies in $\Mor_{\rm{CK} \bf{D}}(V_{n+1}(\bj), V_n(\bj))$. Finally,
\beq
\lambda_{ij} =\lambda_{j-1, j}\lambda_{j-2, j-1} \ldots \lambda_{i, i+1},
\qquad \mu_{ji}=\mu_{i+1, i} \mu_{i+1, i+2} \ldots \mu_{j, j-1}.
\eeq
\noindent
{\it The categories $\rm{CK} \bf{B}$ and $\rm{CK} \bf{C}$}. Here the situation is similar to that of
category $\rm{CK} \bf{D}$, except that we have a freer choice for the lines $l$,
$l'$ and the morphisms $\lambda_{n,n+1}$, $\mu_{n+1,n}$: the natural
requirement is that $l'\neq l$ and that $\lambda_{n,n+1} =S\oplus l$
and $\mu_{n+1,n}=S\oplus l'$.

\vspace{0.3cm}
\noindent
{\it  Properties of ordered Cayley-Klein categories}.
Let $\sigma< \tau$. We define the element $\vartheta_\tau^\sigma$ in
in $\End(V_\tau(\bj))$ by the formula
\be \label{5.8}
\vartheta_\tau^\sigma =\lambda_{\sigma\tau} \mu_{\tau\sigma}.
\ee
It is easily seen that
\be \label{5.9}
(\vartheta_\tau^\sigma)^2=\vartheta_\tau^\sigma, \qquad
\mu_{\tau\sigma} \vartheta_\tau^\sigma=\mu_{\tau\sigma}, \qquad
\vartheta_\tau^\sigma \lambda_{\sigma\tau} =\lambda_{\sigma\tau}.
\ee
 If $\sigma'<\sigma$, then
\be \label{5.10}
\vartheta_\tau^\sigma\vartheta_\tau^{\sigma'} =\vartheta_\tau^{\sigma'}
\vartheta_\tau^\sigma =\vartheta_\tau^{\sigma'}.
\ee
Let us verify the latter equality. Indeed,
\beq
\vartheta_\tau^{\sigma'} \vartheta_\tau^\sigma  &\!\!\!=\!\!\!&
\lambda_{\sigma'\tau} \mu_{\tau\sigma'}
\lambda_{\sigma\tau}\mu_{\tau\sigma} = \lambda_{\sigma'\tau}
\mu_{\sigma\sigma'} \mu_{\tau\sigma} \lambda_{\sigma\tau}
\mu_{\tau\sigma} =\\
 &\!\!\!=\!\!\!& \lambda_{\sigma'\tau} \mu_{\sigma\sigma'}
\mu_{\tau\sigma} = \lambda_{\sigma'\tau} \mu_{\tau\sigma'}
=\vartheta_\tau^{\sigma'},
\eeq
as required.

\begin{lemma} \label{l5.3}
Let $\sigma'<\sigma$ and $\tau'<\tau$. Let $P \in
\Mor(V_{\sigma'(\bj)}, V_{\tau'(\bj)})$.  Then there exists  
$Q\in \Mor(V_\sigma(\bj),
V_\tau(\bj))$ such that
\be \label{5.11}
P=\mu_{\tau\tau'} Q\lambda_{\sigma'\sigma}.
\ee
\end{lemma}

\noindent
{\bf Proof} We have $Q=\lambda_{\tau'\tau} P \mu_{\sigma\sigma'}$.
\hfill{$\Box$}

\begin{lemma} \label{l5.4}
Let $(T,\tau)$ be a (projective) representation
of the ordered category $\rm{CK} \bf{K}$. Let $T(V_\kappa(\bj)) = 0$. Then
$T(V_\sigma(\bj))=0$ for all $\sigma<\kappa$.
\end{lemma}
\noindent
{\bf Proof} Since $T(V_\kappa(\bj))= 0$, it follows that $\tau(Q)=0$
for all $Q\in\End(V_\kappa(\bj))$. In view  of (\ref{5.11}), we have
$\tau(P)=0$ for all $P\in\End(V_\sigma(\bj))$. In
particular, $\tau(1_{V_\sigma(\bj)})=0$. Hence $T(V_\sigma(\bj))=0$.
\hfill{$\Box$}

\begin{lemma} \label{l5.5}
Let $(T,\tau)$ be an irreducible representation.
 Then all the subordinated representations of the Cayley-Klein semigroup
$\End(V(\bj))$ are irreducible.
\end{lemma}

\noindent
{\bf Proof}  Suppose that $T(V(\bj))$ contains an $\End(V(\bj))$-invariant
subrepresentation. Then the cyclic span of A is a non-trivial
subrepresentation of $T$. \hfill{$\Box$}

We note that the converse is false in general. It is, however, true for
ordered Cayley-Klein categories.

\begin{lemma}  \label{l5.6}
Let $(T,\tau)$ be a (projective) representation of the
ordered Cayley-Klein  category $\rm{CK} \bf{K}$. Then the following statements are equivalent:

(a) the representation $(T, \tau)$ is irreducible;

(b) the subordinate, representations of all the semigroups
$\End(V_\tau(\bj))$ are irreducible.
\end{lemma}

\noindent
{\bf Proof} We merely need to prove the implication $(b)\Rightarrow
(a)$. Thus suppose that (b) holds. Let $M$ be
a subrepresentation of $T = (T,\tau)$, and let $N=T/M$. Then for any
$V_\sigma(\bj)$ we have either $M(V_\sigma(\bj)) = 0$ or $N(V_\sigma(\bj)) = 0$. It
now follows immediately from Lemma \ref{l5.5} that either
$M(V_\sigma(\bj))=0$ for all $\sigma$ or $N(V_\sigma(\bj))=0$ for all $\sigma$.
\hfill{$\Box$}

\begin{corollary} \label{c5.1}
The fundamental representations of the Cayley-Klein
categories $\rm{CK} \bf{A}, \rm{CK} \bf{B}, \rm{CK} \bf{C}, \rm{CK} \bf{D}$ are irreducible. \hfill{$\Box$}
\end{corollary}

\noindent
{\it The lowering functor}. Let $T =(T,\tau)$ be a representation of
the category $\rm{CK} \bf{K}$, and let $\alpha<\beta$.

\begin{lemma} \label{l5.7}
Let $P\in \End(V_\alpha(\bj))$. Let
$U_\alpha^\beta(P)=\lambda_{\alpha\beta} P\mu_{\beta\alpha}$.  Then the
map $U_\alpha^\beta$ is an embedding of the semigroup $\End(V_\alpha(\bj))$
into $\End(V_\beta(\bj))$.
\end{lemma}

\noindent
{\bf Proof} We have
\beq
U_\alpha^\beta(P) U_\alpha^\beta(Q) = \lambda_{\alpha\beta} P
\mu_{\beta\alpha}\lambda_{\alpha\beta}Q\mu_{\beta\alpha} = \lambda_{\alpha\beta}
PQ\mu_{\beta\alpha} =U_\alpha^\beta(PQ).
\eeq
On the other hand, $P =\mu_{\beta\alpha}
U_\alpha^\beta(P)\lambda_{\alpha\beta}$, whence it follows that
$U_\alpha^\beta$ is injective. Note that
$U_\alpha^\beta(1)=\vartheta_\beta^\alpha$.  Thus for $\alpha<\beta$
the semigroup $\End(V_\sigma(\bj))$ is embedded in $\End(V_\beta(\bj))$ (one must
bear in mind that $U_\alpha^\beta(\Aut(V_\alpha(\bj)))\not\subset
\Aut(V_\beta(\bj)))$. \hfill{$\Box$}

\begin{proposition} \label{p5.1} 
Let $(T,\tau)$ be a projective representation of
the ordered Cayley-Klein category $\rm{CK} \bf{K}$, and let $\alpha<\beta$. Then the subspace
$\im~~\tau(\vartheta_\beta^\alpha)$ is invariant with respect to
the operators $\tau(U_\alpha^\beta(P))$.  The representation $\tau$ of
the semigroup $\End(V_\alpha(\bj))$ in $T(V_\alpha(\bj))$ is equivalent to the
representation $\tau\circ U_\alpha^\beta$ in
$\im~~\tau(\vartheta_\beta^\alpha)$.
\end{proposition}

\noindent
{\bf Proof} We claim that $\im~~\tau(U_\alpha^\beta(P))\subset
\im~~ \tau(\vartheta_\beta^\alpha)$. In  fact,
\beq
U_\alpha^\beta(P) =\lambda_{\alpha\beta} P \mu_{\beta\alpha}
=\lambda_{\alpha\beta} \mu_{\beta\alpha}\lambda_{\alpha\beta} P
\mu_{\beta\alpha} =\vartheta_\beta^\alpha \lambda_{\alpha\beta} P
\mu_{\beta\alpha}.
\eeq
Furthermore, as is easily seen, the operators
$\tau(\lambda_{\alpha\beta}): T(V_\alpha(\bj))\rightarrow
\im~~\tau(\vartheta_\beta^\alpha)$ and
$\mu_{\beta\alpha}:\im~~\tau(\vartheta_\beta^\alpha) \rightarrow
T(V_\alpha(\bj))$ are $\End(V_\alpha(\bj))$-intertwining operators which are
mutual inverses.  This completes the proof. \hfill{$\Box$}


Let $\alpha<\beta$. We define the {\it lowering functor}
$F_\alpha^\beta$ which associates with each representation $\tau$ of
the semigroup $\End(V_\beta(\bj))$ the representation $\tau\circ
U_\alpha^\beta$ of t semigroup $\End(V_\alpha(\bj))$ in
$\im~~\tau(\vartheta_\beta^\alpha)$.  We omit the trivial verification
of the fact that $F_\alpha^\beta$ is indeed a functor.  It is also
easily seen that $F_\alpha^\beta F_\beta^\gamma = F_\alpha^\gamma$.

\noindent
We now discuss what the representations of the semigroups
$\End_{\rm{CK} \bf{K}}(V_\alpha(\bj))$ look like in the cases $\rm{CK} \bf{K}=\rm{CK} \bf{A}, \rm{CK} \bf{B}, \rm{CK} \bf{C}, \rm{CK} \bf{D}$. We shall
also discuss what the lowering functor looks like in these cases.

\vspace{0.3cm}

{\it Representations of the semigroups $\End_K(V(\bj))$ in the cases
$\rm{CK} \bf{K}=\rm{CK} \bf{A}, \rm{CK} \bf{B}, \rm{CK} \bf{C}, \rm{CK} \bf{D}$}.

\begin{lemma} \label{l5.8}
Let $\rm{CK} \bf{K}=\rm{CK} \bf{A}, \rm{CK} \bf{B}, \rm{CK} \bf{C}, \rm{CK} \bf{D}$. 
Then the group \\$\Aut_K(V(\bj))$ is
dense in $\End_K(V(\bj))$. \hfill{$\Box$}
\end{lemma}

\noindent
We choose in each of the Cayley-Klein categories 
$\rm{CK} \bf{A}, \rm{CK} \bf{B}, \rm{CK} \bf{C}, \rm{CK} \bf{D}$ one object of each
dimension. We denote by $V_n(\bj)$ an object of dimension $n +1$ in the case
of the category $\rm{CK} \bf{A}$, an object of dimension $2n$ in the
cases of 
$\rm{CK} \bf{C}$
and $\rm{CK} \bf{D}$, and an object of dimension $2n+1$ in the case of 
$\rm{CK} \bf{B}$.

\begin{lemma}  \label{l5.9}
The semigroup $\End(V_n(\bj))$ is generated by the group
$\Aut(V_n(\bj))$ and  the linear relation  $\vartheta_n^{n-1}$.
\end{lemma}

\noindent
{\bf Proof} This is done by a straightforward enumeration of cases.
\hfill{$\Box$}


Suppose now that $\tau$ is a holomorphic projective representation of
the semigroup $\End(V_n(\bj))$. Since $\Aut(V_n(\bj))$ is dense in $\End(V_n(\bj))$,
it follows that $\tau$ and the restriction of $\tau$ to $\Aut(V_n(\bj))$
have the same subrepresentations.  In particular, if $\tau$ is
irreducible, then so is its restriction to $\Aut(V_n(\bj))$. Furthermore, it
is well known that projective representations of semisimple groups
and, in particular, the classical groups, can be linearised on their
simply connected covering groups \cite{Neretin}.   Thus
the words ``projective representation of a simply-\-connected semisimple
group $G$'' and ``linear representation of $G$'' essentially mean the same
thing. In particular, we see straight away that the representations of
the semigroup $\End(V_n(\bj))$ are completely reducible.

Let us consider an irreducible representation $\pi=\pi[\ldots,
a_{n-1},a_n]$ of the classical group $\Aut(V_n(\bj))= A_n(\bj), B_n(\bj), C_n(\bj), D_n(\bj)$
with numerical labels $(\ldots, a_{n-1},a_n)$.  We now ask how this
representation can be extended to a projective representation of the
semigroup $\End(V_n(\bj))$. Let $g_j\in \Aut(V_n(\bj))$ and $g_j \rightarrow
\vartheta_n^{n-1}$.  Then for some $\lambda_j\in\Cb^*$
we have $\lambda_j\pi(g_j)=\lambda_j\pi(g_j)\rightarrow
\tau(\vartheta_n^{n-1})$. Therefore either $\tau(\vartheta_n^{n-1})=0$,
or $\tau(\vartheta_n^{n-1})$ is uniquely defined to within
proportionality.  Thus (see Lemma \ref{l5.9}) there exist at most two
extensions of the representation $\pi[\ldots, a_{n-1},a_n]$ to  a
holomorphic representation of $\End(V_n(\bj))$:

(a) {\it The zero extension} $\tau=\pi_0[\ldots, a_{n-1}, a_n]$. Here
$\tau(\vartheta_n^{n-1})=0$ and hence $\tau(P)$ is identically zero on
$\End(V_n(\bj))\backslash \Aut(V_n(\bj))$.

(b) {\it The maximal extension} ~~$\tau=\pi_\max[\ldots, a_{n-1},a_n]$.
Here $\tau(\vartheta_n^{n-1})\neq0$.

The existence of the zero  extension is obvious.  We now explain why
the maximal extension exists.

\noindent
{\it The lowering functor in the cases of the
categories $\rm{CK} \bf{A}, \rm{CK} \bf{B}, \rm{CK} \bf{C}, \rm{CK} \bf{D}$.}

 As before $\rm{CK} \bf{K}$ stands
for one of the categories $\rm{CK} \bf{A}, \rm{CK} \bf{B}, \rm{CK} \bf{C}, \rm{CK} \bf{D}$.

\begin{lemma} \label{l5.10}
Let $(\Pi_\alpha,\pi_\alpha)$ be a fundamental
representation of the Cayley-Klein category $\rm{CK} \bf{K}$.  
Let $h_\alpha(V_n(\bj))$ be a highest
weight vector in $\Pi_\alpha(V_n(\bj))$. Then
\noindent
\beq
\pi_\alpha(\vartheta_n^{n-1}) h_\alpha(V_n(\bj))  &\!\!\!=\!\!\!&
s_nh_\alpha(V_n(\bj)),\;\; 
\pi_\alpha(\lambda_{n-1,n}) h_\alpha(V_{n-1}(\bj)) =t_nh_\alpha(V_n(\bj)),\\
\pi_\alpha(\mu_{n,n-1}) h_\alpha(V_n(\bj))  &\!\!\!=\!\!\!&
p_nh_\alpha(V_{n-1}(\bj)),
\eeq
where $t_n\in\Cb^*$ and $s_n,p_n\in\Cb$ are non-zero if
$\Pi_\alpha(V_{n-1}(\bj))\neq0$.
\end{lemma}

\noindent
{\bf Proof} This is done by enumeration of cases.  \hfill{$\Box$}


Suppose now that $a_\alpha\in\Zb_+$ and $a_\alpha=0$ for sufficiently
large $\alpha$.  Let $q$ be the rightmost non-zero label. Let
\beq
(S,\sigma):= \bigotimes\limits_\alpha (\Pi_\alpha,\pi_\alpha)^{\otimes
a_\alpha}
\eeq
and let $h(V_n(\bj))$ be a highest weight vector in $S(V_n(\bj))$, that is,
\beq
h(V_n)=\bigotimes\limits_\alpha h_\alpha(V_n)^{\otimes a_\alpha}.
\eeq
Note that the space $\Pi_\alpha(V_n(\bj))$ is non-zero if and only if $n\geq
q-1$ (here we need to look at all the fundamental representations and
verify that the above assertion is true for these representations).
Therefore $S(V_n(\bj))$ is non-zero if and only if $n\geq q-1$, and the
vector $h(V_n(\bj))$ is also non-zero for these same $n$.

In each space $S(V_n(\bj))$ we choose the cyclic span $T(V_n(\bj))$ of the vector
$h(V_n(\bj))$ under the action of the group $\Aut(V_n(\bj))$.

\begin{lemma} \label{l5.11}
The set of subspaces $T(V_n(\bj))\subset S(V_n(\bj))$ defines an
irreducible subrepresentation in $(S, \sigma)$.
\end{lemma}

\noindent
{\bf Proof}  We take the cyclic span $H$ of some vector $h(V_n(\bj))$ under
the action of the category $\rm{CK} \bf{K}$. In view of the previous
lemma, 
this
cyclic span contains all the vectors $h(V_n(\bj))$,  this being true for
all  $n$; therefore the $K$-cyclic spans of all the
vectors $h(V_n(\bj))$ coincide. Consequently, for each $k$ the space
$H(V_k(\bj))$ is the $\End(V_k(\bj))$-cyclic span of the vector $h(V_k(\bj))$. But
the group $\Aut(V_k(\bj))$ is dense in $\End(V_k(\bj))$ and therefore
$H(V_k(\bj))=T(V_k(\bj))$ for all $k$.  Furthermore, the cyclic span of
a vector is a subrepresentation, which means that the set of subspaces
$T(V_k(\bj))$ is in fact a subrepresentation in $S$. Its irreducibility
follows from Lemma \ref{l5.4}.  \hfill{}$\Box$


We now consider which representation $\tau_k$ of $\End(V_k(\bj))$ is realized
in $T(V_k(\bj))$. It can be either $\pi_0[\ldots, a_{k-1}, a_k]$, or
$\pi_\max[\ldots, a_{k-1}, a_k]$. Note that the conditions $T(V_{k-1}(\bj))
= 0$ and $\tau(\vartheta_k^{k-1})= 0$ are equivalent (see Proposition
\ref{p5.1}).  Therefore (see above mentioned definition of maximal and zero
extension) we have
\bea \label{5.12}
\tau_{q-1}  &\!\!\!=\!\!\!& \pi_0[\ldots, a_{q-1}],\\ \label{5.13}
\tau_k &\!\!\!=\!\!\!& \pi_\max [\ldots, a_{k-1}, a_k], \qquad {\rm
for} \, k>q-1.
\eea
 We have constructed a set of
representations of the category $\rm{CK} \bf{K}$ and among the subordinate ones, we
have encountered all the representations of all the semigroups
$\End(V_k(\bj))$. Therefore Proposition \ref{p5.1} immediately provides us with
an explicit form of the lowering functor:
\bea \label{5.14}
F_{n-1}^n \pi_0[\ldots, a_{n-1},a_n]&\!\!\!=\!\!\!&0,\nonumber \\
F_{n-1}^n\pi_\max [\ldots, a_{n-1}, 0]&\!\!\!=\!\!\!& \pi_\max [\ldots, a_{n-2},
a_{n-1}],
\eea 
while if $a_n\neq0$ then
\be \label{5.15}
F_{n-1}^n\pi_\max[\ldots, a_{n-1}, a_n] =\pi_0[\ldots, a_{n-2},
a_{n-1}].
\ee
\noindent
{\it Compatible systems}.  Let $\rm{CK} \bf{K}$  be an ordered
Cayley-Klein category. Suppose
that for each object $V_\sigma(\bj)$ we are given an irreducible  projective
representation $\rho_\sigma$ of the semigroup $\End(V_\sigma(\bj))$. We
call this family of representations a compatible system if
$F_{\sigma\rho\tau}^\tau=\rho_\sigma$ for every $\sigma>\tau$.

\begin{proposition}  \label{p5.2}
For any compatible system $\rho_\sigma$ there
exists a unique projective representation $(R, \rho)$ of the Cayley-Klein category
$\rm{CK} \bf{K}$ such that $\rho(P)=\rho_\sigma(P)$ for all $\sigma$
and all $P \in \End(V_\sigma(\bj))$.
\end{proposition}

Note that here the uniqueness is more important than the existence. It
is only the uniqueness that is used in our discussion of the
completeness of the lists of irreducible representations for the
Cayley-Klein categories $\rm{CK} \bf{A}, \rm{CK} \bf{B}, \rm{CK} \bf{C}, \rm{CK} \bf{D}$.

\begin{lemma} \label{l5.12} 
The semigroups $\End(V_\sigma(\bj))$ and the elements
$\lambda_{\sigma\tau}$ and  $\mu_{\tau\sigma}$ generate the entire
groupoid of morphisms of the Cayley-Klein category $\rm{CK} \bf{K}$.
\end{lemma}

\noindent
{\bf Proof} Let $P\in \Mor(V_\sigma(\bj),V_\tau(\bj))$. Let $\kappa>\tau$ and
$\kappa>\sigma$.  Then $P=\mu_{\kappa\tau} P'\lambda_{\sigma\kappa}$,
where $P'=\lambda_{\tau\kappa} P\mu_{\kappa\sigma}\in\End(V_\kappa(\bj))$.
\hfill{$\Box$}


\noindent
{\bf Proof of the proposition}

(a) {\it Uniqueness}. Denote by $R(V_\sigma(\bj))$ the space of the
representation $\rho_\sigma$.  Let $\sigma<\tau$. By (\ref{5.9}), the
operator $\rho(\mu_{\tau\sigma})$ vanishes on the kernel of the
projection $\rho(\vartheta_\tau^0)$.  Furthermore, the operators
\beq
\rho(\lambda_{\sigma\tau}): R(V_\sigma(\bj)) \rightarrow \im~~
\rho_\tau(\vartheta_\tau^\sigma), \qquad
\rho(\mu_{\tau\sigma}):\im~~ \rho_\tau(\vartheta_\tau^0) \rightarrow
R(V_\sigma(\bj))
\eeq
must be $\End(V_\sigma(\bj))$-intertwining operators which are inverses of
each other. Taking the irreducibility of $\rho_\sigma$ into
account, we find that $\lambda_{\sigma\tau}$ and $\mu_{\tau\sigma}$ are
uniquely defined to within multiplication by a constant. The required
uniqueness now follows from the lemma.

(b) {\it Existence}. Suppose then that we are given a compatible
system.  We define the operators
\beq
\rho(\lambda_{\sigma\tau}): R(V_\sigma(\bj))\rightarrow \im~~
\rho_\tau(\vartheta_\tau^\sigma), \qquad
\rho(\mu_{\tau\sigma}):\im~~ \rho_\tau(\vartheta_\tau^\sigma) \rightarrow
R(V_\sigma(\bj))
\eeq
so that they are $\End(V_\sigma(\bj))$-intertwining. 

\noindent
Suppose further that
$P\in \Mor(V_\alpha(\bj), V_\beta(\bj))$.  We choose $\kappa$ such that
$\kappa>\alpha$, $\kappa>\beta$ and we define $\rho(P)$ by the formula
\be \label{5.16}
\rho(P) =\rho(\mu_{\kappa\beta})
\rho_\kappa(\lambda_{\beta\kappa} P \mu_{\kappa\alpha})
\rho(\lambda_{\alpha\kappa}).
\ee
\begin{lemma}  \label{l5.13}
The operator $\rho(P)$ does not depend on the choice
of $\kappa$.
\end{lemma}
\noindent
{\bf Proof} We denote the expression (\ref{5.16}) by $\rho^\kappa(P)$. Let
$\xi>\kappa$. Then
\beq
\rho^\xi(P) = \rho(\mu_{\xi\beta}) \rho_{\xi}(\lambda_{\beta\xi}
P\mu_{\xi\alpha}) \rho(\lambda_{\alpha\xi}).
\eeq
It is clear from the definition of the operators
$\rho(\mu_{\xi\kappa})$ and $\rho(\lambda_{\kappa\xi})$ that for any
$S \in \End(V_\kappa(\bj))$ we have
\beq
\rho_\kappa(S) =s\cdot \rho(\mu_{\xi\kappa})
\rho\xi(\lambda_{\kappa\xi} S\mu_{\xi\kappa}) \rho(\lambda_{\kappa\xi})
\eeq
for some $s\in\Cb^*$. Using the obvious equalities
\beq
\rho(\mu_{\xi\beta}) =s'\cdot \rho(\mu_{\kappa\beta})
\rho(\mu_{\xi\kappa}), \qquad
\rho(\lambda_{\alpha\xi})=s'' \rho(\lambda_{\kappa\xi})
\rho(\lambda_{\alpha\kappa}),
\eeq
where $s', s'' \in \Cb$, we see that for some $t,t',\ldots\in \Cb^*$ we
have
\beq
\rho^\xi(P)  &\!\!\!=\!\!\!&
t\cdot \rho(\mu_{\xi\beta}) \rho_\xi(\lambda_{\beta\xi} P
\mu_{\xi\alpha}) \rho(\lambda_{\alpha\xi}) =\\
 &\!\!\!=\!\!\!&
t'\cdot \rho(\mu_{\kappa\beta}) \rho(\mu_{\xi\kappa})
\rho_\xi(\lambda_{\kappa\xi} \lambda_{\beta\kappa} P
\mu_{\kappa\alpha}\mu_{\xi\kappa}) \rho(\lambda_{\kappa\xi})
\rho(\lambda_{\alpha\kappa})=\\ &\!\!\!=\!\!\!& t'' \cdot
\rho(\mu_{\kappa\beta}) \rho_\kappa(\lambda_{\beta\kappa}
P\mu_{\kappa\alpha}) \rho(\lambda_{\alpha\kappa}) = t'''\cdot
\rho^\kappa(P).
\eeq
Suppose now that $\kappa$ and $\kappa'$ majorize $\alpha$ and $\beta$.
Then there exists $\xi$ majorizing $\kappa$ and $\kappa'$, and therefore
\beq
\rho^\kappa(P) =\rho^\xi(P)=\rho^{\kappa'}(P).
\eeq
It remains to verify the equality
\beq
\rho(P) \rho(Q)=s\cdot \rho(PQ)
\eeq
for any $Q \in \Mor(V_\sigma(\bj), V_\tau(\bj))$, $P\in 
\Mor(V_\tau(\bj),V_\gamma(\bj))$.
Suppose that $\xi$ is greater than $\sigma,\tau,\gamma$. Then for some
$s,s'\ldots \in\Cb^*$ we have
\beq
\rho(P)\rho(Q)  &\!\!\!=\!\!\!&
\rho(\mu_{\xi\gamma}) \rho_\xi(\lambda_{\gamma\xi} P\mu_{\xi\tau})
\rho(\lambda_{\tau\xi}) \rho(\mu_{\xi\tau})\rho_\xi(\lambda_{\tau\xi}
Q\mu_{\xi\sigma}) \rho(\lambda_{\sigma\xi})\\
 &\!\!\!=\!\!\!&
s\cdot \rho(\mu_{\xi\gamma}) \rho_\xi(\lambda_{\gamma\xi}
P\mu_{\xi\tau}) \rho_\xi(\vartheta_\xi^\tau) \rho_\xi(\lambda_{\tau\xi}
Q\mu_{\xi\sigma}) \rho(\lambda_{\sigma\xi})\\
 &\!\!\!=\!\!\!&
s'\cdot \rho(\mu_{\xi\gamma}) \rho_\xi (\lambda_{\gamma\xi}
P[\mu_{\xi\tau}\vartheta_\xi^\tau \lambda_{\tau\xi}] Q\mu_{\xi\sigma})
\rho(\lambda_{\sigma\xi})\\
 &\!\!\!=\!\!\!&
s''\cdot \rho(\mu_{\xi\gamma}) \rho_\xi (\lambda_{\gamma\xi} PQ
\mu_{\xi\sigma}) \rho(\lambda_{\sigma\xi}) =s'''\cdot \rho(PQ).
\eeq 

\vspace{-5mm}
\hfill{$\Box$}

\noindent
{\it Completeness of the lists of irreducible representations for
the Cayley-Klein categories $\rm{CK} \bf{A}, \rm{CK} \bf{B}, \rm{CK} \bf{C}, \rm{CK} \bf{D}$}. It follows immediately from (\ref{5.14}) and
 (\ref{5.15}) that any compatible system for $\rm{CK} \bf{K}=\rm{CK} \bf{A}, \rm{CK} \bf{B}, \rm{CK} \bf{C}, \rm{CK} \bf{D}$ has the
following form. There exists $q\in \Zb_+$ and a finite set [$\ldots,
a_{q-1}, a_q]\in\Zb_+^q$  such that $a_q\neq0$. The
representations $\rho_j$ of the semigroup $\End(V_j(\bj))$ must have the
form
\beq
\rho_j &\!\!\!=\!\!\!&  0 \quad {\rm for} \, j<q-1\\
\rho_{q-1}  &\!\!\!=\!\!\!& \pi_0[\ldots, a_{q-1}],\\
\rho_j &\!\!\!=\!\!\!& \pi_\max [\ldots, a_{q-1}, a_q,0,\ldots, 0]
\quad {\rm for} \, j>q-1.
\eeq
The representations corresponding to these compatible systems have
been constructed above. They have the numerical labels $(\ldots, a_{q-1},
q_q,0,0,\ldots)$. Thus assertion (b) of the
classification theorem is now proved.

\vspace{0.3cm}
\noindent
{\it Complete reducibility}
\begin{theorem} \label{t5.4}
Let $\rm{CK} \bf{K}$ be an ordered Cayley-Klein category. Suppose that the
finite-dimensional projective representations of all the semigroups
$\End_K(V_\sigma(\bj))$ are completely reducible. Then the projective
representations of $\rm{CK} \bf{K}$ are completely reducible.
\end{theorem}
\begin{lemma} \label{l5.14}
Suppose that the conditions of the theorem hold.
Let $(R,\rho)$ be a representation of $\rm{CK} \bf{K}$, and let $H\subset
R(V_\alpha(\bj))$ be an irreducible subrepresentation of a subordinate
representation of the semigroup $\End(V_\alpha(\bj))$. Then the $K$-cyclic
span $S=(S,\sigma)$ of the subspace $H$ is an irreducible
subrepresentation in $R$.
\end{lemma}
\noindent
{\bf Proof of the lemma} For each object $V_\tau(\bj)$ in the category 
$\rm{CK} \bf{K}$
we consider the set $M(V_\tau(\bj))$ of all $h\in S(V_\tau(\bj))$ such that the
cyclic span of $h$ is distinct from $S$.  The condition $h \in
M(V_\tau(\bj))$ is equivalent to the requirement that $\rho(P)= 0$ for all
$P \in \Mor(V_\tau(\bj), V_\alpha(\bj))$.  (If $\rho(P)h\neq0$, then the cyclic
span of $h$ contains $H$ and hence is the whole of $S$.) Suppose that
for some $\tau>\alpha$ the space $M(V_\tau(\bj))\neq 0$.  Let $K(V_\tau(\bj))$ be
a complement of $M(V_\tau(\bj))$ that is invariant with respect to the
semigroup $\End(V_\tau(\bj))$. If $K(V_\tau(\bj))= 0$, then the factor
representation $S/M$ satisfies the condition $(S/M)(V_\tau(\bj))= 0$ and we
obtain a contradiction.  Hence $K(V_\tau(\bj))\neq 0$.  The cyclic span of
any vector $h\in K(V_\tau(\bj))$ under the action of $K$ contains $H$ (since
$h\notin M(V_\tau(\bj))$) and is therefore the whole of $S$. On the other
hand, that part of the cyclic span which lies in $S(V_\tau(\bj))$ coincides
with $\End(V_\tau(\bj))h$, which is the cyclic span of the vector $h$; and
the latter, in turn, is contained in $K(V_\tau(\bj))\neq S(V_\tau(\bj))$.
This is a contradiction. Thus, $M(V_\tau(\bj)) = 0$ for all $\tau>\alpha$.
But it now follows that $M(V_\tau(\bj))= 0$ for all $\tau$. Thus $S$ is
irreducible.  \hfill{$\Box$}

\noindent
{\bf Proof of the theorem} We choose some irreducible representation
$S$ of the representation $R =(R,\rho)$. Let $R(V_\alpha(\bj))\neq
S(V_\alpha(\bj))$.  Let $T(V_\alpha(\bj))$ be an $\End(V_\alpha(\bj))$-invariant
complement of $S(V_\alpha(\bj))$.  We take the $K$-cyclic span $T$ of the
subspace $T(V_\alpha(\bj))$. Clearly for any $V_\beta(\bj)$ we have
$T(V_\beta(\bj))\cap S(V_\beta(\bj))=0$ (otherwise, the set of subspaces
$T(V_\beta(\bj))\cap S(V_\beta(\bj))$  would form a non-trivial subrepresentation
of $S$, which, as we recall, is irreducible). It may turn out,
however, that for some $V_\beta(\bj)$ the equality $S(V_\beta(\bj))\oplus
T(V_\beta(\bj)) =R(V_\beta(\bj))$ does not hold.  There then exists
$\gamma>\alpha$ such that $S(V_\gamma(\bj))\oplus T(V_\gamma(\bj)) \neq
R(V_\gamma(\bj))$ (for this, it suffices to take any $\gamma$ such that
$\gamma>\alpha$ and $\gamma>\beta$, for
$R(V_\gamma(\bj))/[S(V_\gamma(\bj))\oplus T(V_\gamma(\bj))]=0$ would then imply that
$R(V_\beta(\bj))/[S(V_\beta(\bj))\oplus T(V_\beta(\bj))]=0)$. We take a complement
$T'(V\gamma(\bj))$ of $S(V_\gamma(\bj))$  in $R(V_\gamma(\bj))$ that is invariant
with respect to $\End(V_\gamma(\bj))$ and is such that $T'(V_\gamma(\bj))\supset
T(V_\gamma(\bj))$.  Let $T'$ be the cyclic span of $T(V_\gamma(\bj)$).  Clearly,
for all $\mu$ we have $T'(V_\mu(\bj))\supset T(V_\mu(\bj))$ (in fact,
$T'(V_\alpha(\bj))\supset T(V_\alpha(\bj))$, because the cyclic span $N$ of
$T(V_\gamma(\bj))$ contains $T(V_\alpha(\bj))$; indeed, $(T/N)(V_\gamma(\bj)) = 0$,
and hence, $(T/N)(V_\alpha(\bj)) = 0$).

We can now continue this procedure and choose $T''\supset T'$,
$T'''\supset T'',$ ... and so on, and then take their union. If the
partially ordered set $\Sigma$ of indices numbering the objects of the
Cayley-Klein category $\rm{CK} \bf{K}$ has a fairly complicated construction,
 then one has to invoke the standard procedures
relating to Zorn's lemma.   This completes the proof of the theorem.

The theorem on the classification  of holomorphic projective
representations of the Cayley-Klein categories $\rm{CK} \bf{A}, \rm{CK} \bf{B}, \rm{CK} \bf{C}$, and $\rm{CK} \bf{D}$ is now proved in
its entirety.


\end{document}